\documentclass[smallextended]{svjour3}
\usepackage[utf8]{inputenc}
\usepackage{amsmath}
\usepackage{amsfonts}
\usepackage{amssymb}
\usepackage{graphicx}
\usepackage{subfigure}
\usepackage{epsfig}
\usepackage[english]{babel}
\usepackage[left=2cm,right=2cm,top=2cm,bottom=2cm]{geometry}
\usepackage {tikz}
\usetikzlibrary {patterns}
\usepackage{multicol}
\usepackage{color}
\usepackage{ulem}

\DeclareMathOperator{\sign}{sign}

\begin{document}

\author{Roc\'io Balderrama \and Javier Peressutti \and \\ Juan Pablo Pinasco \and Constanza S\'anchez de la Vega \and \\ Federico Vazquez}

\institute{Roc\'io Balderrama \at Departamento de Matem\'atica, Facultad de
Ciencias Exactas y Naturales, Universidad de Buenos Aires \\ Ciudad
Universitaria, Pabell\'on I, C1428EGA Buenos Aires, Argentina.
\and Javier Peressutti \at Instituto de F\'isica de Mar del Plata (IFIMAR) CONICET, UNMDP, Departamento de F\'isica, Universidad Nacional de Mar del Plata, Funes 3350,
7600 Mar del Plata, Argentina.
\and Juan Pablo Pinasco and Constanza S\'anchez de la Vega \at IMAS - CONICET and Departamento de Matem\'atica, Facultad de
Ciencias Exactas y Naturales, Universidad de Buenos Aires, Ciudad
Universitaria, Pabell\'on I, C1428EGA Buenos Aires, Argentina.
\and Federico Vazquez \at Instituto de C\'alculo, FCEN, Universidad de Buenos Aires and CONICET, C1428EGA Buenos Aires, Argentina.}


\title{\bf  Optimal control for a SIR epidemic model with limited quarantine}

\maketitle

\nocite{*}

\begin{abstract}
We study first order necessary conditions for an optimal control problem of a
Susceptible-Infected-Recovered (SIR) model with limitations on the duration of the
quarantine. The control is done by means of the reproduction number, i.e., the number of
secondary infections produced by a primary infection, which represents an external
intervention that we assume time-dependent.  Moreover, the control function can only
be applied over a finite time interval, and the duration of the most strict quarantine
(smallest possible reproduction number) is also bounded.  We consider a maximization
problem where the cost functional has two terms: one is the number of susceptible
individuals in the long-term and the other depends on the cost of interventions.  When the
intervention term is linear with respect to the control, we obtain that the optimal solution is
bang-bang, and we characterize the times to begin and end the strict quarantine.  In the
general case, when the cost functional includes the term that measures the intervention
cost, we analyze the optimality of controls through numerical computations.
\end{abstract}

\section{Introduction}

The Covid-19 pandemic outbreak raises an unprecedented series of decisions in different
countries around the world. Since vaccines and effective pharmaceutical treatments were not
initially available, governments had decided to impose non-pharmaceutical interventions like
social distance, quarantines and total lock-downs as the most effective tools to mitigate the
spread of the disease. Although these kinds of measures are helpful in reducing the virus
transmission and giving time to health systems to adapt, they could be extremely stressful in
terms of economic and social costs, and, in longer periods, tend to have less compliance with
the population.

In this article we consider the classical SIR model introduced by
Kermack  and McKendrick  in \cite{kermack1927contribution} and widely used in epidemiology
\cite{anderson1992infectious,brauer2012mathematical}, where the population is divided in compartments of
Susceptible, Infected and Recovered (or Removed) individuals.  As it is usual in SIR models, we assume that people who have recovered develop immunity
and, therefore, would not be able to get infected nor infect others.  We consider that infection and recovery rates are allowed to change over time, and that are homogeneous among the population, instead of heterogeneous rates as in \cite{buono2013slow,ferrari2021coupling}.
We also assume a mean-field hypothesis that implies random interactions between any pair of agents, unlike the works in \cite{ferreyra2021sir,janson2014,lagorio2011quarantine,velasquez2017interacting,volz2008sir} where interactions are mediated by an underlying network of contacts.

Optimal control problems for a system governed by a SIR or a SEIR model (with the addition of the Exposed compartment) with pharmaceutical interventions as vaccination or treatment were widely studied, see for instance  \cite{dePinho2014} and \cite{Ledzewicz2011}.
On its part, control problems with non-pharmaceutical interventions were less studied. Non-pharmaceutical interventions can range from a mild mitigation policy to a strong
suppression policy. As discussed by Ferguson et al. in \cite{ferguson2020report}, a
suppression policy ``{\it aims to reverse epidemic growth, reducing case numbers to low
levels and maintaining that situation indefinitely}''. Suppression can be achieved by
restricting travels, closing schools and nonessential businesses, banning social gatherings,
and asking citizens to shelter in place. These measures, often referred to as a lock-down, are
highly restrictive on social rights and damaging to the economy. In contrast, a mitigation
policy ``{\it focuses on slowing but not necessarily stopping epidemic spread}''. Mitigation
measures may involve discouraging air travel while encouraging remote working, requiring
companies to provide physical separation between workers, banning large gatherings,
isolating the vulnerable, and identifying and quarantining contagious individuals and their
recent contacts.

 A critical parameter in the SIR model is the basic reproduction number $R_0$,
defined as the quotient between the rates of contagion and recovery (see for instance
\cite{anderson1992infectious}) . At the beginning of the epidemic, when no one in the
population is immune,  infected individuals will infect $R_0$ other people on average. Let us
observe that, for  $R_0<1$, the number of new cases decline, and when $R_0>1$, the
number of new cases grows. However, at any time $t>0$, the effective reproduction number
$R_t$ replaces $R_0$, since the number of contacts between infected and susceptible
agents is reduced due to the interactions with recovered individuals that are immune. Hence,
the epidemic grows
 until a sufficient fraction of the population becomes infected, and after reaching a peak starts to gradually decline.
Following \cite{ferguson2020report}, the suppression phase can achieve $R_0<1$, while the
mitigation measures are unlikely to bring $R_0$ below $1$. Therefore, the number of new
cases are expected to decline during the suppression phase and to start rising again during
the mitigation phase, although at a slower rate than in a non-intervention scenario.

\bigskip

 We assume that  the intervention lasts a preset period of time $T$ (see
\cite{greenhalgh1988some,ketcheson2020optimal}), since it is unrealistic that interventions
can be sustained indefinitely. Also, the lock-down can last at most a given time $\tau<T$,
the maximum time that the population will adhere. Now, there are several interesting
questions related to the implementation of the measures:

\begin{enumerate}
\item When should the suppression policy begin  in $[0,T]$?

\item Is it convenient to split the maximum time $\tau$ into different intervals?

\item Is it better to apply  a strong lock-down followed by mild mitigation measures or not?
\end{enumerate}

In this article we study the previous questions using optimal control tools and numerical
computations (see for instance
\cite{Clarke2013,zbMATH00784953,zbMATH01579035,zbMATH05164170}). The answers clearly can
depend on the goal, which in our case is to maximize the number of people who remains
susceptible in the long-term, that we will call $x_\infty(x(T),y(T),\sigma_0)$ ( see  Section
\ref{se: bang-bang} for the precise details on notation), plus a term that measures the social
and economic cost of interventions. To account for quarantine measures, we consider a
time-dependent contact rate. Using an optimal control approach we derive the necessary
optimality conditions, and we show that the  optimal strategy is of bang-bang type.
Moreover, we characterize the time to start and finish the lock-down during the intervention
phase. Let us remark that these questions make sense also in SIHR models, which include
hospitalized individuals, since it must be necessary to keep the maximum of the hospitalized
group below some threshold (see for instance \cite{kohler2020robust}).

Recently, many works have appeared dealing with these and related issues.
 The timing for start the
suppression measures trying to maximize $x_\infty$ was studied in
\cite{ketcheson2020optimal}, assuming that the control has no extra costs, and the existence
of a bang-bang type control was proved. However, in that work it is assumed that the
lock-down corresponds to a zero reproduction number, something that is impossible to achieve in the real world.
Moreover, it is assumed that the strict lock-down can last during the whole intervention,
which seems to be impracticable. As we will see below, suppression measures can  be
applied at different points of the interval $[0,T]$ depending on the initial proportion of
infected individuals.  On the other hand, in \cite{morris2020optimal} this problem was
analyzed for a different objective function, i.e., minimizing the peak of infected individuals,
and let us cite also \cite{godara2021control}, where the authors do not consider a limit time
for interventions.

 The second question is  suggested by the strategy proposed  in \cite{ferguson2020report}:
    the lock-down must be turned {\it on} and {\it off} several times based on the incidence of the virus in the population.
A control-theoretic approach  was considered in several works (see for instance
    \cite{kohler2020robust,palmer2020optimal,tsay2020modeling}), although no time limits for the interventions
were imposed.  We shall prove that the optimal policy is of bang-bang type, that
is, we can turn {\it on} the quarantine only once, and it must be turned {\it off}
after the maximum allowed time $\tau$.

Finally, the third question involves both suppression and mitigation phases,
and one of the policies was colorfully characterized as  {\it the hammer and the dance} in
\cite{pueyo2020coronavirus}: a strict lock-down, followed by mitigation measures in order to
keep under control the propagation of the disease. Our main result, theorem \ref{te:
control_optimo_sigma1gral},  characterizes the optimal control in a fixed time interval
$[0,T]$ as a bang-bang control, with suppression measures in some interval of maximal
length $\tau$, and mitigation measures applied at the start or at  the end (or both) of the
pandemic. Hence,  it makes sense to apply mild mitigation measures, followed by a strong
lockdown, and then mild mitigation measures again when the hypothesis of item 2 of
the theorem are
satisfied. 
In particular, this hypothesis is satisfied when the lockdown is very strict and the initial fraction of susceptible individuals is larger than $1/R_0$.
There are different ways to justify this policy given in item 2 from theorem 3 using theorem 2: 
the first one is to consider a mitigation policy close to no intervention, and the continuity of solutions with respect to coefficients of the
system of differential equations implies that the solutions obtained from theorem \ref{te:
control_optimo_sigma1gral} are also close, hence if we are under hypothesis given by condition 2, the time to start the strong lockdown could
not be the initial time. 
Another one is to consider a {\it New Normal}, where mild mitigation
measures will remain forever: in this case we can apply directly theorem \ref{te:
control_optimo_sigma1gral} and again we get solutions where the strong lockdown starts
after some interval of time. Finally, numerical simulations in Section \ref{se: simulations} support all these results.

\bigskip
The paper is organized as follows. We set the problem and state the description of the model
in Section \ref{se: model}. In Section \ref{se: bang-bang} we apply Pontryagin's maximum
principle  to our problem and prove that the control is bang-bang in Lemma \ref{le: optimo es
bb}. Also, in Lemma \ref{te: dos_saltos} we show that it jumps at most twice. The main
results that characterize the optimal control are theorem \ref{te: control_optimo_sigma1gral} and
Corollary \ref{co: sigma1es0}, presented in Section \ref{se: caracterization_optimal}, and theorem \ref{th: version_general} in Section
 \ref{se: kappa_posit}.  These two theorems and the corollary are tested in Section \ref{se: simulations} by means of numerical simulations. Finally, in Section \ref{se: conclusions} we discuss our results
including future work.

\section{Description of the model}
\label{se: model}

The basic compartmental model of infectious diseases introduced by Kermack and McKendrik
\cite{kermack1927contribution} divides the population in three compartments with
homogeneous characteristics: Susceptible, Infected, and Recovered or Removed (SIR). In this
model, births and deaths are neglected and the recovered population is assumed to no longer
infect others and also cannot be reinfected. By re-scaling time, we can  get a recovery rate
equal to one, and we call $\sigma_0$ the basic reproduction number
under no interventions.

We address epidemics with no vaccination, where the only possible control is isolation. We
model this non-pharmaceutical intervention via a time dependent reproduction number
$\sigma(t)\in [\sigma_1,\sigma_2]$, where $\sigma_1$ corresponds to a more strict
isolation (hard quarantine) than $\sigma_2$ (soft quarantine), and assume that this intervention can only be applied over a
finite time interval $[0,T]$. After the intervention,
  the restrictions are removed, thus the disease spreads freely and $\sigma(t)=\sigma_0$ for all $t>T$.

Then, we denote by $(x,y)$ the proportion of susceptible and infected individuals, and their
evolution is governed by the following system of coupled nonlinear ordinary differential
equations
\begin{subequations} \label{eq: SIR}
\begin{align}
\label{eq: SIR_x}
x' &= -\gamma \sigma(t) x y, \\
 \label{eq: SIR_iy}
 y' &= \gamma \sigma(t) x y - \gamma y,
\end{align}
\end{subequations}
with $(x(0),y(0))\in {\cal{D}}=\left\{ (x_0,y_0): x_0>0, y_0>0, x_0+y_0\le 1\right\}$ and
$\sigma(t) \in [\sigma_1,\sigma_2]$ for a.e. $t\in [0,T]$, where
$0\le \sigma_1<\sigma_2\le\sigma_0$. 

We also assume that during this time interval $[0,T]$, it is not possible to impose an
extremely restrictive isolation during a long time. Thus, we consider that the more restricted
quarantine (correponding to $\sigma_1$) can lasts at most for a time
pertiod $\tau$, with $\tau\in
(0,T)$. We can describe this restriction with an isoperimetric inequality:
\begin{align} \label{eq: desig_iso_sigma}
\int_0^T \sigma(t) dt \ge \sigma_1 \tau + \sigma_2 (T-\tau).
\end{align}

Note that if $\sigma$ is a control that takes the value $\sigma_1$ for a period of time
$\tilde{\tau}>\tau$ we would have that $\int_0^T \sigma \le \sigma_1  \tilde{\tau} +
\sigma_2(T-\tilde{\tau})<\sigma_1 \tau+ \sigma_2 (T-\tau)$ contradicting inequality \eqref{eq: desig_iso_sigma}.

Our goal is to
maximize the long time limit of the susceptible fraction $x_{\infty}=\lim_{t\to \infty} x(t)$.
Since $x_{\infty}=1/\sigma_0$ only at a single point where $y=0$ and since $y=0$ cannot be
reached at finite time (for instance at time $T$) for $y_0>0$, then any solution satisfies
$x_{\infty}<1/\sigma_0$.

Once the period of intervention is finished at time $T$ we compute $x_{\infty}(x(T),y(T),\sigma_0)=\lim_{t\to \infty} x(t)$
where $(x,y)$ is the solution of system \eqref{eq: SIR} with initial data $(x(T),y(T))$ and constant reproduction number
$\sigma(t)\equiv \sigma_0$ for $t>T$.

Given $(x_0,y_0)\in {\cal D}$, $T,\tau$ fixed satisfying $0<\tau <T$, $0\le \sigma_1<\sigma_2\le \sigma_0$, 
we then consider the following optimal control problem with an objective function $J$ that includes in addition a term that accounts
for the running economic cost of the control:

\begin{quote}
\begin{subequations}
\label{optimal_problem_T_isop}
\begin{align}
\label{eq: funcional_J}
\max \quad &J(x,y,\sigma):=x_{\infty}(x(T),y(T),\sigma_0) + \int_0^T L(x(t),y(t),\sigma(t)) dt \\
\text{s.t.} \quad &x' = -\gamma \sigma(t) x y,  \quad x(0) = x_0, \quad t\in [0,T], \label{eq: SIR_x_initial}\\
&y' = \gamma \sigma(t) x y - \gamma y, \quad y(0)= y_0, \quad t\in [0,T], \label{eq: SIR_y_initial} \\
&\int_0^T \sigma(t) dt \ge \sigma_1 \tau + \sigma_2 (T-\tau)  \label{eq: int_restriction}\\
&\sigma_1 \le \sigma(t) \le \sigma_2.
\label{optimal_problem_rest_sigma_isop}
\end{align}
\end{subequations}
\end{quote}

We add a new state variable given by $v(t)=\int_0^t \sigma(s) ds$ and consider $\sigma:[0,T] \to [\sigma_1,\sigma_2]$ in the class
of Lebesgue-measurable functions in order to prove existence of optimal solution. Thus, we can study the equivalent optimal control problem:
\begin{quote}
\begin{subequations} \label{eq: SIR_iso}
\begin{align}
\label{eq: funcional_J_iso}
\max \quad& J(x,y,v,\sigma):=x_{\infty}(x(T),y(T),\sigma_0) + \int_0^T L(x(t),y(t),\sigma(t)) dt \\
\label{eq: SIR_iso_x}
s.t. \quad &x'(t)= -\gamma \sigma(t) x(t) y(t),   \quad x(0) = x_0, \quad t\in [0,T],\\
 \label{eq: SIR_iso_y}
 &y'(t) = \gamma \sigma(t) x(t) y(t) - \gamma y(t),  \quad y(0)= y_0, \quad t\in [0,T],\\
 \label{eq: SIR_iso_v}
 &v'(t)  =\sigma(t),  \quad v(0)= 0,  \quad t\in [0,T], \\
 \label{eq: sigma_S_T}
 &\sigma(t) \in [\sigma_1,\sigma_2], \quad \text{a.e. } t\in [0,T] \\
  \label{eq: restr_vT} 
  &v(T)  \ge \sigma_2 T +(\sigma_1-\sigma_2) \tau. 
\end{align}
\end{subequations}
\end{quote}
We will refer to a $4$-tuple $(x,y,v,\sigma)$ as an admissible process of the underlying control system if the control $\sigma$ is a measurable function and the state $(x,y,v)$ is an absolutely continuous vector function satisfying \eqref{eq: SIR_iso_x}-\eqref{eq: restr_vT}. 
The optimal control problem consists in finding an optimal admissible process $(x^*,y^*,v^*,\sigma^*)$ that maximizes the cost $J$. In this case, we refer to the control $\sigma^*$ as optimal control.

Next, we give a result on existence of solution for the optimal control problem \eqref{eq: SIR_iso}.
\begin{proposition} \label{pr: existencia}
Assume $L$ is continuous with respect to all its variables, $L(x,y,\cdot)$ is convex for each $(x,y)$ and there exists a constant $\alpha_0$
such that for all $(x,y)$ satisfying $0\le x+y\le 1$ and $\sigma\in [\sigma_1,\sigma_2]$, it holds $L(x,y,\sigma)\ge \alpha_0$. Then, the optimal control problem \eqref{eq: SIR_iso} admits a solution.
\end{proposition}

\begin{proof}
The proof follows directly from theorem 23.11 given in \cite{Clarke2013} which assures the existence of an optimal process.  Since the control space is $[\sigma_1,\sigma_2]$, solutions of system \eqref{eq: SIR_iso_x}-\eqref{eq: SIR_iso_y} satisfy that $0\le x+y \le 1$, the application $x_{\infty}(x,y)$ is continuous and hypothesis on $L$,  it is straightforward to prove conditions a) to f) from theorem 23.11.  Moreover, taking $\sigma(t)\equiv \sigma_2$, we see that the unique solution of system  \eqref{eq: SIR_iso_x}-\eqref{eq: SIR_iso_v} together with $\sigma$ gives an admissible process for which $J$ is finite completing the hypothesis of theorem 23.11.
\end{proof}

In the next section we derive the first order necessary conditions for
the optimal control problem \eqref{optimal_problem_T_isop}. Then,
assuming that $L$ depends only linearly on the control we prove that
an optimal control must be bang--bang  (Lemma \ref{le: optimo es bb}) and
that has at most two jumps (theorem \ref{te: dos_saltos}). Moreover, we
give the main result of the article (theorem \ref{te:
  control_optimo_sigma1gral}) where we characterize these jumps 
 for the case of zero economic cost and $\sigma_2=\sigma_0$.

\section{The optimal control is bang-bang}
\label{se: bang-bang}

We can compute the partial derivatives of $x_{\infty}(x(t),y(t),\sigma_0)$ with respect to
$x(t)$ and $y(t)$ in the same way that it is done in \cite{ketcheson2020optimal}.
\begin{subequations} \label{eq: deriv_xinfty}
\begin{align} \label{eq:deriv_xinfty_x}
\frac{\partial x_{\infty}(x(t),y(t),\sigma_0)}{\partial x(t)}&=\frac{1-\sigma_0 x(t)}{x(t)} \frac{x_{\infty}(x(t),y(t),\sigma_0)}{1-\sigma_0 x_{\infty}(x(t),y(t),\sigma_0)}, \\
 \label{eq:deriv_xinfty_y}
\frac{\partial x_{\infty}(x(t),y(t),\sigma_0)}{\partial y(t)}&=-\frac{\sigma_0 x_{\infty}(x(t),y(t),\sigma_0)}{1-\sigma_0 x_{\infty}(x(t),y(t),\sigma_0)}.
\end{align}
\end{subequations}

We consider the Hamiltonian $H$
\begin{align} \label{Hamiltonian}
H(x,y,v,\sigma,\lambda)= \lambda_0 L(x,y,\sigma)-\lambda_1 (\gamma  \sigma x y) +\lambda_2 (\gamma \sigma x y -\gamma y ) +\lambda_3 \sigma
\end{align}
where $\lambda_0\ge 0$, $\lambda \in \mathbb{R}^3$. Given an optimal solution
$(x^*,y^*,v^*,\sigma^*)$, we denote $$L^*[t]=L(x^*(t),y^*(t),v^*(t),\sigma^*(t)).$$ 
Assuming that $L$ is continuous and admits derivatives $L_x$ and $L_y$ which are themselves continuous in all the variables, the
necessary conditions for a maximum process $(x^*,y^*,v^*,\sigma^*)$ on $[0,T]$ 
are the following (see \cite{Clarke2013,zbMATH00784953}):

\begin{quote}
There exists a real number $\lambda_0 \ge 0$, the adjoint variable $\lambda:[0,T] \to \mathbb{R}^3$ which is absolutely continuous, and $\beta\in \mathbb{R}$ such that $(\lambda_0,\lambda(t),\beta)\neq 0$ for every $t$ and the following conditions hold:
\begin{enumerate}
\item The adjoint variables $\lambda_i(t)$ satisfy  a.e. $t\in [0,T]$
\begin{subequations} \label{eq: lambda_expl}
\begin{align}
\label{eq: lambda_expl1}
\lambda_1'(t)& =-\lambda_0  L^*_x[t]+ (\lambda_1(t) -\lambda_2(t)) \gamma \sigma(t) y(t),\\
\label{eq: lambda_expl2} 
\lambda_2'(t)& = -\lambda_0 L^*_y[t] +(\lambda_1(t)  -\lambda_2(t)) \gamma \sigma(t)x(t)+\gamma \lambda_2(t),\\
\label{eq: lambda_expl3}
\lambda_3'(t)& =-\lambda_0 L^*_v[t],  
\end{align}
\end{subequations}
with final time conditions (using the abbreviation $x_{\infty}$ for $x_{\infty}(x(T),y(T),\sigma_0)$)
\begin{subequations}  \label{eq: lambda_final}
\begin{align}
\lambda_1(T)&=\lambda_0 \frac{\partial x_{\infty}}{\partial x(T)}=\lambda_0 \frac{1-\sigma_0 x(T)}{x(T)} \frac{x_{\infty}}{1-\sigma_0 x_{\infty}}, \label{eq: lambda_final1} \\
 \lambda_2(T)&=\lambda_0 \frac{\partial x_{\infty}}{\partial y(T)}=-\lambda_0\frac{\sigma_0 x_{\infty}}{1-\sigma_0 x_{\infty}}, \label{eq: lambda_final2}\\
 \lambda_3(T)&=\beta \ge 0 \text{  and  } \lambda_3(T) (v(T)-\sigma_2 (T-\tau) - \sigma_1 \tau)=0. \label{eq: lambda_final3}
 \end{align}
\end{subequations}
\item  For a.e. $t\in [0,T]$
\begin{align} \label{eq: optim_cond_sigma}
&\lambda_0 L^*[t] + \sigma^*(t) \left(\gamma x^*(t) y^*(t) (\lambda_2(t)-\lambda_1(t))+\lambda_3(t) \right) \\
&=\max_{\sigma: \sigma_1\leq \sigma \le \sigma_2} \lambda_0 L(x^*(t),y^*(t),\sigma) + \sigma  \left(\gamma x^*(t) y^*(t) (\lambda_2(t)-\lambda_1(t)) +\lambda_3(t) \right). \nonumber
\end{align}
\item There exists a constant $C$ such that for a.e. $t\in [0,T]$
\begin{align} \label{eq: H_es_cte}
\lambda_0 L^*[t] + \sigma^*(t) \left(\gamma x^*(t) y^*(t) (\lambda_2(t)-\lambda_1(t))+\lambda_3(t) \right) -\gamma \lambda_2(t) y^*(t) =C.
\end{align}
\end{enumerate}
\end{quote}

We have the following result:

\begin{lemma}
The optimal control problem is normal.
\end{lemma}
\begin{proof}
Assume $\lambda_0=0$. From \eqref{eq: lambda_expl1}-\eqref{eq: lambda_expl2} with final
time conditions $\lambda_1(T)=\lambda_2(T)=0$, $\lambda_1(t)=\lambda_2(t)=0$ for all
$t\in [0,T]$. Since the multipliers $(\lambda_0,\lambda(t),\beta)\neq 0$, then
$\lambda_3(t)\equiv \beta >0$. Therefore, from the optimality condition given in \eqref{eq:
optim_cond_sigma}, $\sigma^*(t)=\sigma_2$ a.e. $t\in [0,T]$ contradicting the
complementarity condition $v(T)=\sigma_2 (T-\tau)+\sigma_1 \tau$ given in \eqref{eq:
lambda_final3}. Thus, we can assume $\lambda_0=1$, and the proof is finished.
\end{proof}

In what follows we will restrict ourselves to the case where $L$ is a linear function only depending on the control. 
Let $L(x(t),y(t),\sigma(t))=\kappa \sigma(t)$ with $\kappa\ge 0$.
Then  $L$ satisfies the hypothesis of proposition \ref{pr: existencia} showing the existence of an optimal process. Also $L$ satisfies the regularity assumptions needed for the derivation of first order necessary conditions \eqref{eq: lambda_expl}-\eqref{eq: H_es_cte}. In this case the functional $J$ reads
\begin{align}
\label{eq: funcional_J_Llineal}
J(x(t),y(t),\sigma(t),T)&=x_{\infty}(x(T),y(T),\sigma_0) + \kappa \int_0^T\sigma(t) dt.
\end{align}
This functional penalizes the use of a more restrictive quarantine, corresponding to smaller values of $\sigma$.

\begin{lemma} \label{le: optimo es bb}
Let $L(x(t),y(t),\sigma(t))=\kappa \sigma(t)$ with $\kappa\ge 0$ and let $\sigma^*$ be an
optimal control. Then $\sigma^*(t)$ is a bang-bang control.
\end{lemma}

\begin{proof}
From the optimality conditions \eqref{eq: lambda_expl}-\eqref{eq: H_es_cte}, we obtain:
\begin{enumerate}
\item The adjoint variables $\lambda_i(t)$ satisfy
\begin{subequations}  \label{eq: lambda_final_b}
\begin{align}
\lambda_1'(t)& = (\lambda_1(t)-\lambda_2(t)) \gamma \sigma(t) y^*(t)  \text{ a.e. } t\in [0,T],  \label{eq: lambda1_b} \\
\lambda_1(T)&= \frac{1-\sigma_0 x(T)}{1-\sigma_0 x_{\infty}} \frac{x_{\infty}}{x(T)}.  \label{eq: lambda1T_b} \\
\lambda_2'(t)& =(\lambda_1(t) -\lambda_2(t) )\gamma \sigma(t) x^*(t)+ \gamma \lambda_2(t)  \text{ a.e. } t\in [0,T],  \label{eq: lambda2_b} \\
 \lambda_2(T)&=-\frac{\sigma_0 x_{\infty}}{1-\sigma_0 x_{\infty}}.   \label{eq: lambda2T_b}\\
\lambda_3'(t)& =0  \text{ a.e. } t\in [0,T] \rightarrow \lambda_3\equiv \lambda_3(T),\\
 \lambda_3(T)&=\beta \ge 0 \text{  and  } \lambda_3(T) (v(T)-\sigma_2 (T-\tau) - \sigma_1 \tau)=0.   \label{eq: lambda3T_b}
\end{align}
\end{subequations}
\item Define
\begin{align} \label{eq: phi}
\phi(t)=\kappa+\lambda_3(T)+\gamma x^*(t) y^*(t) (\lambda_2(t)-\lambda_1(t)),
\end{align}
then, from \eqref{eq: optim_cond_sigma},
\begin{align} \label{eq: optim_cond_sigma_b}
& \sigma^*(t) \phi(t) =\max_{\sigma: \sigma_1\leq \sigma \le \sigma_2} \sigma  \phi(t). \end{align}
\item For all $t\in [0,T]$
\begin{align} \label{eq: H_es_cte_b}
 \sigma^*(t) \phi(t)-\gamma \lambda_2(t) y^*(t) =C.
\end{align}
\end{enumerate}
Assume $\phi(t)=0$ on an interval $[a,b]\subset [0,T]$, then
$$\kappa+\lambda_3(T)=\gamma x^*(t) y^*(t) (\lambda_1(t)-\lambda_2(t)) \text{  for all  }  t\in [a,b].$$
Computing the derivative we obtain
$$0=-\gamma x^*(t) y^*(t)\lambda_1(t) \text{ for all } t\in (a,b).$$
Then $\lambda_1(t)=\lambda_2(t)=0$ for all $t\in (a,b)$ and thus for all $t\in [0,T]$, contradicting the end point conditions. Therefore, there cannot be singular arcs and the control $\sigma^*$ is given by:
\begin{equation} \label{eq: bang bang}
\sigma^*(t)= \left\{ \begin{array}{ll}
\sigma_2 & \text{ if } \phi(t)>0\\
\sigma_1& \text{ if } \phi(t)<0.
\end{array}\right.
\end{equation}
\end{proof}

\begin{lemma}
Let $(x_0,y_0)$ be given and $(x,y,v,\sigma)$ an admissible process. Then for $t\ge 0$
$$x_{\infty}(x(t),y(t),\sigma_0)\ge x_{\infty}(x_0,y_0,\sigma_0) ,$$
and therefore
\begin{align}
\label{signo_xinf}
x_{\infty}(x_0,y_0,\sigma_0) \le x_{\infty}(x(T),y(T),\sigma_0)< 1/ \sigma_0  .
\end{align}
\end{lemma}
\begin{proof}
  See \cite{ketcheson2020optimal}.
\end{proof}

In the next lemma we will see that the switching function changes sign at most two times,
concluding that an optimal control $\sigma^*$ jumps at most twice.

\begin{lemma} \label{le: dos saltos}
The switching function $\phi$ given in \eqref{eq: phi} changes sign at most twice.
\end{lemma}

\begin{proof}
The proof follows by analysing the phase diagram of $\lambda_1,\lambda_2$.
We begin by noting that a solution $(\lambda_1,\lambda_2)$ of the system \eqref{eq: lambda1_b}-\eqref{eq: lambda2T_b} cannot cross both semilines $\lambda_1=\lambda_2>0$ and $\lambda_1=\lambda_2<0$. This is a consequence of condition \eqref{eq: H_es_cte_b}. Assume there exist $s_1,s_2 \in [0,T]$ such that $\lambda_1(s_1)=\lambda_2(s_1)>0$ and $\lambda_1(s_2)=\lambda_2(s_2)<0$. Evaluating \eqref{eq: phi} on $t=s_i$ for $i=1,2$ we have that $\phi(s_i)=\kappa + \lambda_3(T)\ge 0$ and from \eqref{eq: H_es_cte_b},  $\lambda_2(s_i)=-\frac{C}{\gamma y^*(s_i)}$ if $\kappa+\lambda_3(T)=0$ or, using \eqref{eq: bang bang},
 $\lambda_2(s_i)=\frac{\sigma_2 (\kappa + \lambda_3(T))-C}{\gamma y^*(s_i)}$ if $\kappa+\lambda_3(T)>0$, both cases contradicting that $\lambda_2(s_1)$ and $\lambda_2(s_2)$ had opposite signs.

Since we have end time conditions on $T$ we go backwards from $(\lambda_1(T),\lambda_2(T))$ with $\lambda_1(T)>0$ and $\lambda_2(T)<0$. From \eqref{eq: lambda1_b}, for $\lambda_1<\lambda_2$, $\lambda_1'<0$, thus $\lambda_1$ is decreasing and for the semiplane $\lambda_1>\lambda_2$ we have that $\lambda_1$ is increasing. Also, from \eqref{eq: lambda2_b},  for $\lambda_1=\lambda_2 >0$ we have that $\lambda_2$ is increasing and for $\lambda_1=\lambda_2 <0$, $\lambda_2$ is decreasing. Finally, from \eqref{eq: lambda1_b} and \eqref{eq: lambda2_b}, for $\lambda_2=0$ and $\lambda_1>0$, both $\lambda_1$ and $\lambda_2$ are increasing.

Thus, since the end time conditions are on the region of the phase diagram with $\lambda_1>\lambda_2$ and $\lambda_2<0$ we have that the solution backwards in time moves to the left where $\lambda_1$ decreases and $\lambda_2$ keeps being negative. At some point in time it could cross the semiline $\lambda_1=\lambda_2<0$ (note that there cannot be touch points). If the solution crosses this line it cannot cross the semiline $\lambda_1=\lambda_2>0$ for a previous time and thus it stays on the region $\lambda_1<\lambda_2$ for all previous times.

From the definition of $\phi$ \eqref{eq: phi}, for $\lambda_1<\lambda_2$, we have $\phi>0$. For $\lambda_1=\lambda_2$, $\phi=\kappa + \lambda_3(T) \ge 0$ and for $\lambda_1>\lambda_2$, $\phi$ could become negative.  Since $\phi'(t)=\gamma x^*(t) y^*(t) \lambda_1(t)$, we see that for $t\in [0,T]$ such that $(\lambda_1(t),\lambda_2(t))$ is on the region $\lambda_1>\lambda_2$ the function $\phi$ decreases for such $t's$ with $\lambda_1(t)<0$ and increases for $\lambda_1(t)>0$. Also, let $s_0,s_1 \in [0,T]$ such that $t_0<s_1$, $\lambda_1(s_0)=\lambda_2(s_0)<0$ and $\lambda_2(s_1)<\lambda_1(s_1)=0$, then $\phi(s_0)=\kappa+\lambda_3(T)$, $\phi$ reaches the minimum value on $[t_0,T]$ at $s_1$ and $\phi(T)=\kappa+\lambda_3(T)-\frac{\gamma y^*(T) x_{\infty}}{1-\sigma_0 x_{\infty}} < \phi(s_0)$.

Thus, we conclude that $\phi$ has at most two zeros on $[0,T]$  and the proof is finished.
\end{proof}

\begin{theorem} \label{te: dos_saltos}
Let $(x^*,y^*,v^*,\sigma^*)$ be an optimal process, then:
\begin{equation} \label{eq: sigma_2_saltos}
\sigma^*(t)= \left\{ \begin{array}{ll}
\sigma_2& \text{for }0\le t < t_1 \\
\sigma_1 & \text{ for } t_1 \le t < t_1+\eta\\
\sigma_2 & \text{ for } t_1+\eta \le t \le T
\end{array}\right.
\end{equation}
with $0\le \eta\le \tau$.
\end{theorem}

\begin{proof}
Since the optimal control must be bang-bang satisfying \eqref{eq: bang bang}, from Lemma
\ref{le: dos saltos} it has at most two jumps and from \eqref{eq: int_restriction}, it takes the
value $\sigma_1$ at most for $\tau$ time. Thus the proof is completed. \end{proof}

As a consequence of theorem \ref{te: dos_saltos} we have that if $(x^*,y^*,v^*,\sigma^*)$ is an optimal process, then the optimal control $\sigma^*$ is a piecewise constant function having at most two jumps and therefore its unique associated state $(x^*,y^*,v^*)$ is a piecewise continuously differentiable function. 

\begin{lemma} \label{le: termina_sigma0}
Let  $(x^*,y^*,v^*,\sigma^*)$ be an optimal process on $[0,T]$. Assume that for $x_{\infty}=x_{\infty}(x^*(T),y^*(T),\sigma_0)$ we have $\kappa (1-\sigma_0 x_{\infty})<\gamma y^*(T)x_{\infty}$. If there exists $\delta>0$ such that $\sigma^*(t)=\sigma_2$ for all $t\in [T-\delta,T]$ then $\int_0^T \sigma^*(t)=\sigma_1 \tau + \sigma_2 (T-\tau)$.
\end{lemma}

\begin{proof}
Let $\sigma^*$ be an optimal control on $[0,T]$ such that $\sigma^*(t)=\sigma_2$ for all
$t\in [T-\delta,T]$ for $\delta>0$.  Assume that  $\int_0^T \sigma^*(t)dt > \sigma_1
\tau+(T-\tau)\sigma_2$, then from  \eqref{eq: phi} and using \eqref{eq: lambda1T_b},
\eqref{eq: lambda2T_b} and \eqref{eq: lambda3T_b}, we get
\begin{align} \label{phiTesnegativo}
\phi(T)=\kappa -\gamma  y^*(T) \frac{x_{\infty}}{1-\sigma_0 x_{\infty}} <0,
\end{align}
contradicting \eqref{eq: bang bang}.
\end{proof}

\section{Characterization of the optimal control}
\label{se: caracterization_optimal}

In this section we will give the main theorem of the article, that characterizes the switching
times $t_1$ and {\color{red} $t_1+\eta$} (where $t_1$ is the beginning of the lock-down and $\eta$ is its
duration) for an optimal control.

Let us consider the compact set
$$R=\left\{ (t_1,\eta) \in \mathbb{R}^2: \,0 \le \eta \le \tau, \,\, 0\le t_1 \le T-\eta \right\}. $$
\begin{center}
\begin{figure}[ht!]
	\begin{center}
		\includegraphics[width=5cm]{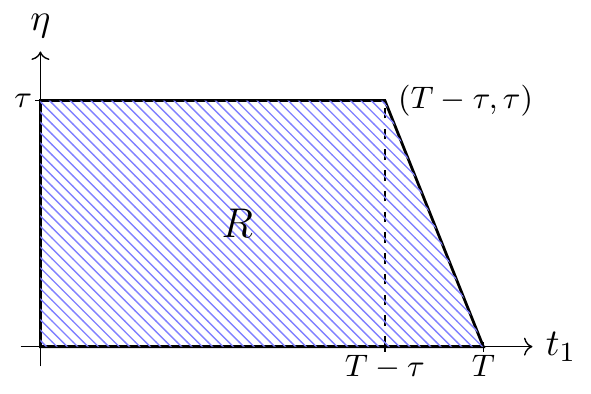}
		\caption{Graphic of set R.}
		\label{trape}
	\end{center}
\end{figure}
\end{center}
Given $(t_1,\eta)\in R$, for simplicity of notation, we will denote $t_0=0, t_2=t_1+\eta$, 
$t_3=T$.

Moreover,  given $(t_1,\eta)\in R$, we will denote $\Psi(s,t,x,y,\sigma)$ the solution of equation \eqref{eq: SIR_iso_x}-\eqref{eq: SIR_iso_y}  for $s\ge t$, $\sigma \in \left\{ \sigma_1,\sigma_2\right\}$ and initial data  $(x,y)\in {\cal D}$ at time $t$ and
\begin{subequations} \label{eq: xi_yi}
\begin{align}
(x_1,y_1)=\Psi(t_1,t_0,x_0,y_0,\sigma_2) \label{eq: x1_y1} \\
(x_2,y_2)=\Psi(t_2,t_1,x_1,y_1,\sigma_1). \label{eq: x2_y2}
\end{align}
\end{subequations}
Then, if we call $(x_{t_1,\eta},y_{t_1,\eta})$, the solution of equation \eqref{eq: SIR_iso_x}-\eqref{eq: SIR_iso_y} associated to the control $\sigma$ given by equation \eqref{eq: sigma_2_saltos} with initial data $(x(0),y(0))=(x_0,y_0)$, we have that
\begin{equation} \label{eq: xy_t_eta}
(x_{t_1,\eta}(s),y_{t_1,\eta}(s))= \left\{ \begin{array}{ll}
\Psi(s,t_0,x_0,y_0,\sigma_2)& \text{for } 0\le s \le t_1 \\
\Psi(s,t_1,x_1,y_1,\sigma_1) & \text{ for } t_1 < s \le t_2\\
\Psi(s,t_2,x_2,y_2,\sigma_2)& \text{ for } t_2< s \le T .
\end{array}\right.
\end{equation}

From theorem \ref{te: dos_saltos} we need to determine the maximum of the function
 \begin{equation} \label{eq: Jteta}
J(t_1,\eta)=x_{\infty}(x_{t_1,\eta}(T),y_{t_1,\eta}(T),\sigma_0) + \kappa (\sigma_1 \eta + \sigma_2 (T-\eta))
\end{equation}
on the compact set $R$.

In order to do that, we need to compute the derivatives of $J$.

After some computations (see \eqref{eq: derivJ_resp_t_inicial} and \eqref{eq: final_para_derivJ_respt} from Supplement) we obtain that
\begin{align} \label{eq: derivJ_resp_t_final}
\frac{\partial J}{\partial t_1}  (t_1,\eta)&=
\left[\frac{ \gamma^2 (\sigma_2-\sigma_1) y_{t_1,\eta}(T) y_{t_1,\eta}(t_1) x_{\infty,t_1,\eta}}{(1-\sigma_0 x_{\infty,t_1,\eta}) } \right] \nonumber \\
& \cdot \left[\int_{t_1}^{t_2} \frac{\sigma_0 x_{t_1,\eta}(r)-1}{y_{t_1,\eta}(r)} dr -\gamma y_2 \int_{t_2}^{T} \frac{\sigma_0 x_{t_1,\eta}(r)-1}{y_{t_1,\eta}(r)} dr  \int_{t_1}^{t_2} \frac{\sigma_2 x_{t_1,\eta}(r)-1}{y_{t_1,\eta}(r)} dr \right]
\end{align}
and
\begin{align} \label{eq: derivJ_resp_eta_final}
\frac{\partial J}{\partial \eta}  (t_1,\eta)= & \frac{\gamma x_{\infty,t_1,\eta}}{1-
\sigma_0 x_{\infty,t_1,\eta}} y_{t_1,\eta}(t_2)(\sigma_2-\sigma_1)
\left(1-(\sigma_0-\sigma_2)y_{t_1,\eta}(T)\gamma \int_{t_2}^T \frac{x_{t_1,\eta}(r)}{y_{t_1,\eta}(r)} dr\right) \nonumber   \\
& -\kappa (\sigma_2-\sigma_1)\nonumber   \\
= &\frac{\gamma  x_{\infty,t_1,\eta}}{1-\sigma_0 x_{\infty,t_1,\eta}} (\sigma_2-\sigma_1)y_{t_1,\eta}(T)\left(1-\gamma y_{t_1,\eta}(t_2)
 \int_{t_2}^T \frac{\sigma_0 x_{t_1,\eta}(r)-1}{y_{t_1,\eta}(r)} dr\right)   \\ & -\kappa (\sigma_2-\sigma_1)\nonumber   \\ \nonumber
\end{align}
where $x_{\infty,t_1,\eta}=x_{\infty}(x_{t_1,\eta}(T),y_{t_1,\eta}(T),\sigma_0)$.
Note that for $(t_1,\eta)\in R$, $\frac{\partial J}{\partial \eta}  (t_1,\eta)>0$ if and only if
\begin{align} \label{eq: cond_kappa}
\frac{ x_{\infty,t_1,\eta}}{1-\sigma_0 x_{\infty,t_1,\eta}}\gamma y_{t_1,\eta}(T)\left(1-\gamma y_{t_1,\eta}(t_2) \int_{t_2}^T \frac{\sigma_0 x_{t_1,\eta}(r)-1}{y_{t_1,\eta}(r)} dr\right)>\kappa.
\end{align}

In the results given in sections \ref{se: sigma2igual0_kappa0} and \ref{se: kappa_posit} we
will assume that $\frac{\partial J}{\partial \eta}(t_1,\eta)>0$ for all $(t_1,\eta)\in R$ and
therefore in the
following remark we analyse the derivatives of $J$ restricted to the superior border of $R$
(see equations \eqref{eq: derivJ_resp_t_final_borde} and \eqref{eq:
derivJ_resp_eta_final_borde}) which will be used later.

\begin{remark} \label{re: cond_max_arriba}
Assume that $\frac{\partial J}{\partial \eta}(t_1,\eta)>0$ for all $(t_1,\eta)\in R$, then the maximum value of $J$ on $R$ must be attained at the superior border
\begin{align} \label{def: P}
P=\left\{ (t_1,\tau), \, t_1 \in [0,T-\tau] \right\} \cup \left\{ (t_1,T-t_1), \, t_1 \in [T-\tau, T] \right\} .
\end{align}
Thus, in this case, for $(t_1,\tau)$ with $t_1\in  [0,T-\tau)$ we have $t_2=t_1+\tau$ and for $(t_1,T-t_1)$ with $t_1\in[T-\tau,T]$ we have $t_2=T$, and the control is as in Figure \ref{fig: bb12}

\bigskip

\begin{figure}[ht!]
	\centering
	\subfigure[{\it \footnotesize For $0\le t_1<T-\tau, \quad t_2=t_1+\tau$. }]{
		\includegraphics[angle=0,width=7cm]{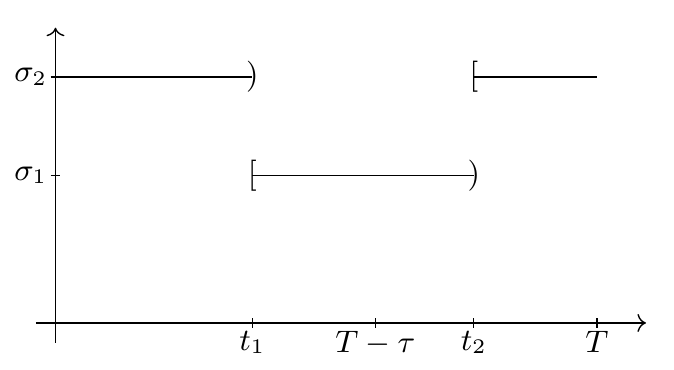}
		\label{fig: bb1}
	}
	\subfigure[{\it \footnotesize For $T-\tau\le t_1<T, \quad t_2=T$.}]{
		\includegraphics[angle=0,width=7cm]{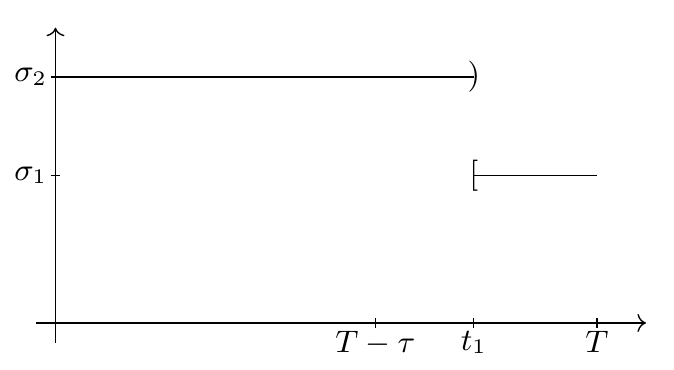}
		\label{fig: bb2}
	}
	\caption{\it \footnotesize Control for $(t_1,\eta)\in P.$  In \ref{fig: bb1} we have a mild-strict-mild quarantine in the intervention interval $[0,T]$. In \ref{fig: bb2} we have a mild-strict quarantine in the intervention interval $[0,T]$.}
	\normalsize\label{fig: bb12}
\end{figure}

Moreover, for $t\in [0,T]$, we define the continuous function
\begin{equation} \label{eq: omega}
w(t)= \left\{ \begin{array}{ll}
\int_{t}^{t+\tau } \frac{\sigma_0 x_{t,\tau}(r)-1}{y_{t,\tau}(r)} dr -\gamma y_{t,\tau}(t+\tau) \int_{t+\tau}^{T} \frac{\sigma_0 x_{t,\tau}(r)-1}{y_{t,\tau}(r)} dr  \int_{t}^{t+\tau } \frac{\sigma_2 x_{t,\tau}(r)-1}{y_{t,\tau}(r)} dr  & \text{for }0\le t \le T-\tau \\
& \\
 \int_{t}^{T} \frac{\sigma_0 x_{t,T-t}(r)-1}{y_{t,T-t}(r)}dr & \text{ for } T-\tau < t \le T,
\end{array}\right.
\end{equation}
and for $t\in[T-\tau,T]$, we define \\
\begin{align}\label{eq: alpha}
\alpha(t)&=
\frac{1}{\gamma y_{t,T-t}(t)}\left( 1-\kappa
 \frac{1-\sigma_0 x_{\infty,t,T-t}}{ \gamma y_{t,T-t}(T) x_{\infty,t,T-t}  } \right).
 \end{align}
Then we have that for $(t_1,\eta)\in P$
\begin{align} \label{eq: derivJ_resp_t_final_borde}
\frac{\partial J}{\partial t_1}  (t_1,\eta)&=  \frac{ x_{\infty,t_1,\eta}}{(1-\sigma_0 x_{\infty,t_1,\eta}) } \gamma^2 (\sigma_2-\sigma_1)  y_{t_1,\eta}(T) y_{t_1,\eta}(t_1) w(t_1),\end{align}
and for $(t_1,T-t_1)\in P$ with $t_1\in [T-\tau,T]$
\begin{align} \label{eq: derivJ_resp_eta_final_borde}
\frac{\partial J}{\partial \eta}  (t_1,T-t_1)&= \frac{ x_{\infty,t_1,\eta}}{1-\sigma_0 x_{\infty,t_1,\eta}}\gamma  (\sigma_2-\sigma_1)y_{t_1,\eta}(T)-\kappa (\sigma_2-\sigma_1)\nonumber   \\
&=\frac{ x_{\infty,t_1,\eta}}{1-\sigma_0 x_{\infty,t_1,\eta}} \gamma^2(\sigma_2-\sigma_1)y_{t_1,\eta}(T)  y_{t_1,\eta}(t_1)  \alpha(t_1).
\end{align}
We consider $\tilde{J}$ the continuous function defined as the restriction of $J(t_1,\eta)$ to $P$, that is
\begin{equation} \label{eq: Jtilde}
\tilde{J}(t_1)= \left\{ \begin{array}{ll}
J(t_1,\tau)  \text{ for } t_1\in [0,T-\tau], \\
& \\
J(t_1,T-t_1) \text{ for } t_1\in [T-\tau,T].
\end{array}\right.
\end{equation}
From \eqref{eq: derivJ_resp_t_final_borde} we have that
\begin{align} \label{eq: derivJ1sigma2}
\tilde{J}' (t_1)&=  \gamma^2 (\sigma_2-\sigma_1) \frac{ x_{\infty,t_1,\tau}}{(1-\sigma_0 x_{\infty,t_1,\tau})} y_{t_1,\tau}(t_1)   y_{t_1,\tau}(T) w(t_1) \quad \text{ for } \quad t_1\in [0,T-\tau),
\end{align}
and from \eqref{eq: derivJ_resp_t_final_borde} and \eqref{eq: derivJ_resp_eta_final_borde}, we get
\begin{align}
 \tilde{J}' (t_1)&= \frac{d J}{dt_1}  (t_1,T-t_1)-\frac{d J}{d\eta}  (t_1,T-t_1) \nonumber \\
& = \gamma^2 (\sigma_2-\sigma_1)  \frac{x_{\infty,t_1,T-t_1}}{1-\sigma_0 x_{\infty,t_1,T-t_1}} y_{t_1,T-t_1}(t_1) y_{t_1,T-t_1}(T)  \left( w(t_1)- \alpha(t_1) \right)\label{eq: J2_alpha_w22} 
\end{align}
for $t_1\in  (T-\tau,T]$.
\end{remark}

In the next subsection, we prove the main result of this article  (theorem \ref{te:
control_optimo_sigma1gral}) for the case $\sigma_2 = \sigma_0$, $\sigma_1\ge 0$ and
$\kappa = 0$. Then, in subsection \ref{se: sigma1igual0_kappa0} we derive the result for
$\sigma_1=0$ (Corollary \ref{co: sigma1es0}) in order to compare our result with the one
obtained in \cite{ketcheson2020optimal}.

\subsection{Case $\sigma_2 = \sigma_0$ and $\kappa = 0$.}
\label{se: sigma2igual0_kappa0}

For $\sigma_2=\sigma_0$, from \eqref{eq: cota_integral} with $i= 2$, we obtain
\begin{align} \label{eq: w_2_0}
w(t) =
\begin{split}
\begin{cases}
\dfrac{y_{t,\tau}(t+\tau)}{y_{t,\tau}(T)}\int_{t}^{t+\tau} \frac{\sigma_0 x_{t,\tau}(r)-1}{y_{t,\tau}(r)} dr  \quad \text{ for } 0\leq t \leq T-\tau\\\\
\int_{t}^{T} \frac{\sigma_0 x_{t,T-t}(r)-1}{y_{t,T-t}(r)} dr \quad \text{ for } T-\tau < t \le T.
\end{cases}
\end{split}
\end{align}
In addition, from \eqref{eq: derivJ_resp_eta_final} and \eqref{eq: cota_integral} and using $\kappa=0$ we have that
\begin{align} \label{eq: derivJ_resp_eta_final_sigma2=0}
\frac{\partial J}{\partial \eta}  (t_1,\eta)&= \frac{ x_{\infty,t_1,\eta}}{1-\sigma_0 x_{\infty,t_1,\eta}}\gamma  (\sigma_0-\sigma_1)y_{t_1,\eta}(t_1+\eta)>0\end{align}
for all $(t_1,\eta)\in R$.

In the next remark we discuss the sign of $w(t)$ for $t\in[0,T-\tau]$ when $\sigma_2=\sigma_0$.

\begin{center}
 \begin{figure}[ht!]
 	\begin{center}
 		\includegraphics[width=7cm]{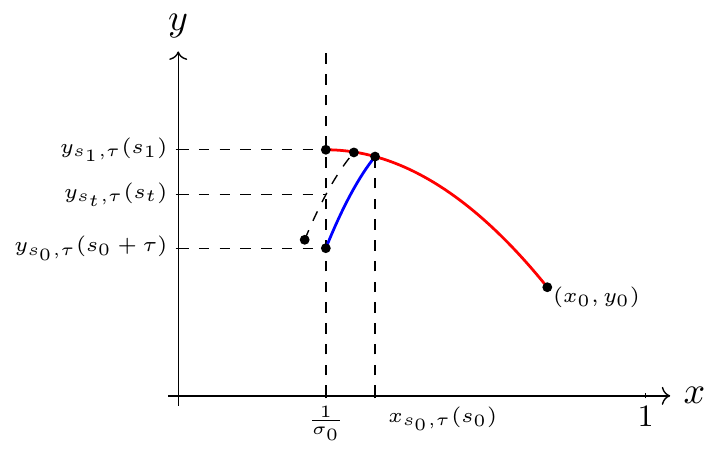}
 		\caption{Trajectories $(x_{t,\tau},y_{t,\tau})$ for $t\in[s_0,s_1]$ when $\sigma_2 = \sigma_0$.}
 		\label{fig: s0s1}
 	\end{center}
 \end{figure}
\end{center}

\begin{remark} \label{re: Gdecrece} Given $(x_0,y_0)\in {\cal D}$ with $x_0>1/\sigma_0$,
assume there exists $s_1\in [0,T-\tau]$ such that the solution $\Psi_1(s_1,0,x_0,y_0,\sigma_2)=1/\sigma_0$, that is $x_{s_1,\tau}(s_1)=1/\sigma_0$ (red line in Fig. \ref{fig: s0s1}). Then, for
$t \in [s_1,T-\tau]$, $x_{t,\tau}(t)\le 1/\sigma_0$ and therefore $x_{t,\tau}(r)<1/\sigma_0$ for $r\in (t,t+\tau]$ implying $w(t)<0$.
Additionally, assume there exists $s_0\in [0,T-\tau] $ such that $x_{s_0,\tau}(s_0+\tau)=1/\sigma_0$ (blue line). In this case, it is clear that  $s_0<s_1$ and also that for all $t\in (s_0,s_1)$ there exists a unique $s_t \in [t,t+\tau]$ such that $x_{t,\tau}(s_t)=1/\sigma_0$.
Moreover we can conclude that for $t \in[0,s_0]$, $x_{t,\tau}(r)\ge 1/\sigma_0$ for all $r\in (t,t+\tau)$ and therefore $w(t)>0$.
If $s_1$ defined before does not exist, that is $x_{s,\tau}(s)>1/\sigma_0$ for all $s\in[0,T-\tau]$, then we take $s_1=T-\tau$. Likewise, if $s_0$ does not exist, that is for all $s\in[0,T-\tau]$, $x_{s}(s+\tau)<1/\sigma_0$, then we take $s_0=0$ and conclude that in either case, the sign of $w(t)$ for $t\ge 0$,  is determined  in the complement of $[s_0,s_1]$.
\end{remark}
In the next lemma we prove that $w$ is a decreasing function on the interval $(s_0,s_1)$
introduced in remark \ref{re: Gdecrece}.

\begin{lemma} \label{le: zdecrece}
For $t\in [0,T-\tau]$, we define
\begin{align} \label{eq: z}
z(t)=\int_{t}^{t+\tau} \frac{\sigma_0 x_{t,\tau}(r)-1}{y_{t,\tau}(r)} dr.
\end{align}
Let $s_0,s_1\in [0,T-\tau]$ be given by remark \ref{re: Gdecrece}. Then $z$ is a decreasing function on $(s_0,s_1)$. Moreover, $w(t)>0$ for $t<s_0$, $w(t)<0$ for $t>s_1$ and consecuently $w$ changes sign at most once on $[0,T-\tau]$.\end{lemma}
\begin{proof}
Since the duration of the strict quarantine is $\tau$ fixed, for simplicity of notation in this proof we neglect the subindex $\tau$ from solutions $x$ and $y$.
Given $t\in (s_0,s_1)$, we define the auxiliary functions
\begin{subequations}  \label{eq: func_aux}
\begin{align}
g_t(s)&=\sigma_0 \sigma_1 x_t(s) y_t(s)+(\sigma_0 x_t(s)-1)(\sigma_1 x_t(s)-1),  \label{eq: g} \\
f_t(s)&=\frac{\sigma_0 x_t(s)-1}{y_t(s)}, \label{eq: f}\\
i_t(s)&=\frac{x_t(s)}{y_t(s)}(\sigma_0 x_t(s)+\sigma_0 y_t(s)-1)+\gamma  g_t(s) \int_t^s \frac{x_t(r)}{y_t(r)}dr. \label{eq: i}
\end{align}
\end{subequations}
By computing the derivative for $s\in(t,t+\tau)$ we obtain
\begin{subequations}  \label{eq: deriv_func_aux}
\begin{align}
g_t'(s)&=-\gamma \sigma_1^2 x_t(s) y_t(s)(\sigma_0 x_t(s)+\sigma_0 y_t(s)-1), \label{eq: derivg} \\
f_t'(s)&=-\gamma \frac{g_t(s)}{y_t(s)}, \label{eq: derivf}\\
i_t'(s)&=-\gamma \sigma_1 x_t(s)\left(\gamma \sigma_1   y_t(s) \int_t^s \frac{x_t(r)}{y_t(r)}dr +1\right)(\sigma_0 x_t(s)+\sigma_0 y_t(s)-1),\label{eq: derivi}
\end{align}
\end{subequations}
and we have that
\begin{align}  \label{eq: z_prima2}
z'(t)=&f_t(t+\tau) - f_t(t) \\
&+ \gamma (\sigma_1-\sigma_0) y_1 \int_{t}^{t+\tau} \left[ \frac{x_t(s)}{ y_t^2(s)}(\sigma_0 x_t(s)+\sigma_0 y_t(s)-1)
+\gamma \left( \int_{t}^s \frac{x_t(r)}{y_t(r)} dr\right)  \frac{g_t(s)}{y_t(s)} \right] ds \nonumber \\
= &f_t(t+\tau) - f_t(t)+ \gamma (\sigma_1-\sigma_0) y_1 \int_{t}^{t+\tau}  \frac{i_t(s)}{ y_t(s)} ds  \nonumber
\end{align}
Note that both $g_t'$ and $i_t'$ have the opposite sign of $(\sigma_0 x_t(s)+\sigma_0
y_t(s)-1)$. First, assume $\sigma_0 x_t(s)+\sigma_0 y_t(s)-1>0$ for all $s\in (t,t+\tau)$, then
from  \eqref{eq: derivg}, $g_t$ is a decreasing function. Moreover, since  $\sigma_1
x_t(t+\tau)-1<0$ and $\sigma_0 x_t(t+\tau)-1<0$ for $t\in (s_0,s_1)$, we deduce that
$g_t(s)>g_t(t+\tau)>0$. In addition, from \eqref{eq: i} we obtain $i_t$ is positive.

 On the
other side, from the fact that $x+y$ is a decreasing function, if we assume that there exists
$s_3\in [t,t+\tau]$ such that $\sigma_0 x_t(s_3)+\sigma_0 y_t(s_3)=1$ we have that for
$s<s_3$, $ \sigma_0 x_t(s)+\sigma_0 y_t(s)-1>0$ and for $s>s_3$, $ \sigma_0
x_t(s)+\sigma_0 y_t(s)-1<0$. Therefore,
$g_t$ and $i_t$ attains a global minimum on $[t,t+\tau]$ at $s_3$ and thus $g_t(s)\ge g_t(s_3)=\sigma_0 y_t(s_3)>0$ and $i_t(s)\ge i_t(s_3)=\gamma g_t(s_3)\int_{t}^{s_3} \frac{x_t(r)}{y_t(r)} dr>0$ for all $s\in [t,t+\tau]$ . Furhtermore, from \eqref{eq: derivf} we have that $f_t(s)$ is also a decreasing function.

Thus, we have proved that for $t\in (s_0,s_1)$, $f_t(s)$ is decreasing on $(t,t+\tau)$ and
$i_t(s)>0$ for all $s\in [t,t+\tau]$, yielding from \eqref{eq: z_prima2} that $z'(t)<0$ for all
$t\in (s_0,s_1)$.

Finally, from remark \ref{re: Gdecrece} we deduce that  $w$ changes sign at most once on $[0,T-\tau]$ (see figure \ref{wdecrece}).
\begin{center}
	\begin{figure}[ht!]
		\begin{center}
			\includegraphics[width=7cm]{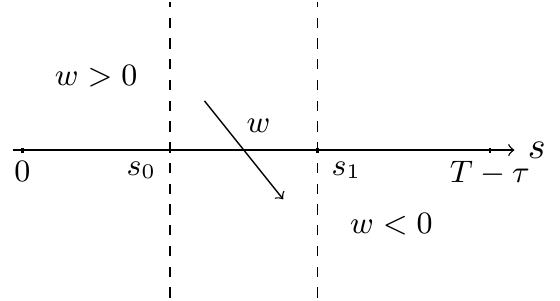}
			\caption{Behaviour of $w(t)$ when $\sigma_0 = \sigma_2$.}
			\label{wdecrece}
		\end{center}
	\end{figure}
\end{center}
\end{proof}

In the next theorem we assume that $x_0<1/\sigma_1$. This condition is always satisfied for $\sigma_1<1$.

\begin{theorem} \label{te: control_optimo_sigma1gral}
Let $0\le \sigma_1<\sigma_2=\sigma_0$ with $\sigma_1<1$, $\kappa = 0$ and $w$ be given by \eqref{eq: w_2_0}. Then the optimal control is unique and is given by
\begin{equation}
\sigma^*(s)= \left\{ \begin{array}{ll}
\sigma_0& \text{for }0\le s <t^*,\\
\sigma_1 & \text{ for } t^* \le s < t^*+\eta,\\
\sigma_0 & \text{ for } t^*+\eta \le s <T,
\end{array}\right.
\end{equation}
where
\begin{enumerate}
\item For $w(0)\le 0$: $t^*=0$ and $\eta=\tau$.
\item For  $w(0)>0$ and   $w(T-\tau)\le 0$: $t^*=\overline{t}$ and $\eta=\tau$ where $\overline{t}$ is the unique value on $[0,T-\tau]$ such that $w(\overline{t})=0$.
\item For $0<w(T-\tau)\le\dfrac{1}{\gamma y_{T-\tau,\tau}(T-\tau)}$:  $t^*=T-\tau$ and $\eta=\tau$.
\item For $w(T-\tau)>\dfrac{1}{\gamma y_{T-\tau,\tau}(T-\tau)}$: $t^*=\tilde{t}$ where $\tilde{t}$ is the unique value on $[T-\tau,T]$ such that $w(\tilde{t})=\dfrac{1}{\gamma y_{\tilde{t},T-\tilde{t}}(\tilde{t})}$ and $\eta=T-\tilde{t}$.
\end{enumerate}
 \end{theorem}

\begin{proof}
From equation \eqref{eq: derivJ_resp_eta_final_sigma2=0} and Remark \ref{re: cond_max_arriba}
the maximum value of $J$ on $R$ must be attained at the superior border $P$ defined in \eqref{def: P}.  Therefore,
from \eqref{eq: derivJ1sigma2} and \eqref{eq: J2_alpha_w22} we have
\begin{align} \label{eq: derivJ1sigma1}
\tilde{J}' (t_1)&=  \gamma^2 (\sigma_0-\sigma_1) \frac{ x_{\infty,t_1,\tau}}{(1-\sigma_0 x_{\infty,t_1,\tau})}  y_{t_1,\tau}(t_1) y_{t_1,\tau}(T) w(t_1) \quad \text{ for } \quad t_1\in [0,T-\tau),
\end{align}
and using that for $\kappa=0$, $\alpha(t_1)=\dfrac{1}{\gamma y_{t_1,T-t_1}(t_1)}$, then
\begin{align}
 \tilde{J}' (t_1)&= \gamma^2 (\sigma_0-\sigma_1)  \frac{x_{\infty,t_1,T-t_1}}{1-\sigma_0 x_{\infty,t_1,T-t_1}} y_{t_1,T-t_1}(t_1)y_{t_1,T-t_1}(T)  \left( w(t_1)- \dfrac{1}{\gamma y_{t_1,T-t_1}(t_1)} \right)\label{eq: J2_alpha_w} \end{align}
for $t_1\in  (T-\tau,T]$.
Note that from \eqref{eq: cota_integral_ij} with $i=1$ and $j=0$, for $t\in [T-\tau,T]$ it holds the identity
\begin{align} \label{eq: factor_de_w2}
\gamma  y_{t_1,T-t_1}(T)  \left( w(t_1)- \dfrac{1}{\gamma y_{t_1,T-t_1}(t_1)} \right) =  \gamma (\sigma_0-\sigma_1)h(t_1)-1
\end{align}
 where
\begin{equation}\label{def: h}
h(t) = y_{t,T-t}(T)\int_{t}^{T}\frac{ x_{t,T-t}(s)}{y_{t,T-t}(s)}ds
\end{equation}
is a decreasing function in $[T-\tau,T]$. In fact,
using that $u_1(s)$ is a positive function (see equation \eqref{eq: u_ti}) and equation \eqref{eq: deriv_x2tini} it is easy to see that
\begin{align*}
\frac{d}{dt_1} \left(\frac{y_{t_1}(T)}{y_{t_1}(s)} \right)<0
\end{align*}
and therefore, $h'(t_1) <0$ for $t_1\in (T-\tau,T)$.
From  \eqref{eq: J2_alpha_w} and \eqref{eq: factor_de_w2} we also have that for $t_1\in [T-\tau,T]$
\begin{align}
 \tilde{J}' (t_1)&=\gamma (\sigma_0-\sigma_1)   \frac{x_{\infty,t_1,T-t_1}}{1-\sigma_0 x_{\infty,t_1,T-t_1}} y_{t_1,T-t_1}(t_1)  \left( \gamma (\sigma_0-\sigma_1)h(t_1)-1\right) \label{eq: J2_alpha_w2}
\end{align}

We consider the following cases:
\begin{enumerate}
\item If $w(0)\le 0$, then from Lemma \ref{le: zdecrece}, $w(t_1)<0$ for all $t_1 \in (0,T-\tau]$. Thus, from \eqref{eq: derivJ1sigma1} we have
\begin{align*}
 \tilde{J}' (t_1) &
\le 0 \text{ for all } t_1\in[0,T-\tau).
\end{align*}
Moreover, using that $w(T-\tau)<0$, the positivity of $\gamma y_{T-\tau,\tau}(T-\tau)$ and \eqref{eq: factor_de_w2}, we obtain that $h(T-\tau)<\dfrac{1}{\gamma(\sigma_0-\sigma_1)}$
and being $h$ a decreasing function, from \eqref{eq: J2_alpha_w2} we deduce that
\begin{align*}
 \tilde{J}' (t_1)&
<0  \text{ for all } t_1\in(T-\tau,T].
\end{align*}
Therefore $t^*=0$.

\item Since $w(0)>0$ and $w(T-\tau)\leq 0$, from Lemma \ref{le: zdecrece} there exists an unique $\overline{t}\in (0,T-\tau]$ such that  $w(\overline{t})=0$, $w(t_1)>0$ for all $t_1\in [0,\overline{t})$ and $w(t_1)<0$  for all $t_1\in (\overline{t},T]$.\\
Moreover, since $w(T-\tau)\leq 0$, in the same way as for the previous item, we have
\begin{align*}
 \tilde{J}' (t_1)
&   <0 \text{ for all } t_1\in(T-\tau,T]
\end{align*}
and from \eqref{eq: derivJ1sigma1}, we obtain
\begin{align*}
 \tilde{J}' (t_1)&\geq 0 \quad \text{ for } t_1<\overline{t}  \quad\text{ and }\quad \tilde{J}' (t_1)< 0 \quad \text{ for } t_1>\overline{t},
\end{align*}
concluding that $t^*=\overline{t}$.

\item
Since
\begin{align*}
0<w(T-\tau)&<\dfrac{1}{\gamma y_{T-\tau,\tau}(T-\tau)},
\end{align*}
from Lemma \ref{le: zdecrece}, $w(t_1)>0$ for all $t_1\in[0,T-\tau]$. On the other side,
since $w(T-\tau) <\dfrac{1}{\gamma y_{T-\tau,\tau}(T-\tau)}$, from \eqref{def: h} we have that $h(T-\tau)\le \frac{1}{\gamma(\sigma_0-\sigma_1)}$ and using that $h$ is a continuous and decreasing function, we obtain that $h(t_1)<\frac{1}{\gamma(\sigma_0-\sigma_1)}$ for all $t_1\in (T-\tau,T]$. Thus,
\begin{align*}  \tilde{J}' (t_1)> 0 \quad \text{ for } t_1\in [0,T-\tau)\quad \text{ and } \quad  \tilde{J}' (t_1)< 0 \quad \text{ for } t_1\in (T-\tau,T].
\end{align*}
 Therefore, $t^*=T-\tau$ and $\eta=\tau$.

\item Since
\begin{align*}
w(T-\tau)&>\dfrac{1}{\gamma y_{T-\tau,\tau}(T-\tau)} >0,
\end{align*}
from Lemma \ref{le: zdecrece}, we have $w(t_1)>0$ for all  $t_1\in[0,T-\tau]$.
On the other side, since $w(T-\tau)>\dfrac{1}{\gamma y_{T-\tau,\tau}(T-\tau)}$, then $h(T-\tau)>\frac{1}{\gamma(\sigma_0-\sigma_1)}$ and using that $h(T)=0$ and $h$ is a continuous and decreasing function, there exists a unique $\tilde{t}$ such that $h(\tilde{t})=\frac{1}{\gamma(\sigma_0-\sigma_1)}$, $h(t)>\frac{1}{\gamma(\sigma_0-\sigma_1)}$ for $t\in[T-\tau,\tilde{t})$ and $h(t)<\frac{1}{\gamma(\sigma_0-\sigma_1)}$ for $t\in(\tilde{t},T]$.
Therefore,
\begin{align*} \tilde{J}' (t_1)> 0  \quad \text{ for }  t_1\in[0,\tilde{t}) \quad \text{ and } \quad  \tilde{J}' (t_1)< 0 \quad \text{ for } t_1\in(\tilde{t},T].
\end{align*}
Consequently, $t^*=\tilde{t}$ and $\eta = T-\tilde{t}$.
\end{enumerate}
\end{proof}

\subsection{Case $\sigma_1= 0 $, $\sigma_2 = \sigma_0$ and $\kappa = 0$.}
\label{se: sigma1igual0_kappa0}

Let $w(t) $ defined as in \eqref{eq: omega}. Note that for $\sigma_1 = 0 $ and $\sigma_2= \sigma_0$, from \eqref{eq: cota_integral_ij} with $i=j = 2$, we obtain

\begin{align}\label{eq: w_0}
w(t) =
\begin{split}
\begin{cases}
\dfrac{y_{t,\tau}(t+\tau)(\sigma_0 x_{t,\tau}(t)-1)}{\gamma y_{t,\tau}(T) y_{t,\tau}(t)}(e^{\gamma \tau}-1) \quad \text{ for } 0\leq t \leq T-\tau\\\\
\dfrac{\sigma_0 x_{t,T-t}(t)-1}{\gamma y_{t,T-t}(t)}(e^{\gamma(T-t)}-1) \quad \text{ for } T-\tau < t \le T.
\end{cases}
\end{split}
\end{align}

It is easy to observe that the sign of $w(t)$ on $[0,T-\tau]$ is given by $(\sigma_0 x_{t,\tau}(t)-1)$.
Moreover, $w$ changes sign at most once on $[0,T-\tau]$, going from positive to negative values.

\begin{corollary} \label{co: sigma1es0}
Let $\sigma_1=0$, $\sigma_2 = \sigma_0>0$ and $k=0$.
Then, the optimal control is unique and is given by
	
\begin{align*}
\sigma^*(s) =
\begin{split}
\begin{cases}
\sigma_0 \quad \text{ for } 0\leq s < t^*,\\
0 \quad \text{ for } t^* \le  s < t^*+\eta,\\
\sigma_0 \quad \text{ for } t^*+\eta \leq s \leq T,
\end{cases}
\end{split}
\end{align*}
where
\begin{enumerate}
\item	For $x_0 \le \frac{1}{\sigma_0}$ : $ t^*= 0$ and $\eta= \tau$.
\item  For $x_0 > \frac{1}{\sigma_0} $ and $x_{T-\tau,\tau}(T-\tau) \le \frac{1}{\sigma_0} $:  $t^* = \overline{t}$ and $\eta = \tau$, where $\overline{t}$ is the unique value on $[0,T-\tau]$ such that $x_{\overline{t},\tau}(\overline{t}) = \dfrac{1}{\sigma_0}$.
\item For $\dfrac{1}{\sigma_0} < x_{T-\tau,\tau}(T-\tau)\le \dfrac{1}{\sigma_0 (1-e^{-\gamma \tau})} $: $t^* = T-\tau$ and $\eta = \tau$.
\item For  $x_{T-\tau,\tau}(T-\tau)> \dfrac{1}{\sigma_0 (1-e^{-\gamma \tau})}$: $t^* = \tilde{t}$ and $\eta = T-\tilde{t}$, where $\tilde{t}$ is the unique value on $[T-\tau, T]$ such that $x_{\tilde{t},T-\tilde{t}}(\tilde{t})= \dfrac{1}{\sigma_0 \left(1-e^{-\gamma(T-\tilde{t}) }\right)}$.
\end{enumerate}
\end{corollary}	

\begin{proof}
The proof follows from theorem \ref{te: control_optimo_sigma1gral} using the fact that
$$\sign(w(t))=\sign(\sigma_0 x_{t,\tau}(t) -1), \quad \text{ for } t\in [0,T-\tau], $$
and
$$w(T-\tau)=\dfrac{1}{\gamma y_{T-\tau,\tau}(T-\tau)}(\sigma_0 x_{T-\tau,\tau}(T-\tau)-1)(e^{\gamma \tau}-1).$$
\end{proof}

Note that in this case if we take $\tau=T$, the corollary is reduced to only two possible
cases: $x_0 \le \dfrac{1}{\sigma_0 (1-e^{-\gamma T})}$ or  $x_0> \dfrac{1}{\sigma_0
(1-e^{-\gamma T})}$,  obtaining the same result  as Ketcheson  \cite{ketcheson2020optimal}
in Theorem 3.

\section{General case}
\label{se: kappa_posit}

In this section we study  the behaviour of optimal solutions for
 the general case when  $0\le\sigma_1<\sigma_2\le \sigma_0$ and $\kappa>0$, that is,
objective function $J$ includes the term that accounts the running cost of the control and
allows us   to account for factors  like the economic cost of intervention or heightened risks
caused by hospital overflow.

\begin{lemma} \label{le: kappagrande}
	Assume $\kappa > \dfrac{ x_{\infty,t_1,\eta}\gamma y_{t_1,\eta}(T)}{1-\sigma_0 x_{\infty,t_1,\eta}}\left(1-\gamma y_2 \int_{t_2}^T \frac{\sigma_0 x_{t,\eta}(r)-1}{y_{t,\eta}(r)} dr\right)$  for all $(t_1,\eta)\in R$, then the optimal control is given by $\sigma^*\equiv \sigma_2$.
\end{lemma}
\begin{proof}
	From \eqref{eq: derivJ_resp_eta_final}, we have that $ \frac{\partial J}{\partial \eta}  (t_1,\eta)<0$ and therefore the maximum value of $J$ on $R$ is attained at the inferior border of $R$ where $\eta=0$ and $J(t_1,0)$ is constant.
\end{proof}

In the next theorem we give a general result including both the economic cost of intervention ($\kappa >0$) and a mitigation phase different from the no intervention one, that is $\sigma_2<\sigma_0$.
In section \ref{se: simulations} we give numerical simulations supporting this result. When $\sigma=\sigma_0$ and $\kappa=0$ we recover theorem \ref{te: control_optimo_sigma1gral} proved in section \ref{se: sigma2igual0_kappa0}.
\begin{theorem} \label{th: version_general}
	Let $0\le \sigma_1<\sigma_2\le\sigma_0$ with $\sigma_1<1$, $\kappa$ satisfying \eqref{eq: cond_kappa} for all $(t_1,\eta)\in R$ and let $w$ and $\alpha$ defined as in \eqref{eq: omega} and \eqref{eq: alpha} respectively , then the optimal control is unique and is given by
	\begin{equation}
	\sigma^*(s)= \left\{ \begin{array}{ll}
	\sigma_2& \text{for }0\le s <t^*,\\
	\sigma_1 & \text{ for } t^* \le s < t^*+\eta,\\
	\sigma_2 & \text{ for } t^*+\eta \le s <T,
	\end{array}\right.
	\end{equation}
	where
	\begin{enumerate}
		\item For $w(0)\le 0$: $t^*=0$ and $\eta=\tau$.
		\item For  $w(0)>0$ and   $w(T-\tau)\le 0$: $t^*=\overline{t}$ and $\eta=\tau$ where $\overline{t}$ is the unique value on $[0,T-\tau]$ such that $w(\overline{t})=0$.
		\item For $0<w(T-\tau)\le\alpha(T-\tau)$:  $t^*=T-\tau$ and $\eta=\tau$.
		\item For $w(T-\tau)>\alpha(T-\tau)$: $t^*=\tilde{t}$ where $\tilde{t}$ is the unique value on $[T-\tau,T]$ such that $w(\tilde{t})=\alpha(\tilde{t})$ and $\eta=T-\tilde{t}$.
	\end{enumerate}
\end{theorem}

\section{Numerical simulations}
\label{se: simulations}

In this section we check numerically Corollary \ref{co: sigma1es0},
theorem \ref{te: control_optimo_sigma1gral} and theorem \ref{th: version_general}, which summarize the main results for the different
case scenarios.  We integrate the system of equations~\eqref{eq: SIR}
for different time periods $\tau$ of the hard quarantine ($\sigma=\sigma_1$) in the bang-bang control, starting at time $t$ and constrained to a maximum ending time $T$, which is the control period
($\sigma(t)=\sigma_0$ for $t > T$).  That is, the value of $\sigma(t)$ adopts the following form, depending on $t$, $\tau$ and $T$:
\begin{eqnarray}
  \sigma(r) =
  \begin{cases}
    \sigma_2 & \mbox{for $0 \le r  < t$}, \\
    \sigma_1 & \mbox{for $t \le r < t + \eta$},  \\
    \sigma_2 & \mbox{for $t+\eta \le r  < T$}, \\
    \sigma_0 & \mbox{for $r \ge T$},
  \end{cases}
  \label{sigma-t}
\end{eqnarray}
where
\begin{eqnarray}
	\eta =
	\begin{cases}
		\tau & \mbox{for $t+\tau \le T$ ~ and} \\
		T-t & \mbox{for $t+\tau > T$}.
	\end{cases}
	\label{Delta-t}
\end{eqnarray}

We start by testing the simplest case scenario $\kappa=0$,
$\sigma_2=\sigma_0$ and $\sigma_1 \ge 0$, and we then test the most general case $\kappa>0$ and $0 < \sigma_1 < \sigma_2 < \sigma_0$.

\subsection{Case $\kappa=0$, $\sigma_2 = \sigma_0$ and $\sigma_1= 0$ (Corollary \ref{co: sigma1es0})}

We first analyze the case $\kappa=0$ with a bang-bang control in the interval $[0,T]$ that consists of a soft quarantine $(\sigma_2=1.5)$ and an extremely hard and unrealistic quarantine ($\sigma_1=0$) during which there are no infections.  The other parameters in the simulations are $\gamma=0.01$ and $T=2600$, together with the initial condition $y_0=y(t=0)=10^{-6}$ and $x_0=x(t=0)=1-10^{-6}$.  We can see from Corollary \ref{co: sigma1es0} that the optimum initial time of the hard quarantine is $t^*=0$ for $x_0 \le 1/\sigma_0$ ($w(0) \le 0$), while for $x_0 > 1/\sigma_0$ ($w(0)>0$) is given by
\begin{eqnarray}
t^* =
\begin{cases}
		\overline{t} & \mbox{for $0 \le \tau \le \overline{\tau}$, where $x_{T-\tau,\tau}(T-\tau) \le \dfrac{1}{\sigma_0}$ and
			$\eta=\tau$}, \\
		T-\tau & \mbox{for $\overline{\tau} < \tau < \tilde{\tau}$, where $\dfrac{1}{\sigma_0} < x_{T-\tau,\tau}(T-\tau)  \le \frac{1}{\sigma_0(1-e^{-\gamma \tau})}$ and $\eta=T-\tau$},  \\
		\tilde{t} & \mbox{for $\tau \ge \tilde{\tau}$, where $x_{T-\tau,\tau}(T-\tau) > \frac{1}{\sigma_0(1-e^{-\gamma \tau})}$, and $\eta=T-\tilde{t}$}. \\
\end{cases}
\label{tini-opt}
\end{eqnarray}
where $\overline{t} = T-\overline{\tau} \in [0,T-\tau]$ is the unique value,
independent from $\tau\in [0,\overline{\tau}]$, such that
$x_{\overline{t},\tau}(\overline{t})=\frac{1}{\sigma_0}$. Likewise, $\tilde{t}=T-\tilde{\tau} \in [T-\tau,T]$ is the unique value, independent from $\tau\in [\tilde{\tau},T]$, such that $x_{\tilde{t},T-\tilde{t}}(\tilde{t})= \frac{1}{\sigma_0(1-e^{-\gamma (T-\tilde{t})})}$.

\begin{figure}[ht!]
  \begin{center}
    \centerline{\includegraphics[width=7cm]{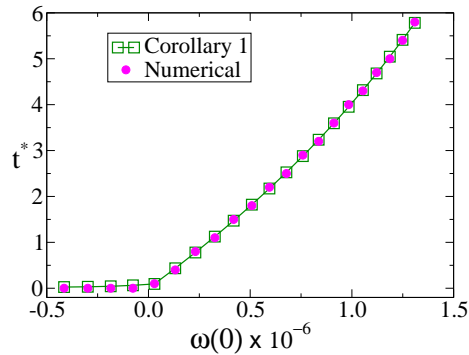}}
    \caption{Optimal initial time of the hard quarantine $t^*$ vs $w(0)$ for $\tau=100$ and $\kappa=0$.  The other parameters are $\gamma=0.01, \sigma_2=\sigma_0=1.5, \sigma_1=0$ and $T=2600$.}
    \label{tini-w0}
  \end{center}
\end{figure}

\begin{figure}[ht!]
  \begin{center}
    \centerline{\includegraphics[width=10cm]{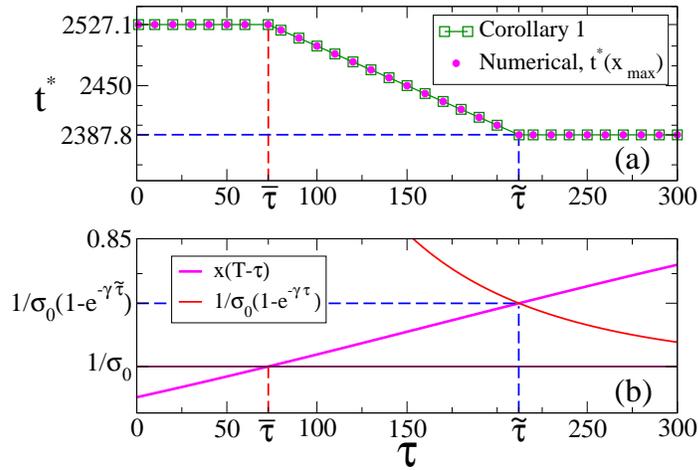}}
    \caption{(a) Optimal initial time $t^*$ vs hard quarantine length
      $\tau$ for $\kappa=0$ and $x_0>1/\sigma_0$.  (b) Graphical
      determination of the times $\overline{\tau}$ and
      $\tilde{\tau}$ that define the three regions for the
      different behaviours of $t^*$.  The parameters are
      $\gamma=0.01, \sigma_2=\sigma_0=1.5, \sigma_1=0$ and
      $T=2600$.  The initial condition corresponds to $x(0)=1-10^{-6},
      y(0)=10^{-6}$ ($R_0=\sigma_0 x(0) >1$).  The optimum times of
      the first and last regions are $\overline{t} \simeq 2527.1$ and
      $\tilde{t} \simeq 2387.8$, respectively, determined by
      $\overline{\tau} \simeq 72.9$ and $\tilde{\tau} \simeq 212.2$.}
    \label{tini-R}
  \end{center}
\end{figure}

\begin{figure}[ht!]
  \begin{center}
    \includegraphics[width=8.3cm]{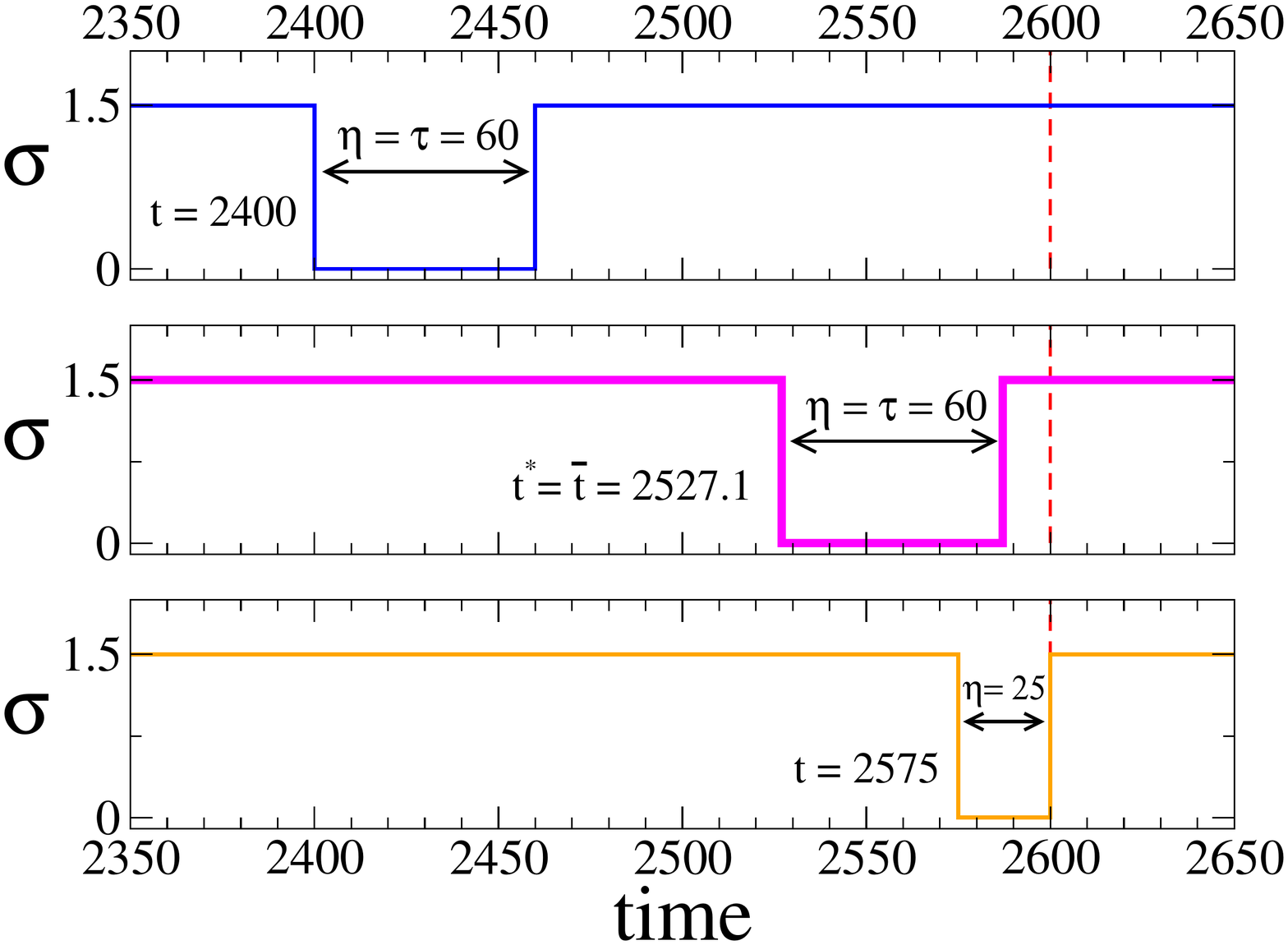}
    \includegraphics[width=8.3cm]{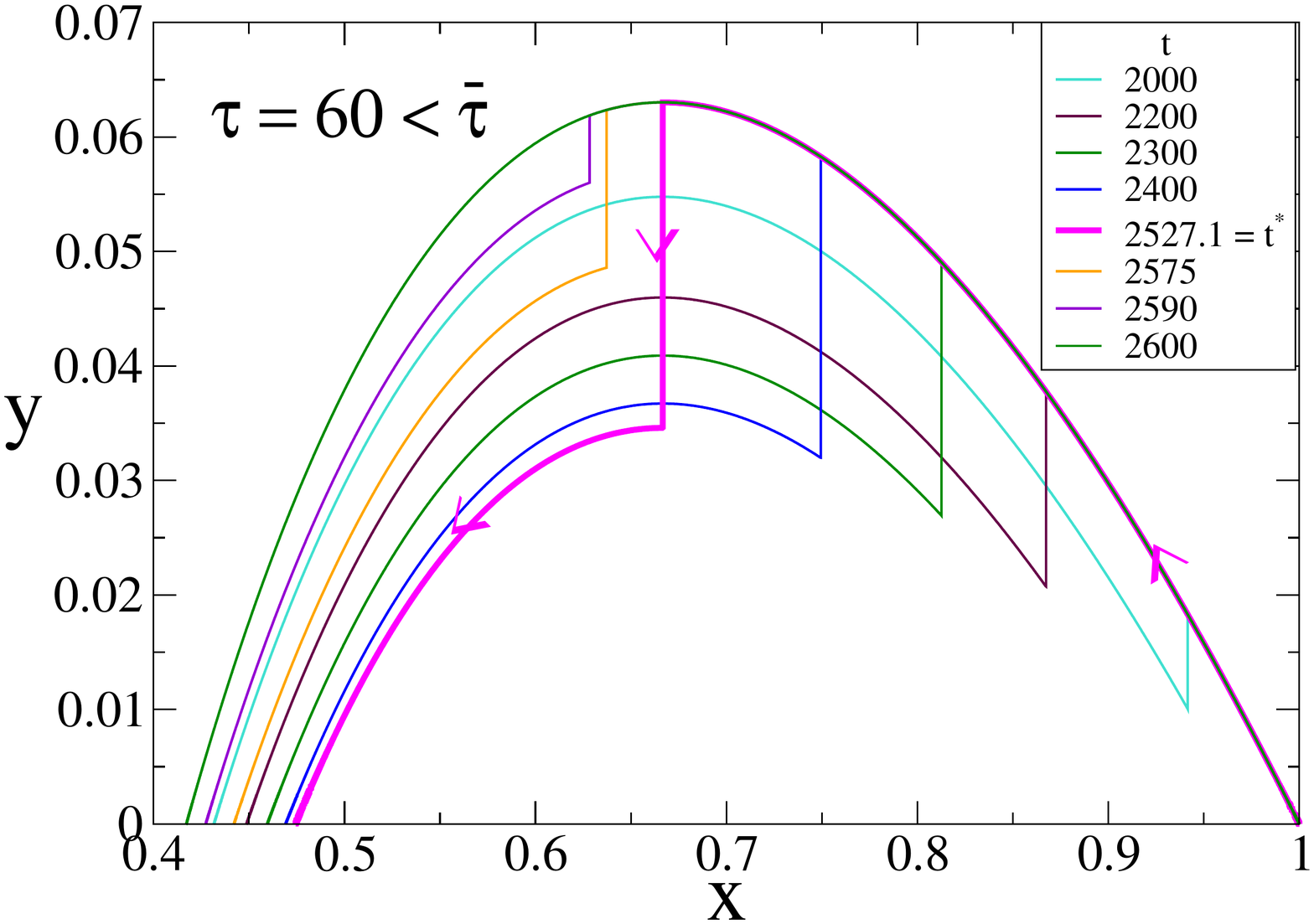}
    \includegraphics[width=8.3cm]{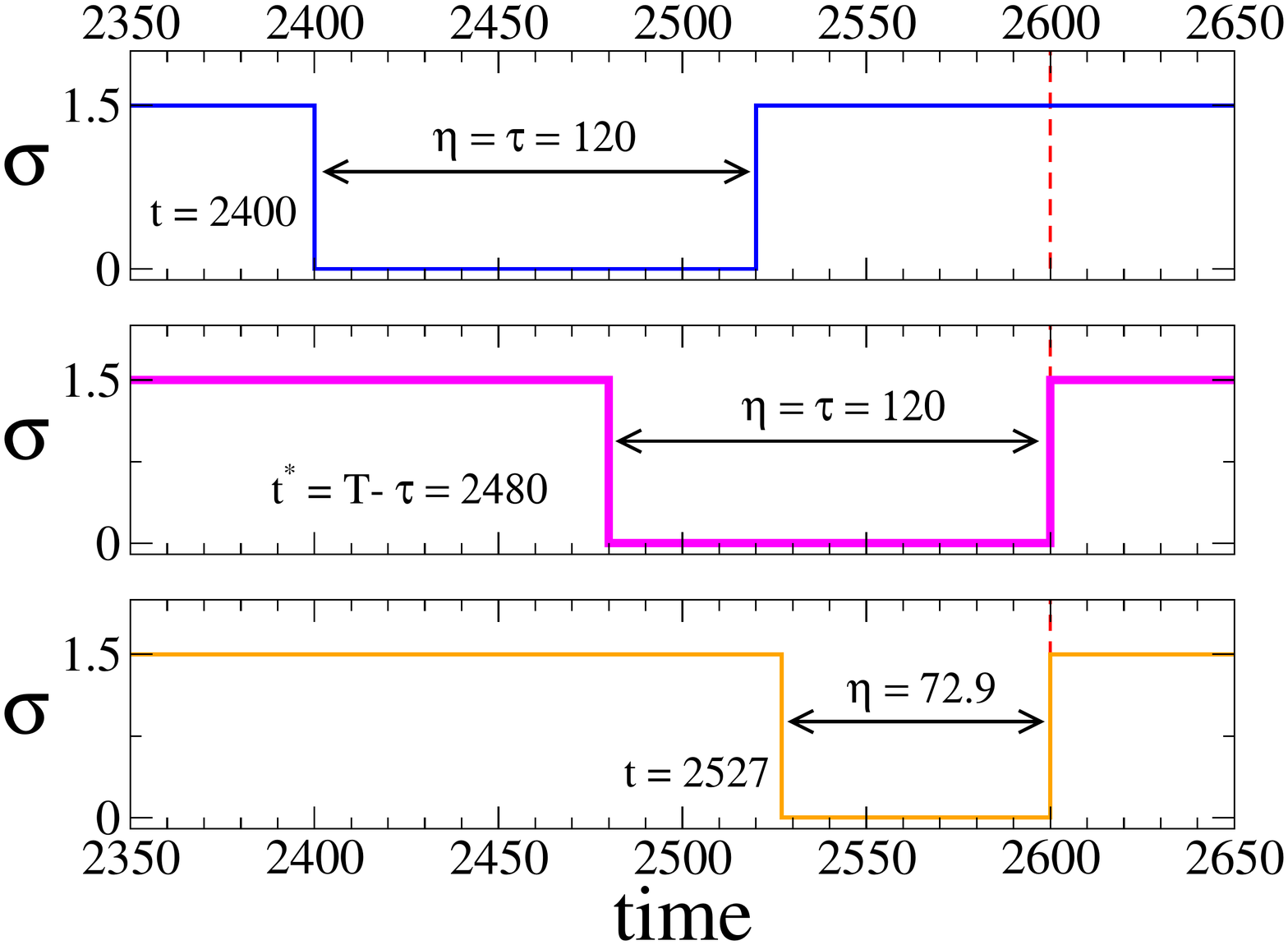}
    \includegraphics[width=8.3cm]{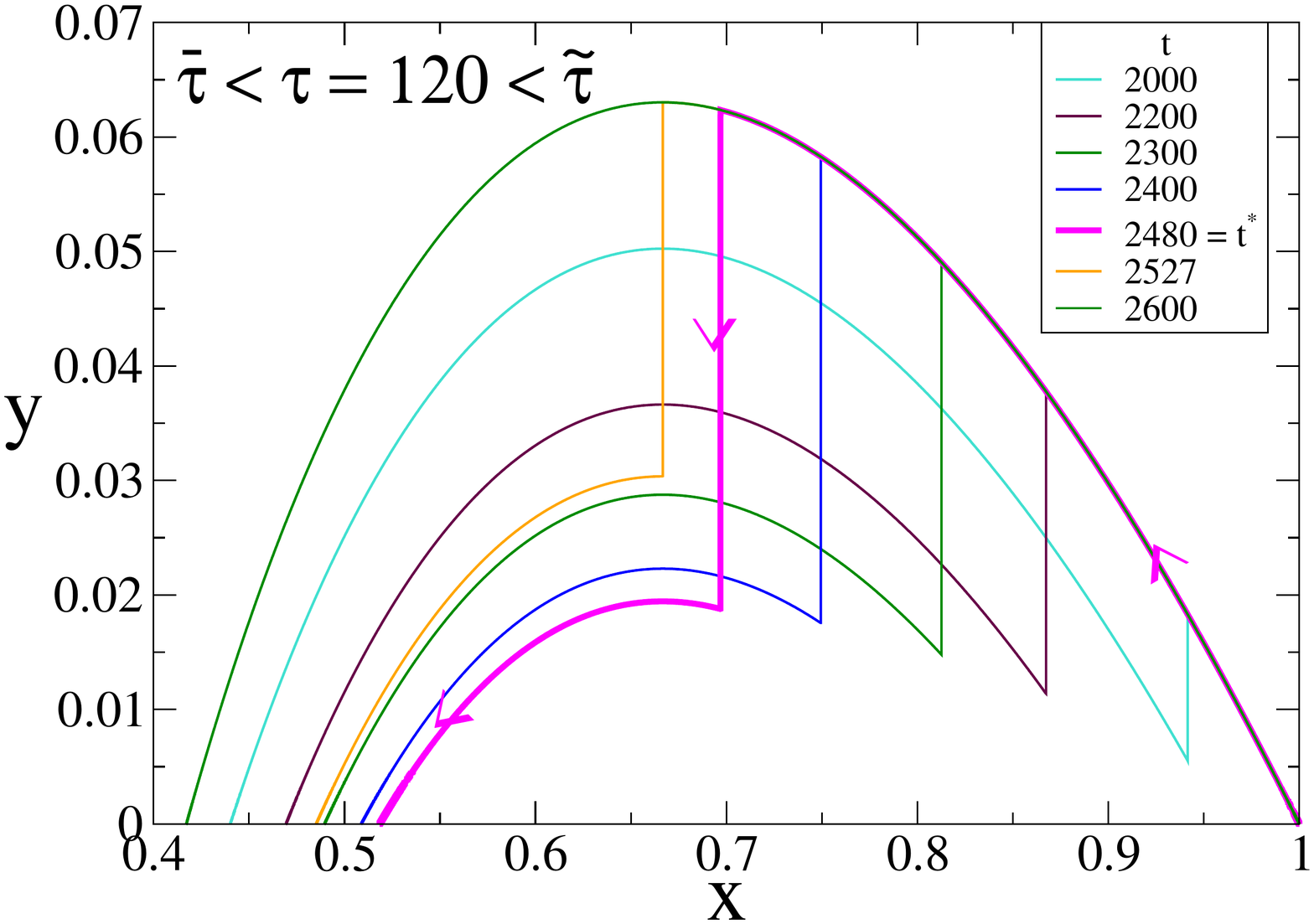}
    \includegraphics[width=8.3cm]{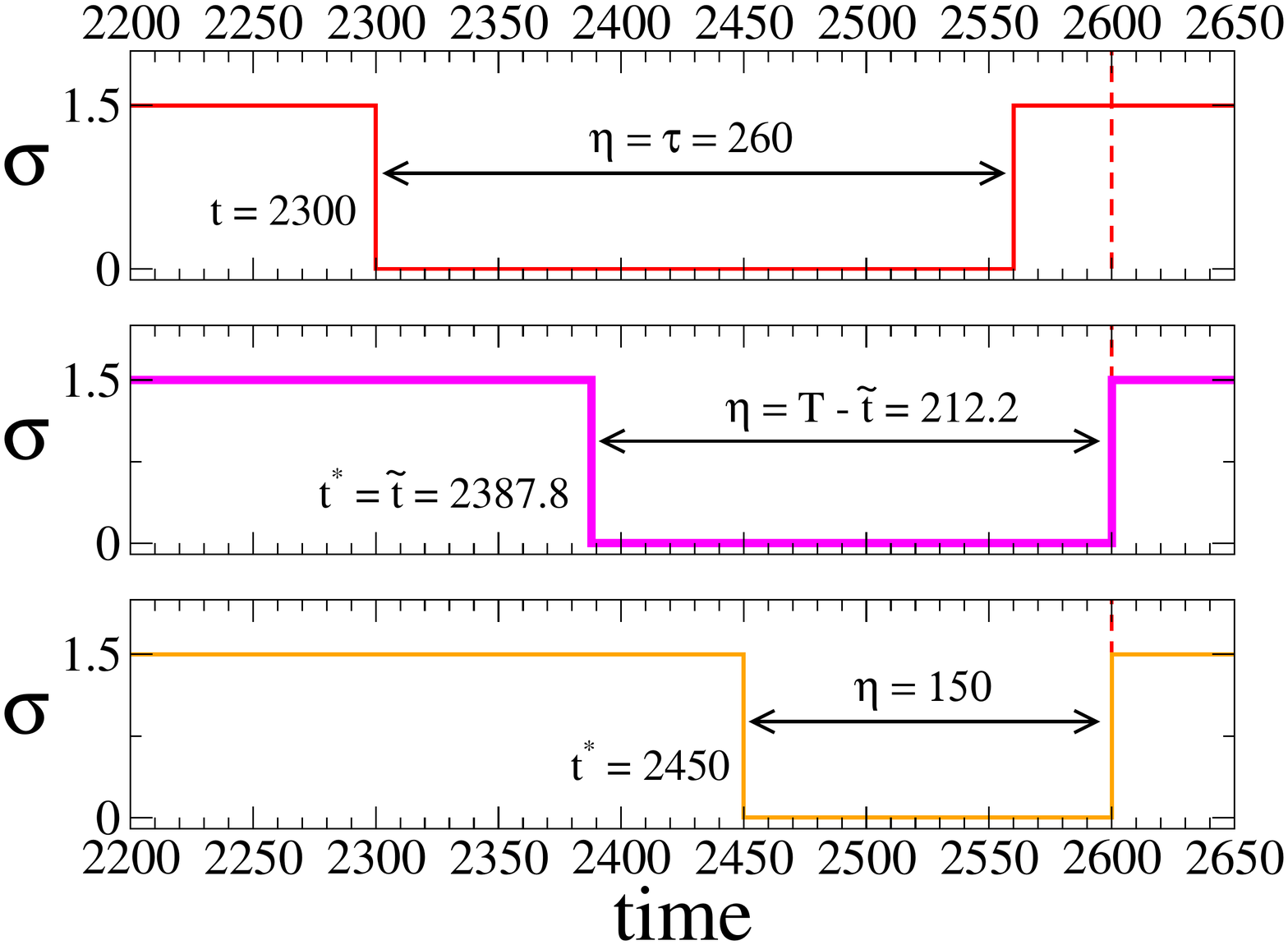}
    \includegraphics[width=8.3cm]{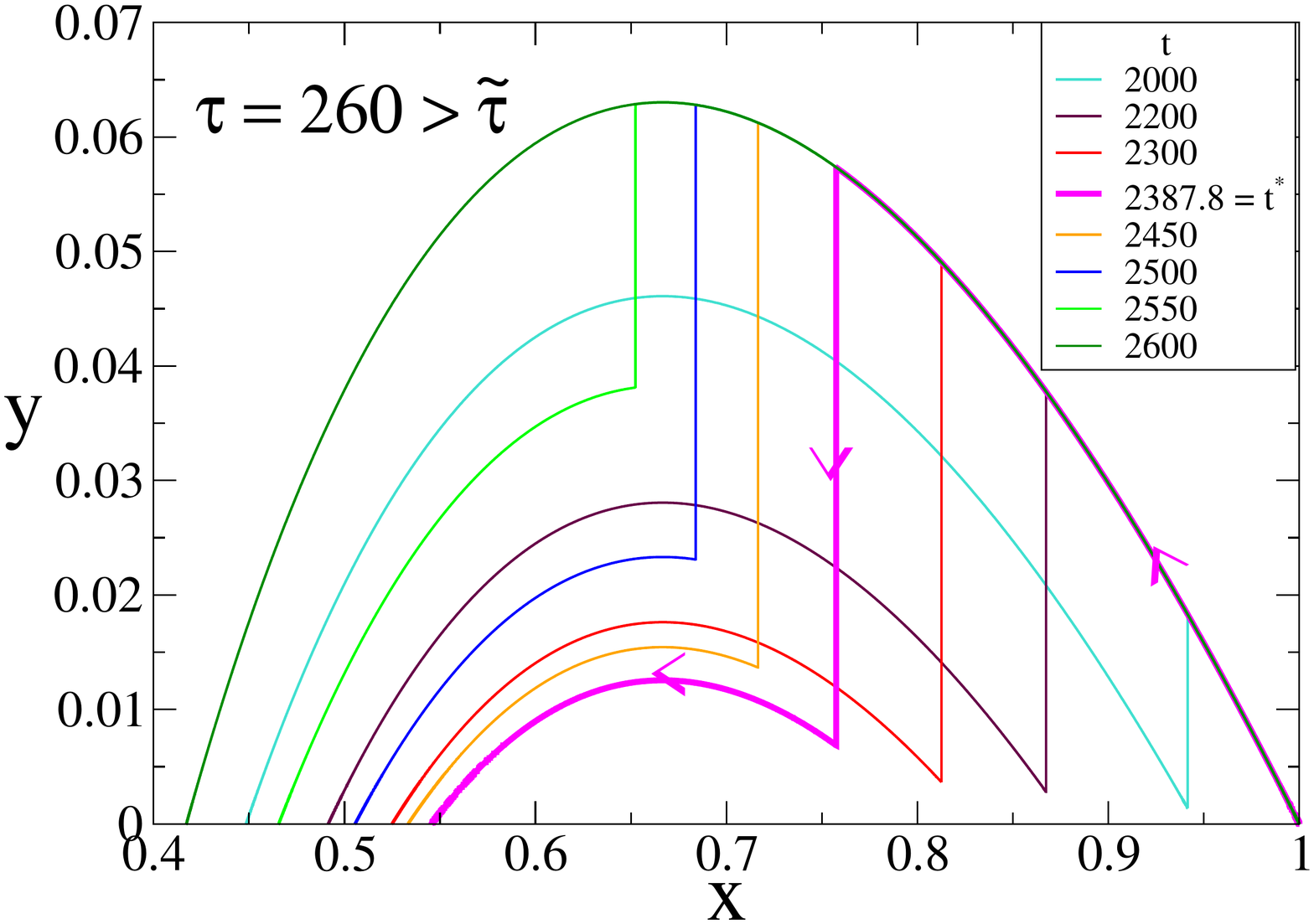}
    \caption{System's trajectory in the $x-y$ phase space (right
      panels), for $\gamma=0.01$, $\sigma_2=\sigma_0=1.5$,
      $\sigma_1=0$ and $T=2600$, and the three values of $\tau$ indicated
      in the legends corresponding to the different regimes of the
      optimum time $t^*$ (pink lines).  Left panels show the time
      evolution of $\sigma$ for three different initial times $t$ of
      the hard quarantine in each case. The optimum times are
      $t^*=\overline{t} \simeq 2527.1$ for $\tau=60$ (top panels),
      $t^*=T-\tau=2480$ for $\tau=120$ (middle panels) and $t^*=\tilde{t} \simeq
      2387.8$ for $\tau=260$ (bottom panels).}
    \label{x-y-2600-00}
  \end{center}
\end{figure}

As $\kappa=0$, $x_{\infty}$ reaches a maximum value when the hard quarantine starts at the optimal time $t^*$ [see \eqref{eq: funcional_J_Llineal}].
In Fig.~\ref{tini-w0} we plot $t^*$ vs $w(0)$ for $\tau=100$, calculated from \eqref{eq: w_0} (squares) and by estimating the maximum of $x_{\infty}$ (circles).  We can see that $t^*$ takes values close to zero for $w(0) \le 0$.  In the rest of this section we consider the case $w(0)>0$.

The behaviour of $t^*$ from \eqref{tini-opt} for $x_0>1/\sigma_0$ ($w(0)>0$) is tested in
Fig.~\ref{tini-R}(a), where we compare numerical results (circles) with
that obtained from \eqref{tini-opt} (squares, Corollary \ref{co:
  sigma1es0}).  We observe that the agreement between numerics and
Corollary \ref{co: sigma1es0} is very good. Figure~\ref{tini-R}(b) is
an auxiliary plot that shows how to obtain graphically the optimum
times $\overline{t} \simeq 2527.1$ and $\tilde{t} \simeq 2387.8$ that
define the three different regimes of $t^*$ defined in
\eqref{tini-opt}.  These times are obtained by estimating the values of $\tau$ for which the curve $x_{t,T-t}(T-\tau)$ crosses the lines $1/\sigma_0$ and $1/\left[ \sigma_0(1-e^{-\gamma \tau}) \right]$, which happens at $\overline{\tau} \simeq 72.9$ and $\tilde{\tau} \simeq 212.2$, respectively.

\begin{remark}
  The effective reproductive number $R_t^{\sigma} \equiv \sigma \,
  x_{\sigma}(t)$ represents the mean number of individuals that an
  agent infects during its infectious period, at time $t$.  It is
  interesting to note that the optimal time from \eqref{tini-opt} can be rewriten in terms of $R_t^{\sigma}$ as
  \begin{eqnarray}
    t^* =
    \begin{cases}
      \overline{t} & \mbox{for $0 \le \tau \le \overline{\tau}$, where $R_{T-\tau}^{\sigma_0} \le 1$ and
        $\eta=\tau$}, \\
      T-\tau & \mbox{for $\overline{\tau} \le \tau \le \tilde{\tau}$, where $1 < R_{T-\tau}^{\sigma_0} \le \frac{1}{1-e^{-\gamma \tau}}$ and $\eta=T-\tau$},  \\
      \tilde{t} & \mbox{for $\tau > \tilde{\tau}$, where $R_{T-\tau}^{\sigma_0} > \frac{1}{1-e^{-\gamma \tau}}$, and $\eta=T-\tilde{t}$}. \\
    \end{cases}
    \label{tini-opt2}
  \end{eqnarray}
  Here $R_{T-\tau}^{\sigma_0}=\sigma_0 \, x_{\sigma_0}(T-\tau)$,
  $\overline{\tau}=T-\overline{t}$ and
  $\tilde{\tau}=T-\tilde{t}$, where $\overline{t}$ and $\tilde{t}$ are
  determined from the relations
  \begin{eqnarray}
    R_{\overline{t}}^{\sigma_0} = 1 ~~~~~ \mbox{and} ~~~~~ R_{\tilde{t}}^{\sigma_0} =
    \frac{1}{1-e^{-\gamma(T-\tilde{t})}}.
  \end{eqnarray}
\end{remark}

In Fig.~\ref{x-y-2600-00} we show the evolution of the system in the
$x-y$ phase space for a given $\tau$ and various $t$ (right panels),
together with the evolution of $\sigma(t)$ (left panels), which
describe the three different behaviours of $t^*$.  All curves start at
$(x_0,y_0)=(0.999999,0.000001)$ and follow the top curve with soft
quarantine ($\sigma=\sigma_0$) until the hard quarantine starts at $t$
($\sigma=\sigma_1=0$), vertically falling down up to a lower level
curve when the soft quarantine starts, and finally following this
curve until the fixed point $(x_{\infty},0)$ is asymptotically
reached.  The vertical trajectory describes the evolution within the
hard quarantine where $x(t)$ remains constant, given that
$\sigma(t)=\sigma_1=0$ in that period.  The optimum time $t^*$ that
leads to the maximum of $x_{\infty}$ corresponds to the time for which
$y(t)$ drops to the lowest level curve in the interval $[t,t+\eta]$
(pink curve).  For $\tau=60 < 72.9 = \overline{\tau}$
(Fig.~\ref{x-y-2600-00} top panels) we
see that the maximum of $x_{\infty}$ is reached starting the hard
quarantine at $t^*=\overline{t}=2527.1$, where the effective
reproduction number is
$R_{\overline{t}}^{\sigma_0}=R_{\overline{t}+\eta}^{\sigma_0}=1$, and
thus there is no new outbreak when the hard quarantine is released
($\frac{dy}{dt} |_{t+\eta}=0$).  In this case the entire quarantine
period $\eta=\tau$ is used.  For $\overline{\tau} < \tau=120 <
\tilde{\tau}=212.2$ (Fig.~\ref{x-y-2600-00} middle panels) the optimum initial time is $t^*=T-\tau=2480 < \overline{t}$, obtained by still using the entire hard quarantine period but starting earlier than $\overline{t}$.  Finally, for
$\tau=260 > \tilde{\tau}$ (Fig.~\ref{x-y-2600-00} bottom panels) the
optimum is $t^*=\tilde{t}=2387.8>T-\tau=2340$, where it turns more
effective to use the hard quarantine for a shorter time $T-t^*<
\tau$.  Notice that implementing a shorter but later hard quarantine
is more efficient than using a longer and earlier hard quarantine, as
we can see by comparing $\sigma(t)$ for $\eta=260$ and $\eta=212.2$ in
the bottom panels for the $\tau=260$ case.

\subsection{Case $\kappa=0$, $\sigma_2 = \sigma_0$ and
  $\sigma_1>0$ (theorem  \ref{te: control_optimo_sigma1gral})}

We now analyze the case $\kappa=0$, $\sigma_2=\sigma_0=1.5$ and $\sigma_1=0.3>0$,
with $\gamma=0.01$ and $T=2600$.  This corresponds to a hard
quarantine that is softer than in the previous case $\sigma_1=0$, and
during which there are infections. Initially is $y_0=y(t=0)=10^{-6}$
and $x_0=x(t=0)=1-10^{-6}$.  We can
see from theorem \ref{te: control_optimo_sigma1gral} that the optimum initial time of the hard quarantine $t^*$ for $w(0)>0$ is
\begin{eqnarray}
t^* =
\begin{cases}
		\overline{t} & \mbox{for $0 \le \tau \le \overline{\tau}$, where $w(0)\ge 0, w(T-\tau) \le 0$ and
			$\eta=\tau$}, \\
		T-\tau & \mbox{for $\overline{\tau} \le \tau \le \tilde{\tau}$, where $0 < w(T-\tau)  \le \frac{1}{\gamma y_{T-\tau,\tau}(T-\tau)}$ and $\eta=T-\tau$},  \\
		\tilde{t} & \mbox{for $\tau > \tilde{\tau}$, where $w(T-\tau) >\frac{1}{\gamma y_{T-\tau,\tau}(T-\tau)}$, and $\eta=T-\tilde{t}$}, \\
\end{cases}
\label{tini-opt1}
\end{eqnarray}
where $\overline{t}\in [0,T-\tau]$ is the unique value, depending on
$\tau\in [0,\overline{\tau}]$, such that $w(\overline{t})=0$. On the
other hand, $\tilde{t}\in [T-\tau,T]$ is the unique value, independent
from $\tau\in [\tilde{\tau},T]$, such that $w(\tilde{t})=
\frac{1}{\gamma y_{\tilde{t},T-\tilde{t}}(\tilde{t})}$. The dependence
and  independence of $w(t)$ on $\tau$ for $t\in[0,T-\tau]$ and
$t\in[T-\tau,T]$, respectively, can be seen from the definition of
$w(t)$ in \eqref{eq: w_2_0}.

\begin{figure}[ht!]
  \begin{center}    \centerline{\includegraphics[width=10cm]{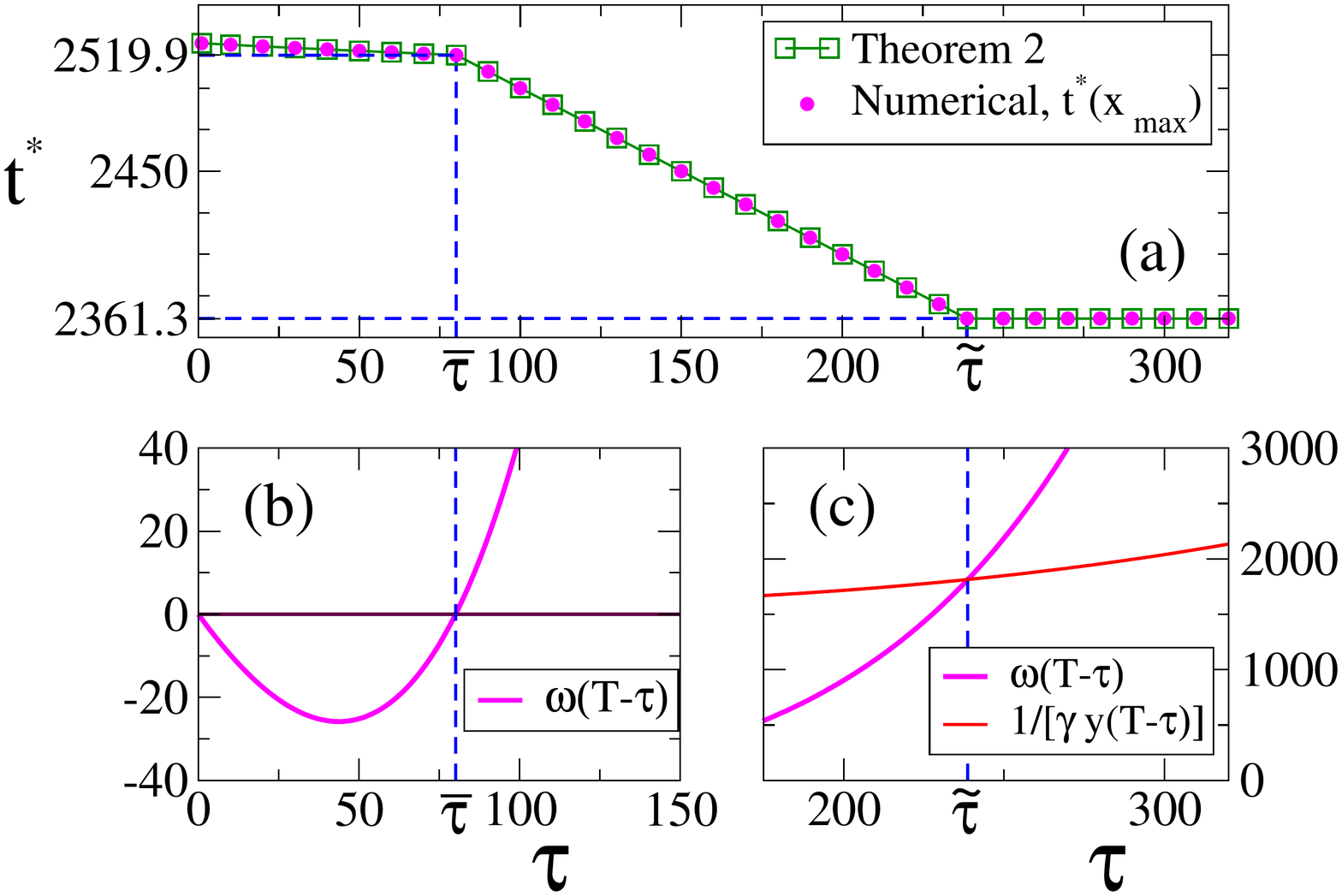}}
    \caption{(a) Optimal initial time $t^*$ vs hard quarantine length
      $\tau$ for $\kappa=0$ and $w(0)>0$.  (b) and (c) Graphical determination of the times
      $\overline{\tau}$ and $\tilde{\tau}$, respectively, which define
      the three regions for the different behaviours of $t^*$.  The
      parameters are $\gamma=0.01, \sigma_0=1.5, \sigma_1=0.3,
      T=2600$.  The initial condition corresponds to $x(0)=1-10^{-6},
      y(0)=10^{-6}$. The optimum time for the region $\tau >
      \tilde{\tau} \simeq 238.7$ is $\tilde{t} \simeq 2361.3$, while $t^*$ has a
      slight dependence on $\tau$ for $\tau \le \overline{\tau} \simeq
      80.1$.}
    \label{tini-R2}
  \end{center}
\end{figure}

In Fig.~\ref{tini-R2}(a) we compare numerical results (circles) with
results from \eqref{tini-opt1} (squares, theorem \ref{te:
  control_optimo_sigma1gral}), where we see a very good agreement.  At $t^*$, $x_{\infty}$ reaches a maximum.
Unlike the $\sigma_1=0$ case, for $\sigma_1=0.3>0$ the optimal time
$t^*$ in the $0 \le \tau \le \overline{\tau} \simeq 80.1$ interval
depends on $\tau$, that is, $t^*=\overline{t}(\tau)$, while for $\tau
> \tilde{\tau} \simeq 238.7$ is $t^*=\tilde{t} \simeq 2361.3$ independent of $\tau$.
Figures~\ref{tini-R2}(b) and (c) show that the optimal times
$\overline{t}$ and $\tilde{t}$ are estimated, respectively, as the values of
$t=T-\tau$ for which the curve $w(T-\tau)$ crosses the horizontal line $0$
and the curve $1/\left[\gamma y_{t,T-t}(T-\tau) \right]$.

\begin{figure}[ht!]
  \begin{center}
    \includegraphics[width=8.3cm]{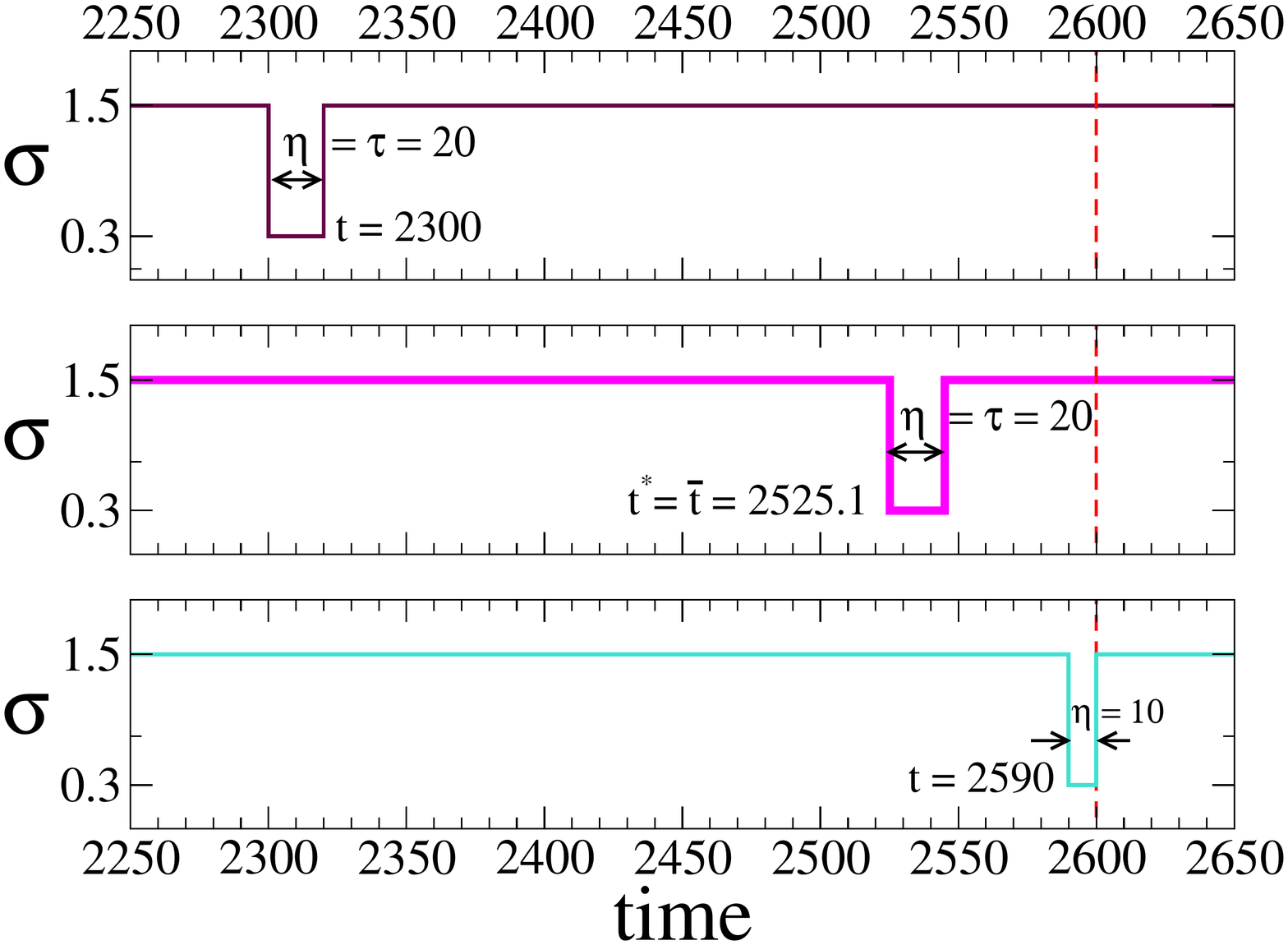}
    \includegraphics[width=8.3cm]{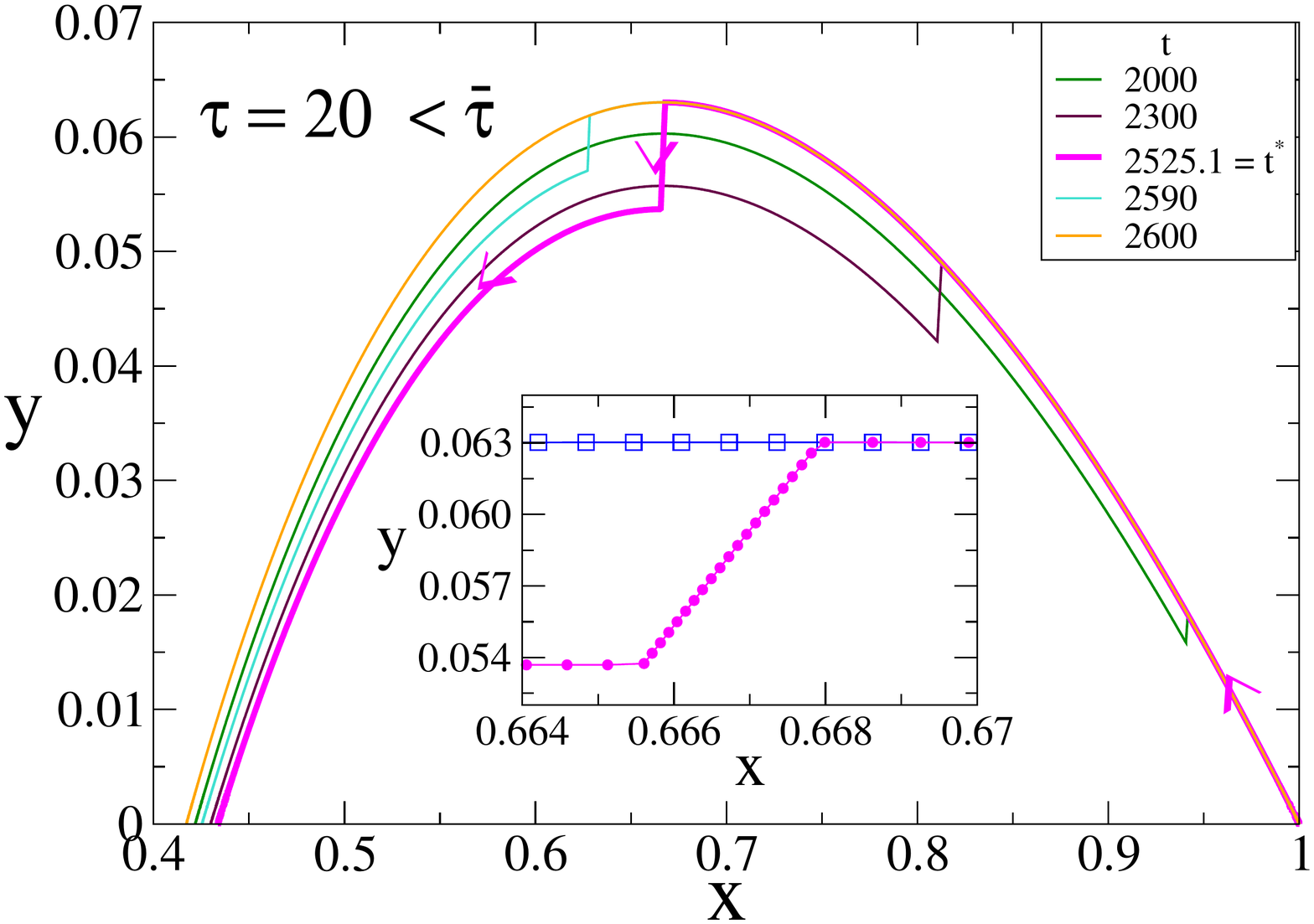}
    \includegraphics[width=8.3cm]{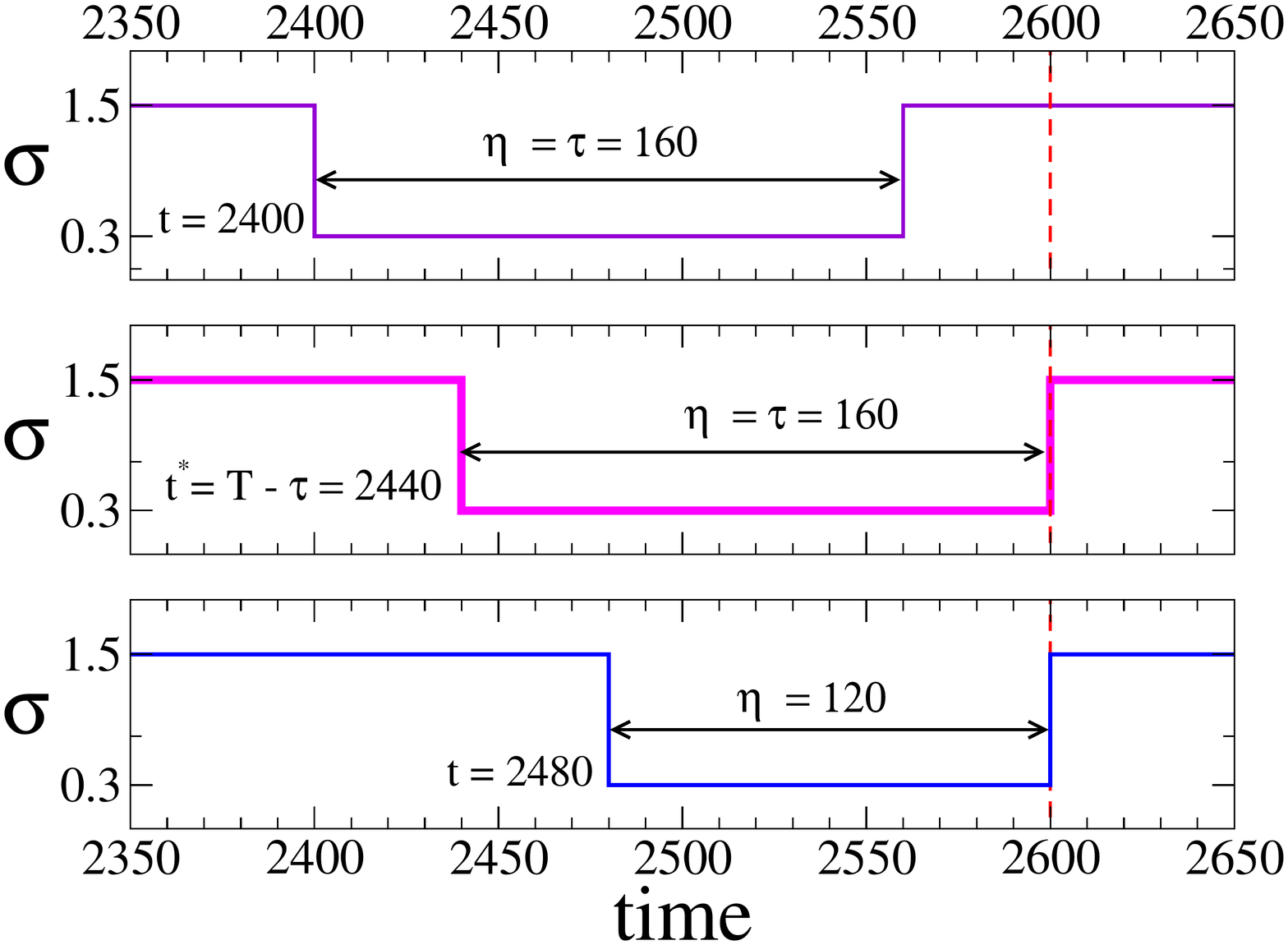}
    \includegraphics[width=8.3cm]{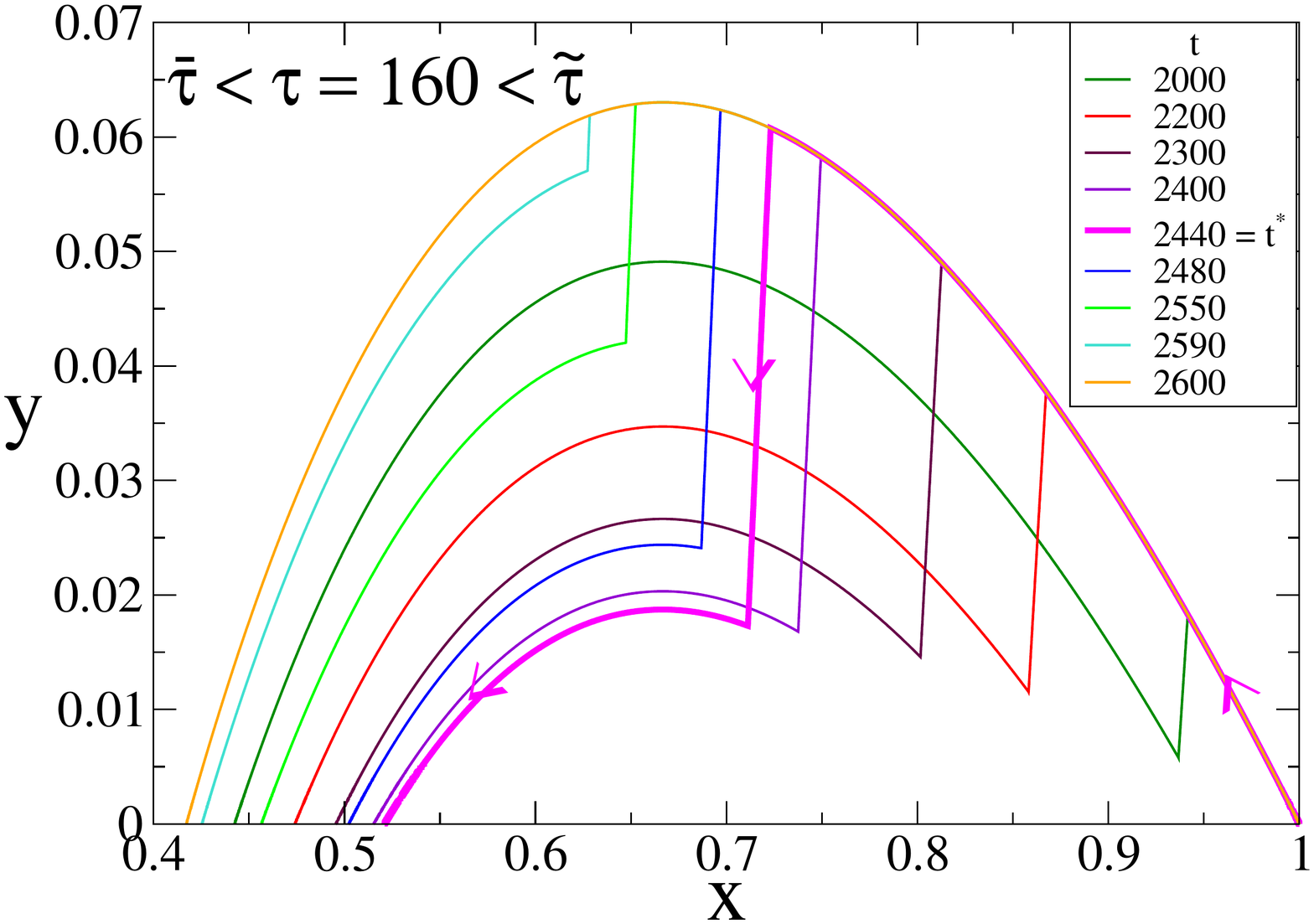}
    \includegraphics[width=8.3cm]{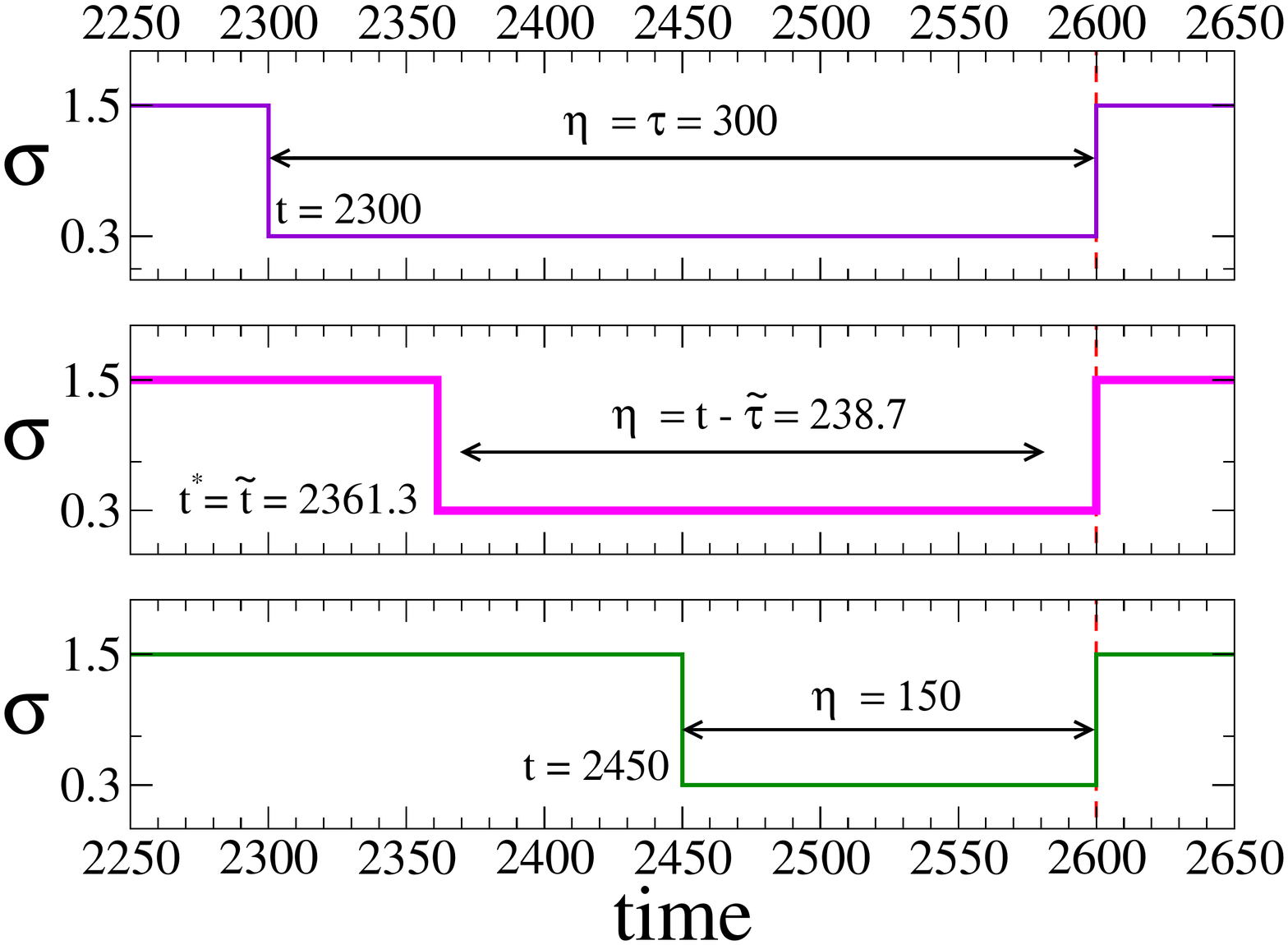}
    \includegraphics[width=8.3cm]{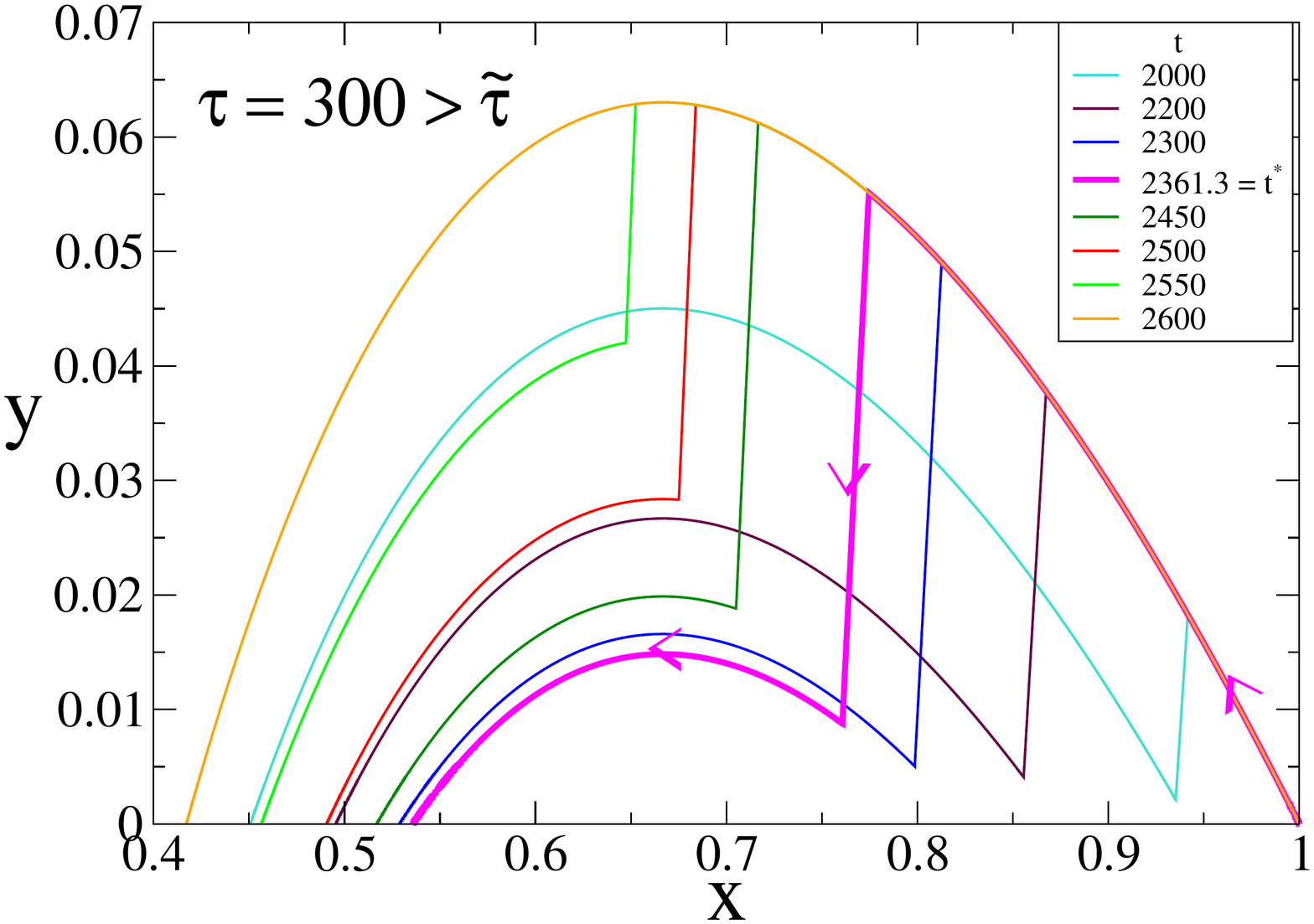}
    \caption{System's trajectory in the $x-y$ phase space (right
      panels), for $\gamma=0.01$, $\sigma_2=\sigma_0=1.5$,
      $\sigma_1=0.3$ and $T=2600$, and the three values of $\tau$ indicated
      in the legends corresponding to the different regimes of the
      optimum time $t^*$ (pink lines).  Left panels show the time
      evolution of $\sigma$ for three different initial times $t$ of
      the hard quarantine in each case. The optimum times are
      $t^*=\overline{t} \simeq 2525.1$ for $\tau=20$ (top panels),
      $t^*=T-\tau=2440$ for $\tau=160$ (middle panels) and $t^*=\tilde{t} \simeq
      2361.3$ for $\tau=300$ (bottom panels).}
    \label{x-y-2600-03}
  \end{center}
\end{figure}

Figure~\ref{x-y-2600-03} is analogous to Fig.~\ref{x-y-2600-00} for
the $\sigma_1=0$ case, and depicts the three different behaviours of
$t^*$.  Curves are similar to those of $\sigma_1=0$, where the main
difference is that for $\sigma_1=0.3>0$ the trajectory of the system
within the hard quarantine in the $x-y$ space is described by a
diagonal line (see inset of top-right panel), given that $\sigma_1$ is
larger than zero and thus $x(t)$ decreases in this period.  At the
optimum time $t^*$, $y(t)$ drops to the lowest level curve in the
interval $[t,t+\eta]$ (pink curves).

\subsection{General case $\kappa>0$ and $0 < \sigma_1 < \sigma_2 <  \sigma_0$ (theorem 3)}

In this section we analyze the most general case $\kappa=10^{-5}>0$, with a soft quarantine ($\sigma_2=1.5$) together with a hard quarantine that is not very strict ($\sigma_1=0.3$) during the control interval $t \in [0,T]$, and a case scenario that simulates no restrictions or quarantine ($\sigma(t)=\sigma_0=2.2$) after the control period $t >T$.  The rest of the parameters are the same as those in the previous studied cases.  Then, from theorem \ref{th: version_general}  the optimum initial time $t^*$ is given by
\begin{eqnarray}
t^* =
\begin{cases}
		\overline{t} & \mbox{for $0 \le \tau \le \overline{\tau}$, where $w(0)\ge 0, w(T-\tau) \le 0$ and  $\eta=\tau$}, \\
		T-\tau & \mbox{for $\overline{\tau} \le \tau \le \tilde{\tau}$, where $0 < w(T-\tau)  \le \alpha(T-\tau)$ and $\eta=T-\tau$},  \\
		\tilde{t} & \mbox{for $\tau > \tilde{\tau}$, where $w(T-\tau) >\alpha(T-\tau)$, and $\eta=T-\tilde{t}$}, \\
\end{cases}
\label{tini-opt3}
\end{eqnarray}
where $\overline{t}\in [0,T-\tau]$ is a unique value that depends on
$\tau\in [0,\overline{\tau}]$ and satisfies $w(\overline{t})=0$, while $\tilde{t}\in [T-\tau,T]$ is a unique value independent of $\tau\in [\tilde{\tau},T]$ that satisfies $w(\tilde{t})= \alpha(\tilde{t})$.  Here $\alpha(t)$ is given by \eqref{eq: alpha}, whereas the dependence and independence of $w(t)$ on $\tau$ for $t\in[0,T-\tau]$ and $t\in[T-\tau,T]$, respectively, is seen in the definition of
$w(t)$ in \eqref{eq: omega}.

\begin{figure}[t]
  \begin{center}    \centerline{\includegraphics[width=10cm]{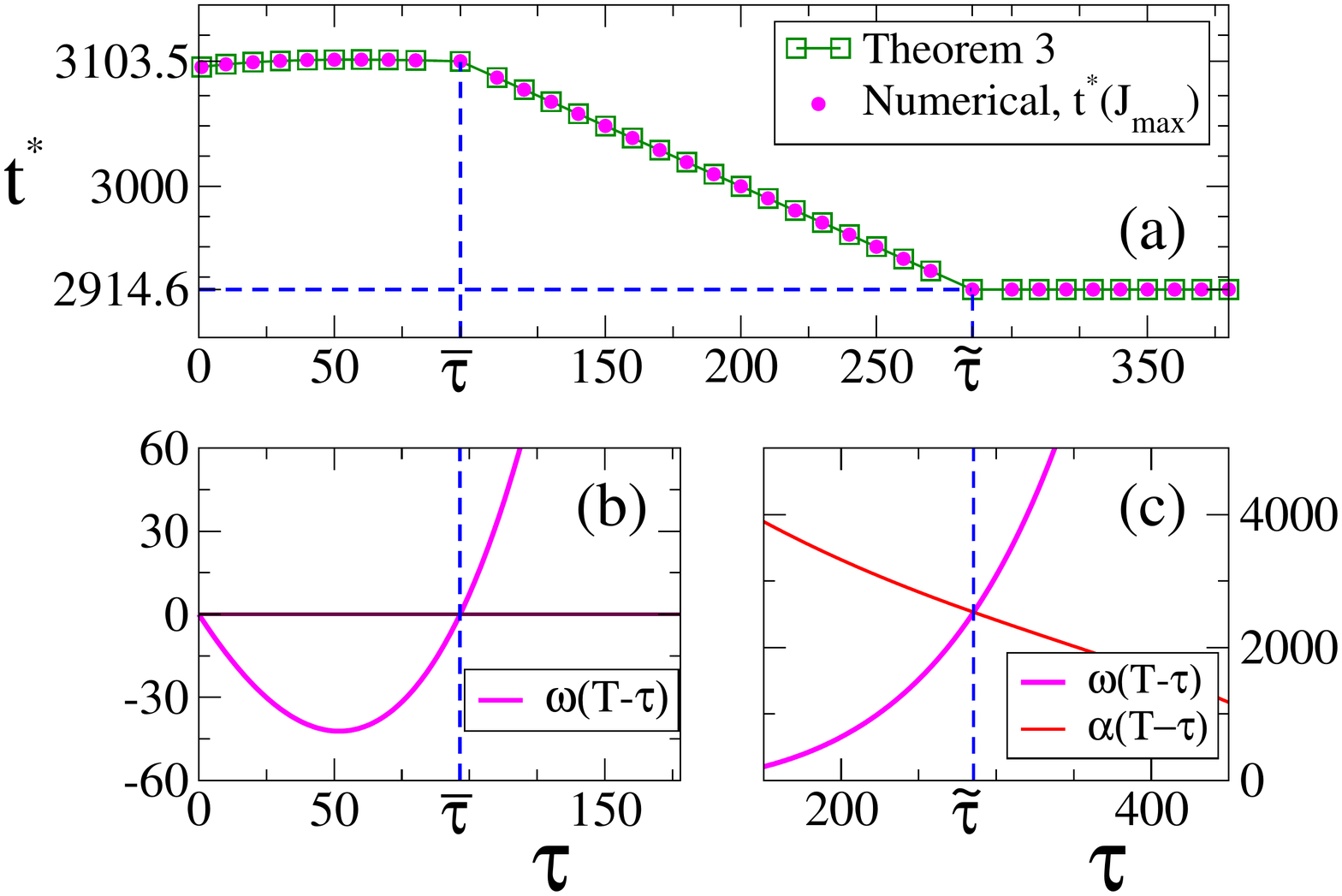}}
    \caption{(a) Optimal initial time $t^*$ vs hard quarantine length
      $\tau$ for $\kappa=10^{-5}$ and $w(0)>0$.  (b) and (c) Graphical
      determination of the times $\overline{\tau}$ and $\tilde{\tau}$,
      respectively, which define the three regions for the different
      behaviours of $t^*$.  The parameters are $\gamma=0.01$,
      $\sigma_0=2.2$, $\sigma_1=0.3$, $\sigma_2=1.5$ and $T=3200$.
      The initial condition corresponds to $x(0)=1-10^{-6},
      y(0)=10^{-6}$. The optimum time for the region
      $\tau > \tilde{\tau} \simeq 285.4$ is $\tilde{t} \simeq 2914.6$,
      while $t^*$ has a
      slight dependence on $\tau$ for $\tau \le \overline{\tau} \simeq
      96.5$.}
    \label{tini-R3}
  \end{center}
\end{figure}

Given that we consider here $\kappa>0$, $J$ reaches a maximum at the optimum time $t^*$ (see \eqref{eq: funcional_J_Llineal}).  Figure~\ref{tini-R3}(a) shows the behaviour of $t^*$ as a function of $\tau$ for $\kappa=10^{-5}$, where we observe a very good agreement between numerical results (circles) and theorem \ref{th: version_general} (squares).  We also see that $t^*$ depends slightly on $\tau$ in the $0 \le \tau \le \overline{\tau} \simeq 96.5$ interval, while $t^*=\tilde{t} \simeq 2914.6$ for $\tau > \tilde{\tau} \simeq 285.4$.  The optimal times $\overline{t}$ and $\tilde{t}$ are estimated as the values of $t=T-\tau$ for which the curve $w(T-\tau)$ crosses the horizontal line $0$ and the curve $\alpha(T-\tau)$, respectively (Figs.~\ref{tini-R3}(b) and (c)).

\begin{figure}[ht!]
  \begin{center}
    \includegraphics[width=8.3cm]{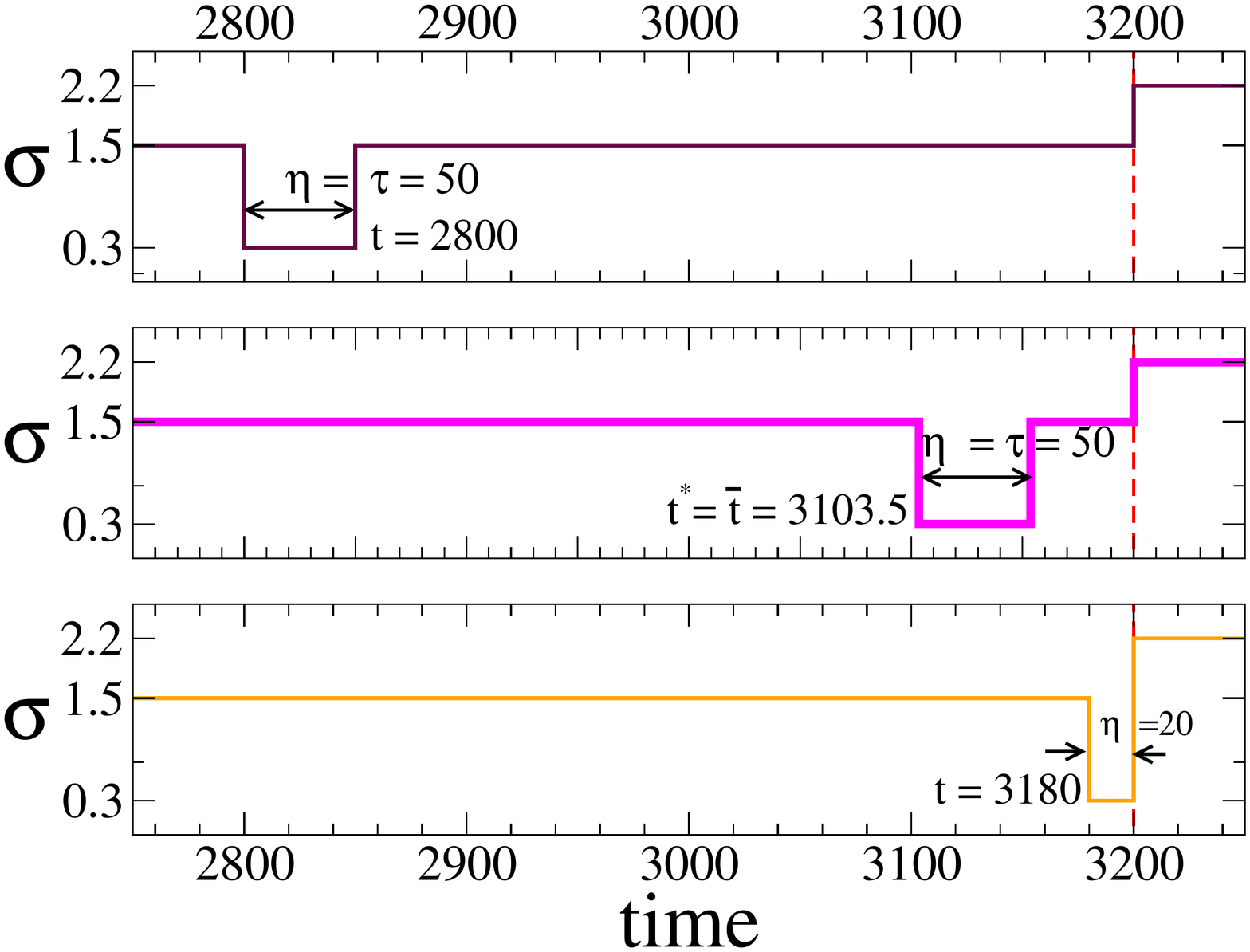}
    \includegraphics[width=8.3cm]{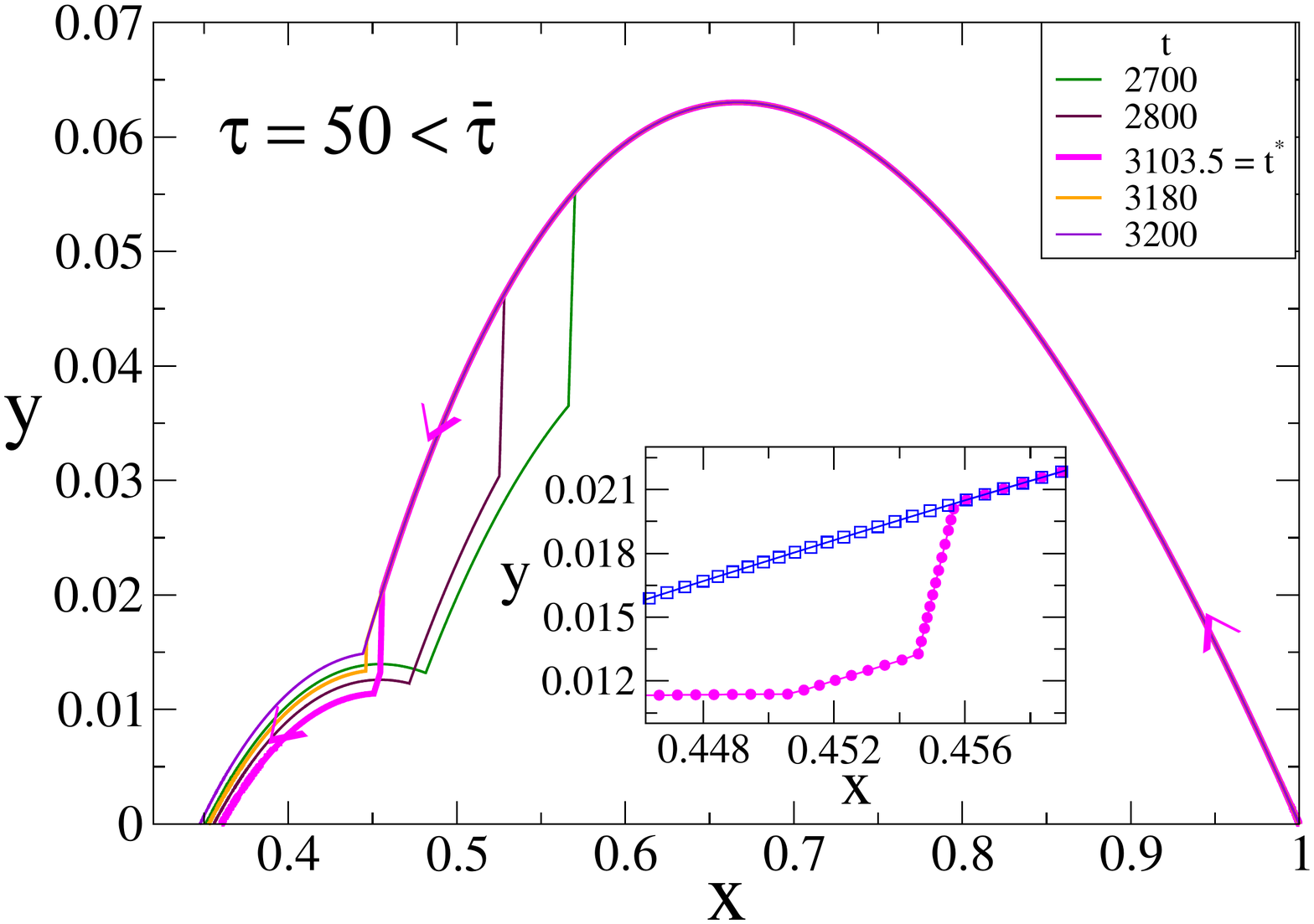}
    \includegraphics[width=8.3cm]{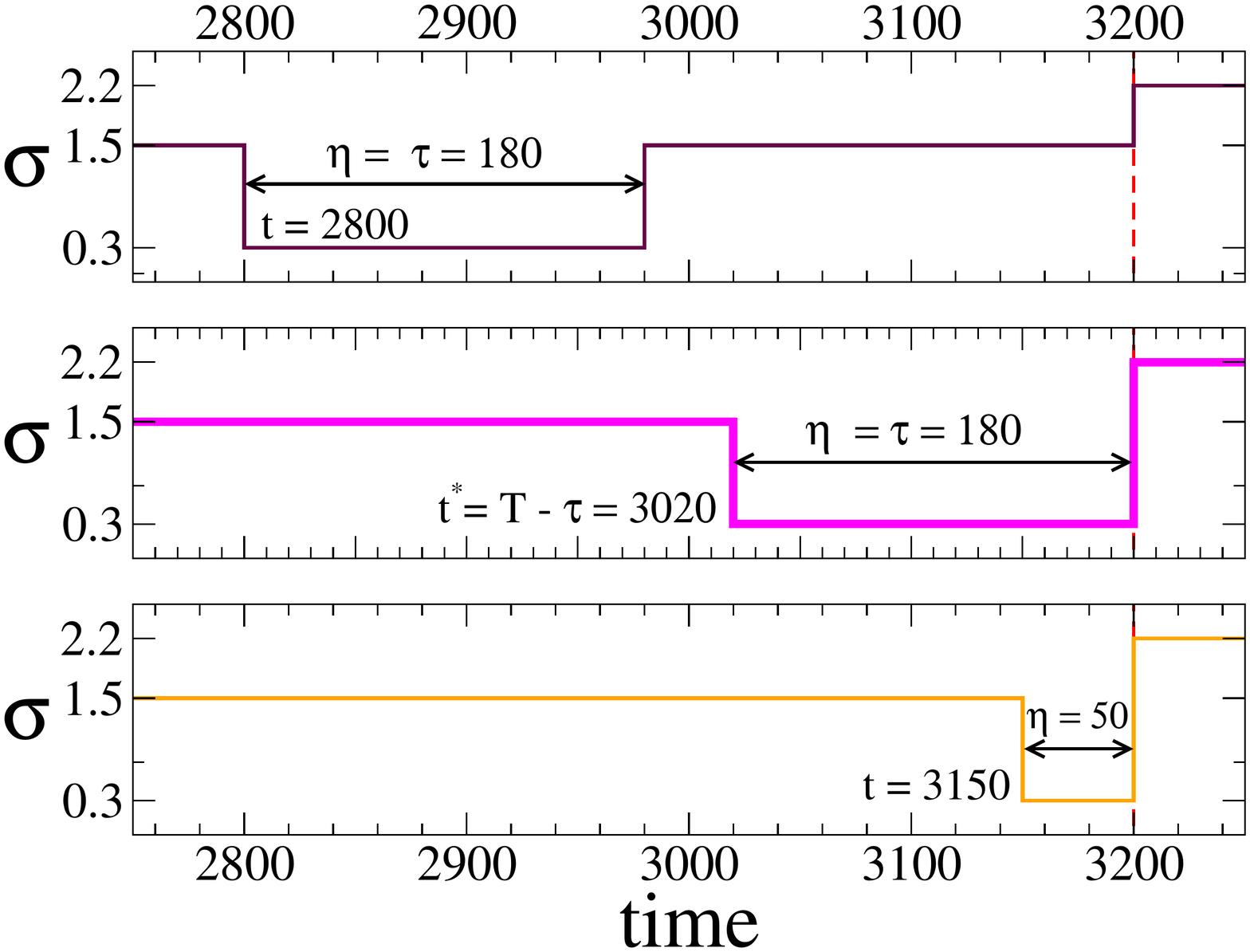}
    \includegraphics[width=8.3cm]{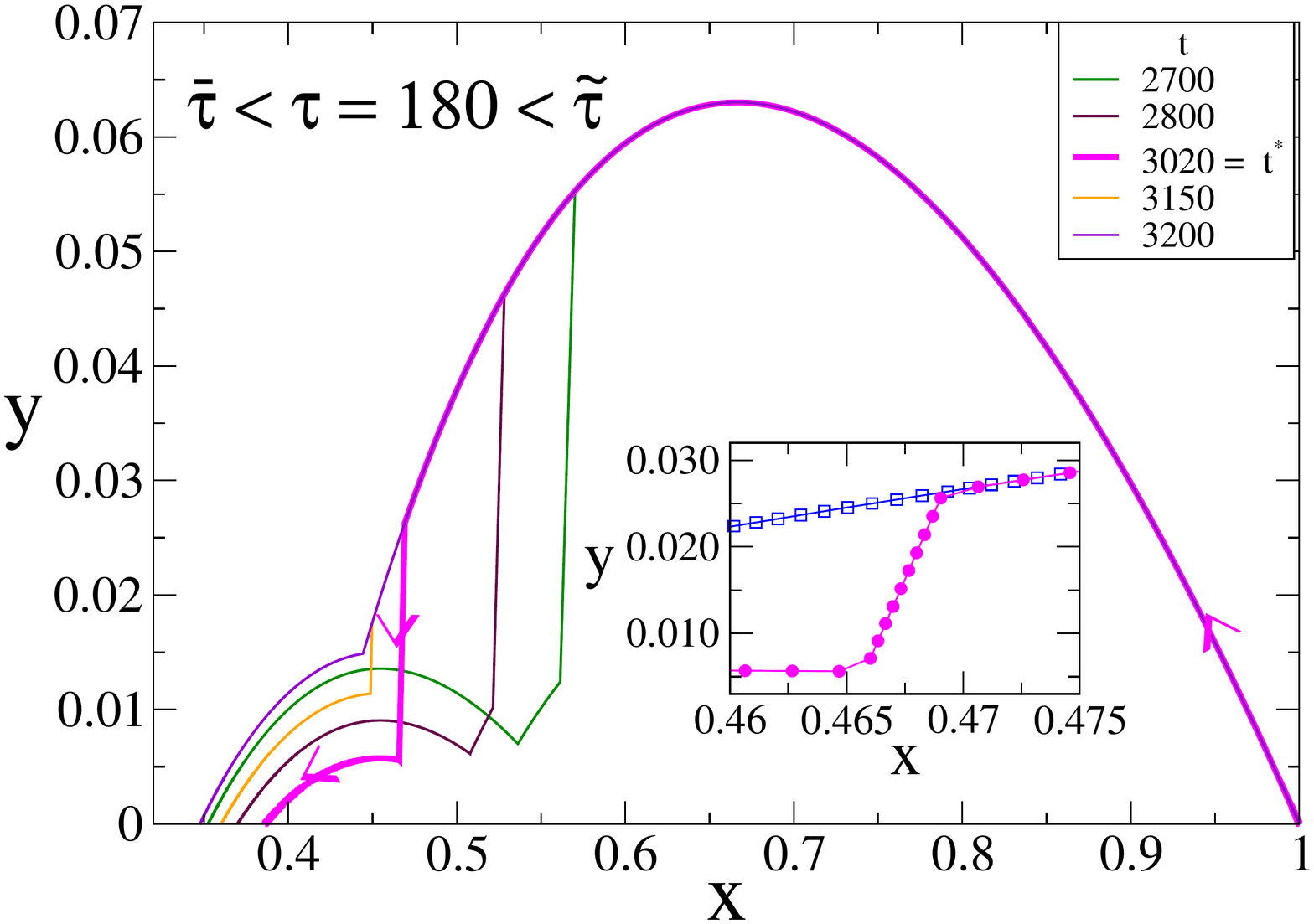}
    \includegraphics[width=8.3cm]{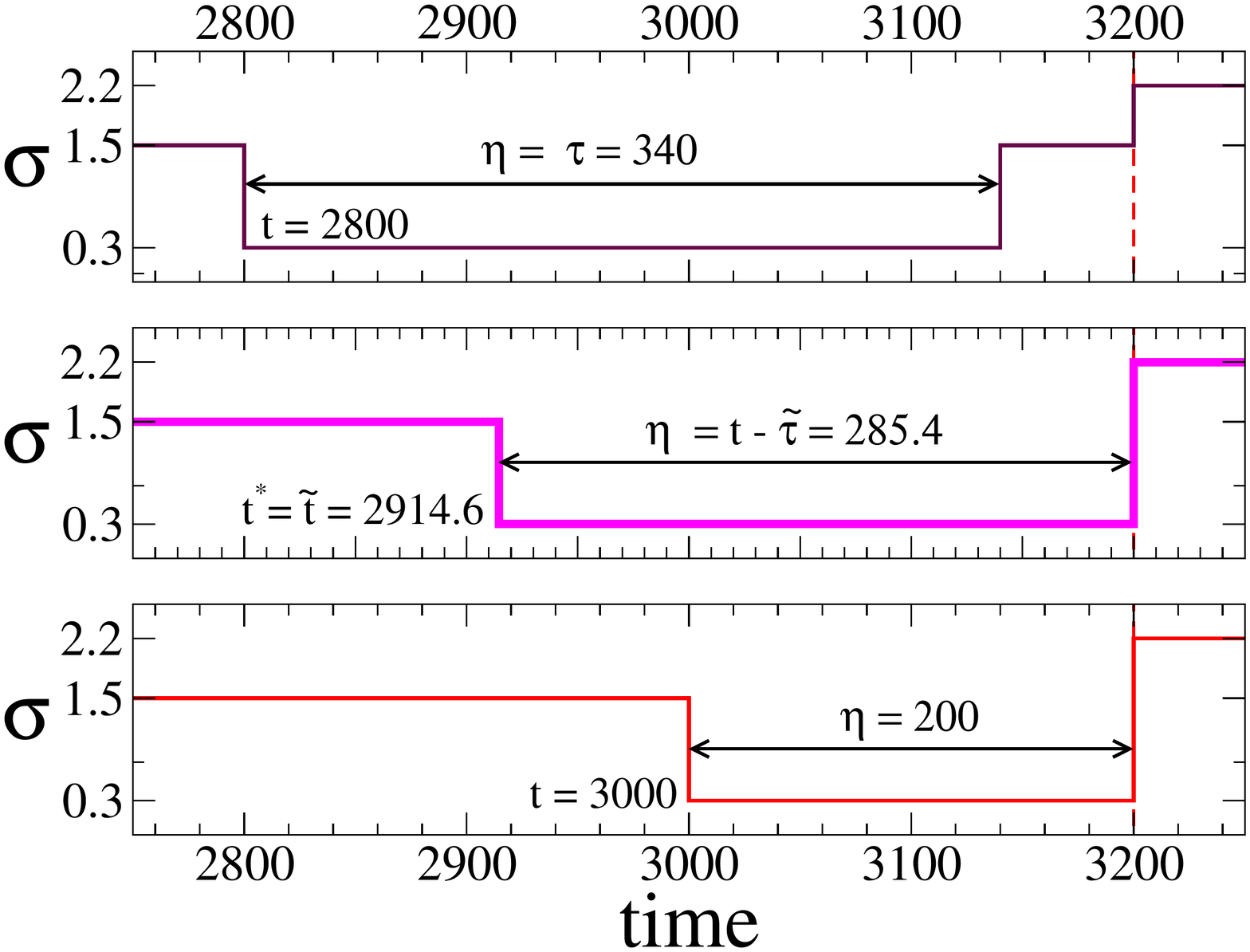}
    \includegraphics[width=8.3cm]{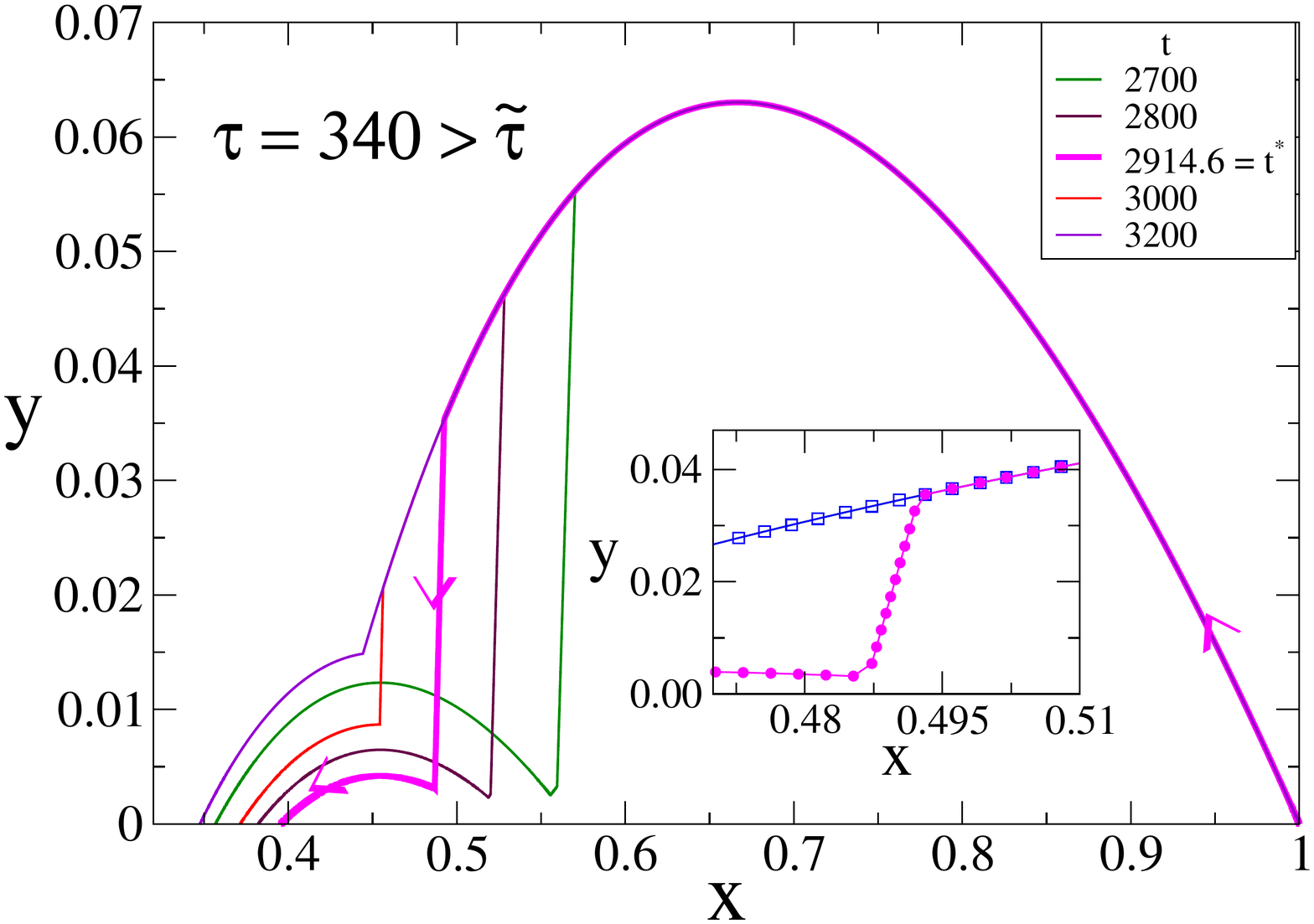}
    \caption{System's trajectory in the $x-y$ phase space (right
      panels), for $\gamma=0.01$, $\sigma_0=2.2$, $\sigma_1=0.3$,
      $\sigma_2=1.5$ and $T=3200$, and the three values of $\tau$ indicated
      in the legends corresponding to the different regimes of the
      optimum time $t^*$ (pink lines).  Left panels show the time
      evolution of $\sigma$ for three different initial times $t$ of
      the hard quarantine in each case. The optimum times are
      $t^*=\overline{t} \simeq 3103.5$ for $\tau=50$ (top panels),
      $t^*=T-\tau=3020$ for $\tau=180$ (middle panels) and $t^*=\tilde{t} \simeq
      2914.6$ for $\tau=340$ (bottom panels).}
    \label{x-y-3200-03}
  \end{center}
\end{figure}

In the right panels of Figure~\ref{x-y-3200-03} we show the system's evolution in the $x-y$ space for three different values of $\tau$ corresponding to the different behaviour of $t^*$.  Unlike the previously studied cases where $\sigma_2=\sigma_0$ (Figs.~\ref{x-y-2600-00} and \ref{x-y-2600-03}), here we observe that the curves $[x(t),y(t)]$ may exhibit up to three different regimes within the control period $T$, which is due to the fact that $\sigma$ jumps three times in that interval: from $\sigma_2$ to $\sigma_1$ at time $t$, from $\sigma_1$ to $\sigma_2$ at $t+\tau$ and from $\sigma_2$ to $\sigma_0$ at $T$.  This can be clearly seen in the $t^*=3103.5$ curve for $\tau=50<\overline{\tau}$ (inset of top-right panel of Fig.~\ref{x-y-3200-03}).  For $\tau$ in the other two regions ($\tau=180$ and $340$), the hard quarantine ends at $T$ for $t^*$, and thus $\sigma$ jumps twice and $[x(t),y(y)]$ exhibits two regimes in $[0,T]$ (insets of middle-right and bottom-right panels).  As in the previously studied cases, $y(t)$ drops to the lowest level curve in the interval $[t,t+\eta]$ for the optimum time $t^*$ (pink curves).

\section{Conclusions}
\label{se: conclusions}

In this paper, we have studied an optimal control problem on a SIR dynamics, with a control on the reproduction number $\sigma(t)$ and a limitation in the duration of the intervention $T$ and quarantine.  Based on the
Pontryagin's maximum principle, we have given first order necessary conditions with a cost that
takes into account both the maximization of the susceptible population in the long term and a
penalization of the lock-down. We also point out that we have employed a novel proof to establish
our analytical results. Moreover, some numerical examples have been provided to show the
validity of our theoretical results.

Given a fixed time of intervention $T$ where control strategies can be applied, and a strict quarantine period $\tau<T$ that represents
the maximum time lapse for the stronger intervention, we proved that the  optimal strategy is bang-bang when the second term of the cost functional is linear with respect to the control. More
precisely, the optimal solution consists of switching at most twice between a mild and a stronger quarantine, where the latter lasts at most a time period $\tau$.

Although some studies have supported the idea that a too soon or too late intervention may
not  minimize the total mortality, we found a broader scenario.  This is because the optimal solution takes the value
$\sigma=\sigma_1$ corresponding to the lock-down on an interval $[t^*,t^*+\eta]\subseteq [0,T]$,
with $t^*$ and $\eta\le \tau$ depending on the initial data $x_0, y_0, \gamma, \tau,
\sigma_0, \sigma_1$  and $T$. In fact, we showed that, in some cases, the optimal strategy
consists of taking $t^*= 0 $ or $t^*+\eta = T$ (see theorem \ref{te:
control_optimo_sigma1gral} items 1 and 3-4 respectively). 
However, for an initial condition that corresponds to a real-life case scenario in which the percentage of the population that is infected is small when non-pharmaceutical interventions start, we obtained that the optimal strategy consists on delaying the beginning of the lock-down (items 2-4 from theorem  \ref{te:
control_optimo_sigma1gral}).  
For the particular case $\tau \ll T$,  this optimum consists in applying a mild mitigation policy at the beginning of the intervention, followed by a strong suppression policy and then a mild mitigation again (mild--strong--mild strategy); while for $\tau \lesssim T$ the optimum corresponds to a mild--strong mitigation strategy.

We have also studied the possibility of implementing intermittent quarantines, and the
possibility of applying suppression measures first, followed by mitigation measures. In both
cases, if the total duration of measures is limited, we have shown that they are not
optimal in order to maximize the fraction of susceptible individuals at the end of the pandemic.

A major concern with respect to the current COVID-19 crisis is the possibility  of an overload
of available treatment resources. Since the hospitalized individuals are a fraction of the infected
population, a natural objective is to keep the number of infected individuals below some
threshold for all times. In a future work we intend to extend our analytic results including a
running state constraint that takes this restriction into account. Also, we aim to study the role
of an underlying network of contacts, and  changes in contact rates
 due to individual measures triggered by  fear of contagion.



\section{Supplement}

We begin by computing the derivatives of
$x_{t_1,\eta}(T)=\Psi_1(T,t_1+\eta,x_2,y_2,\sigma_2)$ and
$y_{t_1,\eta}(T)=\Psi_2(T,t_1+\eta,x_2,y_2,\sigma_2)$ with respect to $t_1$.

We recall two properties for the solutions of ordinary differential equations. First, the relation
between the derivative with respect to initial time and the derivatives with respect to initial
data give us the equation
\begin{align} \label{eq: deriv_sol_resp_tiempo_inicial}
\frac{\partial \Psi_j(s,t,x,y,\sigma)}{\partial t}&=\frac{\partial \Psi_j(s,t,x,y,\sigma)}{\partial x} \gamma \sigma x y - \frac{\partial \Psi_j(s,t,x,y,\sigma)}{\partial y} \gamma y(\sigma x -1)
\end{align}
for $\Psi$ defined at the begining of section \ref{se: caracterization_optimal}, with $s\ge t$, $\sigma \in \left\{ \sigma_1,\sigma_2\right\}$, initial data $(x,y)\in {\cal D}$ at time $t$ and $j=1,2$.

Second, the dependence of the solution $\Psi(s,t,x,y,\sigma)$ with respect to initial data
$x,y$ is given by the following known equations. For simplicity of notation, when there is no
risk of confusion, we will denote $\Psi(s)$ for $\Psi(s,t,x,y,\sigma)$,

\begin{equation}
\left( \begin{array}{ll}
\frac{\partial \Psi_1}{\partial x} & \frac{\partial \Psi_1}{\partial y} \\
\frac{\partial \Psi_2}{\partial x}&\frac{\partial \Psi_2}{\partial y}
\end{array}\right)' (s)=
\left( \begin{array}{ll}
-\gamma \sigma \Psi_2 (s) & -\gamma \sigma \Psi_1(s)\\
\gamma \sigma \Psi_2(s)&\gamma ( \sigma \Psi_1(s)-1)
\end{array}\right) .
\left( \begin{array}{ll}
\frac{\partial \Psi_1}{\partial x} & \frac{\partial \Psi_1}{\partial y} \\
\frac{\partial \Psi_2}{\partial x}&\frac{\partial \Psi_2}{\partial y}
\end{array}\right)(s) .
\end{equation}
with initial data
\begin{equation}
\left( \begin{array}{ll}
\frac{\partial \Psi_1}{\partial x} & \frac{\partial \Psi_1}{\partial y} \\
\frac{\partial \Psi_2}{\partial x}&\frac{\partial \Psi_2}{\partial y}
\end{array}\right)(t)=Id .
\end{equation}
Then, for $i=1,2$, we call
\begin{subequations}  \label{eq: def_u_v}
\begin{align}
u_i(s)=u(s,t_i,x_i,y_i,\sigma_i)=\frac{\partial \Psi_1}{\partial x_i}(s,t_i,x_i,y_i,\sigma_i) - \frac{\partial \Psi_1}{\partial y_i}(s,t_i,x_i,y_i,\sigma_i),  \label{eq: def_u} \\
v_i(s)=v(s,t_i,x_i,y_i,\sigma_i)=\frac{\partial \Psi_2}{\partial x_i}(s,t_i,x_i,y_i,\sigma_i) - \frac{\partial \Psi_2}{\partial y_i}(s,t_i,x_i,y_i,\sigma_i),  \label{eq: def_v}
\end{align}
\end{subequations}
for $s\in [t_i,t_{i+1}]$,
and we have the system of equations on $u_i$ and $v_i$
\begin{equation} \label{eq: derivparc_u_v}
\left( \begin{array}{l}
u_i'(s) \\
v_i'(s)
\end{array}\right)=
\left( \begin{array}{ll}
-\gamma \sigma_i \Psi_2(s) & -\gamma \sigma_i \Psi_1(s)\\
\gamma \sigma_i \Psi_2(s)&\gamma ( \sigma_i \Psi_1(s)-1)
\end{array}\right) .
\left( \begin{array}{l}
u_i(s) \\
v_i(s)
\end{array}\right) .
\end{equation}
with initial data
\begin{equation}
\left( \begin{array}{ll}
u_i(t_i)\\
v_i(t_i)
\end{array}\right)=
\left( \begin{array}{ll}
1\\
-1
\end{array}\right).
\end{equation}
Therefore, after some computations and using \eqref{eq: deriv_sol_resp_tiempo_inicial} and \eqref{eq: def_u_v} we obtain for $s\in (t_1,t_1+\eta)$
\begin{subequations} \label{eq: deriv_xey2tini_tini}
\begin{align}  \label{eq: deriv_x2tini}
\frac{d x_{t_1,\eta}}{dt_1}(s)&= -\gamma(\sigma_2- \sigma_1) x_1 y_1 u_1(s),   \\
 \label{eq: deriv_y2tini}
 \frac{d y_{t_1,\eta}}{dt_1}(s)&=-\gamma(\sigma_2- \sigma_1) x_1 y_1 v_1(s),
\end{align}
\end{subequations}
and for $s\in (t_1+\eta,T]$
\begin{subequations} \label{eq: deriv_xeytini_tini_3}
\begin{align}
\frac{d x_{t_1,\eta}}{dt_1}(s)&=\gamma (\sigma_2-\sigma_1) x_2 y_2 u_2(s)  \nonumber\\
&-\gamma (\sigma_2-\sigma_1) x_1 y_1\left(\frac{\partial \Psi_1(s,t_2,x_2,y_2,\sigma_2)}{\partial x_2}
u_1(t_2)  +\frac{\partial \Psi_1(s,t_2,x_2,y_2,\sigma_2)}{\partial y_2}v_1(t_2)\right),\label{eq: deriv_xeytini_tini_3a}\\
\frac{d y_{t_1,\eta}}{dt_1}(s)&=\gamma (\sigma_2-\sigma_1) x_2 y_2 v_2(s) \nonumber\\
&-\gamma (\sigma_2-\sigma_1) x_1 y_1\left(\frac{\partial \Psi_2(s,t_2,x_2,y_2,\sigma_2)}{\partial x_2}u_1(t_2) +\frac{\partial \Psi_2(s,t_2,x_2,y_2,\sigma_2)}{\partial y_2}v_1(t_2) \right).\label{eq: deriv_xeytini_tini_3b}
\end{align}
\end{subequations}

Moreover, using that for any $(x_i,y_i)\in {\cal D}$, $\Psi(s,t_i,x_i,y_i,\sigma_i)$, satisfies for $s\in [t_i,t_{i+1}]$
\begin{align}
\Psi_1(s,t_i,x_i,y_i,\sigma_i) e^{-\sigma_i (\Psi_1(s,t_i,x_i,y_i,\sigma_i)+\Psi_2(s,t_i,x_i,y_i,\sigma_i))} =x_i e^{-\sigma_i (x_i+y_i)},
\end{align}
we compute the derivatives with respect to $x_i$ and $y_i$ and using $\Psi_1(s)=\Psi_1(s,t_i,x_i,y_i,\sigma_i)$, we obtain for $s\in [t_i,t_{i+1}]$
\begin{align}
\frac{\partial \Psi_1}{\partial x_i}(s)-\sigma_i \Psi_1(s) \left( \frac{\partial \Psi_1}{\partial x_i}(s)+\frac{\partial \Psi_2}{\partial x_i}(s)\right)  &=(1-\sigma_i x_i) \frac{\Psi_1(s)}{x_i}, \label{eq: deriv_Psi_1}  \\
\frac{\partial \Psi_1}{\partial y_i}(s)-\sigma_i \Psi_1(s) \left( \frac{\partial \Psi_1}{\partial y_i}(s)+\frac{\partial \Psi_2}{\partial y_i}(s)\right) &=-\sigma_i \Psi_1(s). \label{eq: deriv_Psi_2}
\end{align}
Then, substracting the last two equations
\begin{align} \label{eq: expres_uts}
 u_i(s) -\sigma_i  \Psi_1(s,t_i,x_i,y_i,\sigma_i) (  u_i(s)+ v_i(s))=\frac{ \Psi_1(s,t_i,x_i,y_i,\sigma_i)}{x_i}
\end{align}
for $s\in[t_i,t_{i+1}]$ and therefore using \eqref{eq: expres_uts}, from \eqref{eq:
derivparc_u_v} we have that $u_i$ satisfies the ordinary differential equation
\begin{align*}
 u_i'(s)&=\gamma \left(\sigma_i \Psi_1(s,t_i,x_i,y_i,\sigma_i)-\sigma_i \Psi_2(s,t_i,x_i,y_i,\sigma_i) -1\right) u_i(s)+\gamma  \frac{\Psi_1(s,t_i,x_i,y_i,\sigma_i)}{x_i}, \\
u_i(t_i)&= 1.
\end{align*}
In the rest of this section, for simplicity of notation we will denote $x$ and $y$ for $x_{t_1,\eta}$ and $y_{t_1,\eta}$, defined in \eqref{eq: xy_t_eta},  respectively.

Thus,  for $s\in[t_i,t_{i+1}]$ when $i=1,2$, we obtain
\begin{align}
u_i(s)&=\frac{x(s)y(s)}{x_i y_i}+\gamma \frac{x(s)y(s)}{x_i} \int_{t_i}^{s} \frac{1}{y(r)}dr \nonumber\\
& = \frac{x(s)}{x_i }+\gamma \sigma_i \frac{x(s)y(s)}{x_i} \int_{t_i}^{s} \frac{x(r)}{y(r)}dr,  \label{eq: u_ti} \\
v_i(s)&=-\frac{x(s)}{x_i}+\frac{(1-\sigma_ix(s))y(s)}{x_i}\gamma  \int_{t_i}^{s} \frac{x(r)}{y(r)}dr, \label{eq: v_ti} \\
u_i(s)+v_i(s)&=\frac{y(s)}{x_i}\gamma  \int_{t_i}^{s} \frac{x(r)}{y(r)}dr. \label{eq: umasv_ti}
\end{align}
Also, from \eqref{eq: derivparc_u_v} and \eqref{eq: deriv_Psi_2}
we can prove for $\Psi(T)=\Psi(T,t_2,x_2,y_2,\sigma_2)$  that
\begin{subequations} \label{eq: wywtilde}
\begin{align}
\frac{\partial \Psi_1}{\partial x_2}(T)+\frac{\partial\Psi_2}{\partial x_2}(T) =\frac{y(T)}{y_2}+ (1-\sigma_2 x_2)\left(u_2(T)+v_2(T) \right), \label{eq: wT}\\
\frac{\partial \Psi_1}{\partial y_2}(T)+\frac{\partial\Psi_2}{\partial y_2}(T) =\frac{y(T)}{y_2}- \sigma_2 x_2{x_2} \left(u_2(T)+v_2(T) \right). \label{eq: wTtilde}
\end{align}
\end{subequations}

Analogously, 
\begin{subequations} \label{eq: deriv_xeytini_eta}
\begin{align}
\frac{d x_{t_1,\eta}}{d\eta}(T)&=\gamma (\sigma_2-\sigma_1) x_2 y_2 u_2(T), \\
\frac{d y_{t_1,\eta}}{d \eta}(T)&=\gamma (\sigma_2-\sigma_1) x_2 y_2 v_2(T) .
\end{align}
\end{subequations}
We can now compute the derivatives of $J(t_1,\eta)$ given by \eqref{eq: Jteta}. From \eqref{eq: deriv_xinfty}, 

\begin{align} \label{eq: derivJ_resp_t_inicial}
\frac{\partial J}{\partial t_1}  (t_1,\eta)&= \frac{d x_{\infty}(x_{t_1,\eta}(T),y_{t_1,\eta}(T),\sigma_0)}{d t_1} \nonumber   \\
&= \frac{ x_{\infty,t_1,\eta}}{1-\sigma_0 x_{\infty,t_1,\eta}}
\left(\frac{1-\sigma_0 x_{t_1,\eta}(T)}{x_{t_1,\eta}(T)} \frac{d x_{t_1,\eta}(T)}{d t_1}-\sigma_0 \frac{d y_{t_1,\eta}(T)}{d t_1}\right) \nonumber\\
&= \frac{ x_{\infty,t_1,\eta}}{(1-\sigma_0 x_{\infty,t_1,\eta}) x_{t_1,\eta}(T)}
\left( (1-\sigma_0x_{t_1,\eta}(T)) \frac{d x_{t_1,\eta}(T)}{d t_1}-\sigma_0 x_{t_1,\eta}(T) \frac{d y_{t_1,\eta}(T)}{d t_1}\right),
\end{align}
and, from \eqref{eq: deriv_xeytini_tini_3}
\begin{align*}
&(1-\sigma_0 x_{t_1,\eta}(T)) \frac{d x_{t_1,\eta}(T)}{d t_1}-\sigma_0 x_{t_1,\eta}(T)\frac{d y_{t_1,\eta}(T)}{d t_1}\\
&=\gamma (\sigma_2-\sigma_1)x_2y_2 \left( (1-\sigma_0 x_{t_1,\eta}(T)) u_2(T)-\sigma_0 x_{t_1,\eta}(T)v_2(T) \right) \\
&-  \gamma (\sigma_2-\sigma_1) x_1 y_1 \left(  u_1(t_2) \left((1-\sigma_0 x_{t_1,\eta}(T)) \frac{\partial \Psi_1(T,t_2,x_2,y_2,\sigma_2)}{\partial x_2}+\sigma_0 x_{t_1,\eta}(T)\frac{\partial \Psi_2(T,t_2,x_2,y_2,\sigma_2)}{\partial x_2} \right) \right. \\
& \left. +v_1(t_2) \left( (1-\sigma_0 x_{t_1,\eta}(T))\frac{\partial \Psi_1(T,t_2,x_2,y_2,\sigma_2)}{\partial y_2}  +\sigma_0 x_{t_1,\eta}(T) \frac{\partial \Psi_2(T,t_2,x_2,y_2,\sigma_2)}{\partial y_2}\right) \right),
\end{align*}
using \eqref{eq: deriv_Psi_1}, \eqref{eq: deriv_Psi_2}, \eqref{eq: expres_uts} and \eqref{eq: wywtilde} we obtain
\begin{align*}
&(1-\sigma_0 x_{t_1,\eta}(T)) \frac{d x_{t_1,\eta}(T)}{d t_1}-\sigma_0 x_{t_1,\eta}(T)\frac{d y_{t_1,\eta}(T)}{d t_1}\\
&=\gamma (\sigma_2-\sigma_1) x_{t_1,\eta}(T)\left[ (1+(\sigma_2-\sigma_0)(u_2(T)+v_2(T))) (y_2-y_1+
(\sigma_2-\sigma_1)x_1y_1 (u_1(t_2)+v_1(t_2)))\right. \\
&\left. - (\sigma_2-\sigma_0)x_1y_1\frac{y_{t_1,\eta}(T)}{y_2}(u_1(t_2)+v_1(t_2)) \right].
\end{align*}
Also, using that for $i=1,2$
\begin{align} \label{eq: cota_integral}
\frac{1}{y_{i}}-\frac{1}{y_{i+1}}=-\int_{t_i}^{t_{i+1}}  \left( \frac{1}{y}\right)'(r) dr=\gamma \int_{t_i}^{t_{i+1}}   \frac{\sigma_i x(r)-1}{y(r)} dr,
\end{align}
we have that
\begin{align} \label{eq: cota_integral_ij}
y_{i+1}-y_i+\gamma y_i y_{i+1} (\sigma_j-\sigma_i) \int_{t_i}^{t_{i+1}} \frac{x(r)}{y(r)} dr=\gamma y_i y_{i+1} \int_{t_i}^{t_{i+1}}   \frac{\sigma_j x(r)-1}{y(r)} dr.
\end{align}
Therefore, from \eqref{eq: umasv_ti} and \eqref{eq: cota_integral_ij}
\begin{align}
&\frac{(1-\sigma_0 x_{t_1,\eta}(T)) \dfrac{d x_{t_1,\eta}(T)}{d t_1}-\sigma_0 x_{t_1,\eta}(T)\dfrac{d y_{t_1,\eta}(T)}{d t_1}}{\gamma^2 (\sigma_2-\sigma_1) x_{t_1,\eta}(T) y_{t_1,\eta}(T)y_1} \nonumber \\\nonumber\\
&= \left(1-\gamma y_2 \int_{t_2}^{T} \frac{\sigma_0 x(r)-1}{y(r)} dr \right) \int_{t_1}^{t_2} \frac{\sigma_2 x(r)-1}{y(r)} dr -(\sigma_2-\sigma_0) \int_{t_1}^{t_2} \frac{x(r)}{y(r)} dr \nonumber \\
&= \int_{t_1}^{t_2} \frac{\sigma_0 x(r)-1}{y(r)} dr -\gamma y_2 \int_{t_2}^{T} \frac{\sigma_0 x(r)-1}{y(r)} dr  \int_{t_1}^{t_2} \frac{\sigma_2 x(r)-1}{y(r)} dr  \label{eq: final_para_derivJ_respt}.
\end{align}

Then, replacing in \eqref{eq: derivJ_resp_t_inicial} we obtain

\begin{align}
\frac{\partial J}{\partial t_1}  (t_1,\eta)&= \frac{ \gamma^2 (\sigma_2-\sigma_1) y_{t_1,\eta}(T) y_1 x_{\infty,t_1,\eta}}{1-\sigma_0 x_{\infty,t_1,\eta}} . \nonumber \\
& \left[\int_{t_1}^{t_2} \frac{\sigma_0 x(r)-1}{y(r)} dr -\gamma y_2 \int_{t_2}^{T} \frac{\sigma_0 x(r)-1}{y(r)} dr  \int_{t_1}^{t_2} \frac{\sigma_2 x(r)-1}{y(r)} dr \right].
\end{align}

Note that for $\sigma_2=\sigma_0$ we have from \eqref{eq: cota_integral} for $i=2$ that
\begin{align*}
&(1-\sigma_0 x_{t_1,\eta}(T)) \frac{d x_{t_1,\eta}(T)}{d t_1}-\sigma_0 x_{t_1,\eta}(T)\frac{d y_{t_1,\eta}(T)}{d t_1}\\
&= \gamma^2 (\sigma_2-\sigma_1) x_{t_1,\eta}(T) y_1y_2 \int_{t_1}^{t_2} \frac{\sigma_2 x(r)-1}{y(r)} dr,
\end{align*}
and then,
\begin{align}
\frac{\partial J}{\partial t_1}  (t_1,\eta)&= \frac{  x_{\infty,t_1,\eta}\gamma^2 (\sigma_2-\sigma_1)  y_1y_2}{1-\sigma_0 x_{\infty,t_1,\eta}}
  \int_{t_1}^{t_2} \frac{\sigma_2 x(r)-1}{y(r)} dr.
\end{align}
On the other hand,
\begin{align} \label{eq: derivJ_resp_eta}
\frac{\partial J}{\partial \eta}  (t_1,\eta)&= \frac{d x_{\infty}(x_{t_1,\eta}(T),y_{t_1,\eta}(T),\sigma_0)}{d \eta}  -\kappa (\sigma_2-\sigma_1)\nonumber   \\
&= \frac{ x_{\infty,t_1,\eta}}{(1-\sigma_0 x_{\infty,t_1,\eta})x_{t_1,\eta}(T)}
\left((1-\sigma_0 x_{t_1,\eta}(T)) \frac{d x_{t_1,\eta}(T)}{d \eta}-x_{t_1,\eta}(T)\sigma_0 \frac{d y_{t_1,\eta}(T)}{d \eta}\right)  -\kappa (\sigma_2-\sigma_1),\nonumber   \\
\end{align}
where $x_{\infty,t,\eta}=x_{\infty}(x_{t,\eta}(T),y_{t,\eta}(T),\sigma_0)$.

From \eqref{eq: def_u}, \eqref{eq: def_v} and \eqref{eq: deriv_xeytini_eta}, we obtain
\begin{align} \label{eq: derivJ_resp_eta2}
&\frac{1}{x_{t_1,\eta}(T)}\left((1-\sigma_0 x_{t_1,\eta}(T)) \frac{d x_{t_1,\eta}(T)}{d \eta}-\sigma_0 x_{t_1,\eta}(T) \frac{d y_{t_1,\eta}(T)}{d \eta}\right) \nonumber   \\
&= \gamma  y_2(\sigma_2-\sigma_1)\frac{x_2}{x_{t_1,\eta}(T)}\left((1-\sigma_0 x_{t_1,\eta}(T))u_2(T)-\sigma_0 x_{t_1,\eta}(T) v_2(T)\right)\nonumber   \\
& = \gamma  y_2(\sigma_2-\sigma_1)\frac{x_2}{x_{t_1,\eta}(T)}\left((1-\sigma_2 x_{t_1,\eta}(T))u_2(T)-\sigma_2 x_{t_1,\eta}(T) v_2(T)-(\sigma_0-\sigma_2)x_{t_1,\eta}(T)(u_2(T)+v_2(T))\right),\nonumber   \\
\end{align}
and using \eqref{eq: expres_uts}, \eqref{eq: umasv_ti} and \eqref{eq: cota_integral_ij}
\begin{align}
& = \gamma  y_2(\sigma_2-\sigma_1)\frac{x_2}{x_{t_1,\eta}(T)}\left(\frac{x_{t_1,\eta}(T)}{x_2}-(\sigma_0-\sigma_2)x_{t_1,\eta}(T)(u_2(T)+v_2(T))\right)\nonumber   \\
& = \gamma  y_2(\sigma_2-\sigma_1)\frac{x_2}{x_{t_1,\eta}(T)}\left(\frac{x_{t_1,\eta}(T)}{x_2}-(\sigma_0-\sigma_2)x_{t_1,\eta}(T)\frac{y_{t_1,\eta}(T)}{x_2}\gamma \int_{t_2}^T \frac{x(r)}{y(r)} dr\right)\nonumber   \\
& = \gamma  y_2(\sigma_2-\sigma_1)\left(1-(\sigma_0-\sigma_2)y_{t_1,\eta}(T)\gamma \int_{t_2}^T \frac{x(r)}{y(r)} dr\right)\nonumber   \\
& = \gamma  (\sigma_2-\sigma_1)\left(y_{t_1,\eta}(T)-\gamma y_2 y_{t_1,\eta}(T) \int_{t_2}^T \frac{\sigma_0 x(r)-1}{y(r)} dr\right).\nonumber   \\
\end{align}
Thus, we have
\begin{align} \label{eq: derivJ_resp_eta3}
\frac{\partial J}{\partial \eta}  (t_1,\eta)&
= \frac{ x_{\infty,t_1,\eta}}{1-\sigma_0 x_{\infty,t_1,\eta}}\gamma  (\sigma_2-\sigma_1)y_{t_1,\eta}(T)\left(1-\gamma y_2\int_{t_2}^T \frac{\sigma_0 x(r)-1}{y(r)} dr\right)-\kappa (\sigma_2-\sigma_1).\nonumber   \\
\end{align}

\bibliographystyle{plain}

\begin{thebibliography}{10}

\bibitem{anderson1992infectious}
Roy~M Anderson and Robert~M May.
\newblock {\em Infectious diseases of humans: dynamics and control}.
\newblock Oxford University Press, 1992.

\bibitem{dePinho2014}
M.~H.~A. {Biswas}, L.~T. {Paiva}, and M.d.~R. {de Pinho}.
\newblock {A SEIR model for control of infectious diseases with constraints}.
\newblock {\em {Math. Biosci. Eng.}}, 11(4):761--784, 2014.

\bibitem{brauer2012mathematical}
Fred Brauer, Carlos Castillo-Chavez, and Carlos Castillo-Chavez.
\newblock {\em Mathematical models in population biology and epidemiology},
  volume~2.
\newblock Springer, 2012.

\bibitem{buono2013slow}
Camila Buono, Federico Vazquez, Pablo~Alejandro Macri, and LA~Braunstein.
\newblock Slow epidemic extinction in populations with heterogeneous infection
  rates.
\newblock {\em Physical Review E}, 88(2):022813, 2013.

\bibitem{Clarke2013}
Francis {Clarke}.
\newblock {\em {Functional analysis, calculus of variations and optimal
  control}}, volume 264.
\newblock London: Springer, 2013.

\bibitem{ferguson2020report}
Neil Ferguson, Daniel Laydon, Gemma Nedjati~Gilani, Natsuko Imai, Kylie
  Ainslie, Marc Baguelin, Sangeeta Bhatia, Adhiratha Boonyasiri, ZULMA
  Cucunuba~Perez, Gina Cuomo-Dannenburg, et~al.
\newblock Report 9: Impact of non-pharmaceutical interventions (npis) to reduce
  covid19 mortality and healthcare demand.
\newblock {\em Imperial College London}, 2020.

\bibitem{ferrari2021coupling}
Carlo~Giambiagi Ferrari, Juan~Pablo Pinasco, and Nicolas Saintier.
\newblock Coupling epidemiological models with social dynamics.
\newblock {\em Bulletin of Mathematical Biology}, 83(7):1--23, 2021.

\bibitem{ferreyra2021sir}
Emanuel~Javier Ferreyra, Matthieu Jonckheere, and Juan~Pablo Pinasco.
\newblock {S}{I}{R} dynamics with vaccination in a large configuration model.
\newblock {\em arXiv preprint arXiv:1912.12350}, 2019.

\bibitem{godara2021control}
Prakhar Godara, Stephan Herminghaus, and Knut~M Heidemann.
\newblock A control theory approach to optimal pandemic mitigation.
\newblock {\em PloS one}, 16(2):e0247445, 2021.

\bibitem{greenhalgh1988some}
David Greenhalgh.
\newblock Some results on optimal control applied to epidemics.
\newblock {\em Mathematical Biosciences}, 88(2):125--158, 1988.

\bibitem{zbMATH00784953}
Richard~F. {Hartl}, Suresh~P. {Sethi}, and Raymond~G. {Vickson}.
\newblock {A survey of the maximum principles for optimal control problems with
  state constraints}.
\newblock {\em {SIAM Rev.}}, 37(2):181--218, 1995.

\bibitem{zbMATH01579035}
Herbert~W. {Hethcote}.
\newblock {The mathematics of infectious diseases}.
\newblock {\em {SIAM Rev.}}, 42(4):599--653, 2000.

\bibitem{janson2014}
Svante Janson, Malwina Luczak, and Peter Windridge.
\newblock Law of large numbers for the {S}{I}{R} epidemic on a random graph
  with given degrees.
\newblock {\em Random Structures \& Algorithms}, 45(4):726--763, 2014.

\bibitem{kermack1927contribution}
William~Ogilvy Kermack and Anderson~G McKendrick.
\newblock A contribution to the mathematical theory of epidemics.
\newblock {\em Proceedings of the royal society of london. Series A, Containing
  papers of a mathematical and physical character}, 115(772):700--721, 1927.

\bibitem{ketcheson2020optimal}
David~I Ketcheson.
\newblock Optimal control of an sir epidemic through finite-time
  non-pharmaceutical intervention.
\newblock {\em J. Math. Biol.}, 83(7), 2021.

\bibitem{kohler2020robust}
Johannes K{\"o}hler, Lukas Schwenkel, Anne Koch, Julian Berberich, Patricia
  Pauli, and Frank Allg{\"o}wer.
\newblock Robust and optimal predictive control of the {COVID}-19 outbreak.
\newblock {\em Annual Reviews in Control}, 2020.

\bibitem{lagorio2011quarantine}
C~Lagorio, Mark Dickison, F~Vazquez, Lidia~A Braunstein, Pablo~A Macri,
  MV~Migueles, Shlomo Havlin, and H~Eugene Stanley.
\newblock Quarantine-generated phase transition in epidemic spreading.
\newblock {\em Physical Review E}, 83(2):026102, 2011.

\bibitem{Ledzewicz2011}
Urszula {Ledzewicz} and Heinz {Sch\"attler}.
\newblock {On optimal singular controls for a general SIR-model with
  vaccination and treatment}.
\newblock {\em {Discrete Contin. Dyn. Syst.}}, 2011:981--990, 2011.

\bibitem{zbMATH05164170}
Suzanne {Lenhart} and John~T. {Workman}.
\newblock {\em {Optimal control applied to biological models}}.
\newblock Boca Raton, FL: Chapman \& Hall/CRC, 2007.

\bibitem{morris2020optimal}
Dylan~H Morris, Fernando~W Rossine, Joshua~B Plotkin, and Simon~A Levin.
\newblock Optimal, near-optimal, and robust epidemic control.
\newblock {\em arXiv preprint arXiv:2004.02209}, 2020.

\bibitem{palmer2020optimal}
Aaron~Z Palmer, Zelda~B Zabinsky, and Shan Liu.
\newblock Optimal control of {COVID}-19 infection rate with social costs.
\newblock {\em arXiv preprint arXiv:2007.13811}, 2020.

\bibitem{pueyo2020coronavirus}
Tomas Pueyo.
\newblock Coronavirus: the hammer and the dance.
\newblock {\em Medium [Internet]}, 2020.

\bibitem{tsay2020modeling}
Calvin Tsay, Fernando Lejarza, Mark~A Stadtherr, and Michael Baldea.
\newblock Modeling, state estimation, and optimal control for the {US}
  {COVID}-19 outbreak.
\newblock {\em Scientific reports}, 10(1):1--12, 2020.

\bibitem{velasquez2017interacting}
F{\'a}tima Vel{\'a}squez-Rojas and Federico Vazquez.
\newblock Interacting opinion and disease dynamics in multiplex networks:
  discontinuous phase transition and nonmonotonic consensus times.
\newblock {\em Physical Review E}, 95(5):052315, 2017.

\bibitem{volz2008sir}
Erik Volz.
\newblock {S}{I}{R} dynamics in random networks with heterogeneous
  connectivity.
\newblock {\em Journal of Mathematical Biology}, 56(3):293--310, 2008.

\end{thebibliography}

\end{document}